\documentclass[12pt]{amsart}
\usepackage{amssymb}
\usepackage[dvips]{graphics}
\ifx\mathrm\undefined\let\mathrm\rm\fi
\ifx\mathbf\undefined\let\mathbf\bf\fi
\ifx\mathfrak\undefined\let\mathfrak\frak\fi
\ifx\mathcal\undefined\let\mathcal\cal\fi
\ifx\mathbb\undefined\let\mathbb\Bbb\fi
\ifx\emph\undefined\let\emph\it\fi
\newcommand{\g}{{{\mathfrak g}\,}}

\newcommand{\h}{{{\mathfrak h\,}}}
\newcommand{\Id}{{\operatorname{Id}}}
\newcommand{\Z}{{\mathbb Z}}

\newcommand{\R}{{\mathbb R}}
\newcommand{\C}{{\mathbb C}}

\newcommand{\Ref}[1]{{(\ref{#1})}}
\newcommand{\be}{\begin{displaymath}}
\newcommand{\ee}{\end{displaymath}}
\newcommand{\bea}{\begin{eqnarray*}}
\newcommand{\eea}{\end{eqnarray*}}

\newcommand{\Mu}{{\mathrm{M}}}

\newenvironment{example}
{\noindent{\bf Example\/}:}{\par}

\let\dl\delta \let\Dl\Delta
 \let\eps\varepsilon \let\epsilon\eps

\let\al\alpha

\let\la\lambda \let\La\Lambda
\let\om\omega \let\Om\Omega
 \let\phi\varphi

 \let\Si\Sigma

\newcommand{\half}{\frac12}
\newcommand{\bean}{\begin{eqnarray}}
\newcommand{\eean}{\end{eqnarray}}

\newcommand{\res}{{\operatorname{res}}}

\newcommand{\vs}{\vspace{.5\baselineskip}}
\newtheorem%
{thm}{Theorem}
\newtheorem%
{proposition}[thm]{Proposition}
\newtheorem%
{lemma}[thm]{Lemma}
\newtheorem%
{lemmadef}[thm]{Lemma-Definition}
\newtheorem%
{corollary}[thm]{Corollary}
\newtheorem%
{conjecture}[thm]{Conjecture}
\newcommand{\End}{{\operatorname{End\,}}}

\newcommand{\Sym}{{\operatorname{Sym\,}}}

\renewcommand{\Im}{{\operatorname{Im\,}}}

\title[Resonance Relations]
{Resonance Relations for Solutions of the Elliptic QKZB Equations,
Fusion Rules, and
Eigenvectors of Transfer Matrices of Restricted Interaction-round-a-face Models}

\author[G. Felder and A. Varchenko]
{G.\ Felder$^{\,\star}$ and
A.\ Varchenko$^{\,\diamond,1}$}
\thanks{$^1$ Supported in part by NSF grant DMS-9801582}

\begin{document}
\maketitle

\begin{center}
{\it
$^\star$Department of Mathematics, ETH Zentrum, CH-8092 Z\"urich, Switzerland,

felder@math.ethz.ch

\medskip

$^\diamond$Department of Mathematics, University of North Carolina
at Chapel Hill, 

Chapel Hill, NC 27599-3250, USA,

av@math.unc.edu}
\end{center}



\centerline{January, 1999}

\begin{abstract}
 Conformal blocks for the WZW model on tori can be represented by
 vector valued Weyl anti-symmetric theta functions on the Cartan
 subalgebra satisfying vanishing conditions on root hyperplanes. We
 introduce a quantum version of these vanishing conditions in the $sl_2$
 case.  They are compatible with the qKZB equations and are obeyed by
 the hypergeometric solutions as well as by their critical level
 counterpart, which are Bethe eigenfunctions of IRF row-to-row transfer
 matrices. In the language of IRF models the vanishing conditions turn
 out to be equivalent to the $sl_2$ fusion rules defining restricted
 models.
\end{abstract}


\thispagestyle{empty}

\section{Introduction}

This paper has two motivations: to understand quantization of vanishing conditions
for conformal blocks on tori and to construct eigenvectors of row-to-row transfer
matrices of restricted interaction-round-a-face models in statistical mechanics.

The Knizhnik-Zamolodchikov-Bernard equations are a system of differential equations
arising in conformal field theory on Riemann surfaces. 
For each $g, n\in \Z_{\geq 0}$, a simple complex Lie algebra $\g$, $n$ highest weight
$\g$-modules $V_i$ and a complex parameter $\kappa$, we have such a system of equations.
In the case of genus $g=1$, they have the form
\bean\label{KZB}
\kappa \partial_{z_j}v=-\sum_\nu h^{(j)}_\nu\partial_{\la_\nu}v
+\sum_{l,\,l\neq j}r(z_j-z_l,\la,\tau)^{(j,l)}v,
\eean
\bean\label{heaT}
4\pi i \kappa \partial_\tau v=\Delta_\la v+{1\over 2}\sum_{i,j}s(z,\la,\tau)^{(i,j)}v.
\eean
The unknown function $v$ takes values in the zero weight space $V[0]=
\cap_{x\in \h}\text{Ker}(x)$ of the tensor product $V=V_1\otimes \cdots \otimes V_n$
with respect to the Cartan subalgebra $\h$ of $\g$. It depends on variables $z_1,...,z_n\in \C$,
modulus $\tau$ of the elliptic curve and $\la=\sum \la_\nu h_\nu \in \h$, where $(h_\nu)$ is an orthonormal basis
of $\h$ with respect to a fixed invariant bilinear form. The notation $x^{(i)}$
for $x\in \End (V_i)$ or $x\in \g$  means $1\otimes \cdots \otimes x \otimes \cdots \otimes 1$.
Similarly $x^{(i,j)}$ denotes the action on the $i$-th and $j$-th factors of $x\in \End (V_i\otimes V_j)$.
In the equation, $r,s\in \g\otimes \g$ are suitable given tensor valued functions.

The KZB equations can be viewed as equations for horizontal sections of a connection with fiber
over $z_1,...,z_n,\tau$ being the space of functions of $\la$. 
If $\kappa$ is an integer not less than the dual Coxeter number $h^\vee$ of $\g$, 
$V_i$ are irreducible representations with highest weights $\La_i,\,
(\theta,\La_i)\leq \kappa - h^\vee$ with $\theta$ being the highest root, then
the connection has an invariant finite dimensional
subbundle of conformal blocks coming from conformal field theory. 
The fiber over $(z_1,...,z_n,\tau)$ of the subbundle of conformal blocks consists of
holomorphic functions of $\la$ obeying the following three conditions \cite{FG}, \cite{FW},
\cite{EFK}.
\begin{enumerate}
\item[I.]
The functions $v(\la)$ are $V[0]$-valued theta functions.
Namely, $v(\la)$ are periodic  with respect to the coroot lattice
$Q^\vee$ and
\bean\label{theta-property}
v(\la+q\tau)=e^{-\pi i \kappa(q,q)-2\pi i \kappa(q,\la)+2\pi i \sum_j
q^{(j)}z_j}v(\la), \qquad \forall q\in Q^\vee,
\notag
\eean
\item[II.] The functions $v(\la)$ are anti-symmetric with respect to the standard action of the Weyl group
of $\g$ on $V[0]$-valued functions of $\la \in \h$,
\item[III.] The functions $v(\la)$ satisfy the vanishing conditions.
Namely, for all roots $\al$, integers $r,s,l,\,l\geq 0$, and $\xi\in
\g_\al$, where $\g_\al$ is the root space corresponding to $\al$, we have
\bean
 (\,\sum_{j=1}^n e^{2\pi i s z_j}\xi^{(j)}\,)^lv\,=\,O((\al(\la)-r-s\tau)^{l+1}),
\notag
\eean
as $\al(\la)\to r+s\tau$.
\end{enumerate}

In \cite{F}, an elliptic quantum group $E_{\tau,\eta}(\g)$ was introduced and a difference version
of the KZB equations \Ref{KZB} was defined. Later a difference version of the KZB heat equation \Ref{heaT}
was suggested in \cite{FV3}. For the case $\g=sl_2$ considered in this paper, 
the quantum Knizhnik-Zamolodchikov-Bernard equations are linear difference equations
with step $p\in\C$ on a $V[0]$-valued function
$v(z_1,...,z_n,\la,\tau)$ of the form
\bean\label{quantum-kzb}
{}
\\
v(z_1,...,z_j+p,z_n,\la,\tau)\,=\,K_j(z_1,...,z_n,\tau,p,\eta)\,v(z_1,...,z_n,\la,\tau),
\qquad j=1,...,n,
\notag
\eean
where the linear operators $K_j$ are defined in terms of R-matrices of the elliptic quantum group
$E_{\tau,\eta}(sl_2)$. More precisely, the space $V$ is endowed with
an $E_{\tau,\eta}(sl_2)$-module structure depending on $(z_1,...,z_n)$
and denoted $V_1(z_1)\otimes ... \otimes V_n(z_n)$,  
and the operators are defined in terms of this module structure.
The difference version on the KZB heat equation
\Ref{heaT} has the form
\bean\label{q-heat}
v(z_1,...,z_n,\tau)\,=\,T(z_1,...,z_n,\tau,p,\eta)\,v(z_1,...,z_n,\la,\tau+p)
\eean
for a suitable linear operator $T$, but we do not consider the qKZB heat equation in this paper.

In the semiclassical limit $p,\eta \to 0$ with $p/\eta=-2\kappa$  fixed, the qKZB equations and
the qKZB heat equation turn into the corresponding differential equations.

In \cite{FTV2}, \cite{MV}, hypergeometric solutions of the qKZB equations \Ref{quantum-kzb}
we constructed. In \cite{FV4}, under certain conditions, it was shown that 
all solutions of the qKZB equations are sums of the  hypergeometric solutions. Under the same conditions
it was shown that the hypergeometric solutions also satisfy the qKZB heat equation \Ref{q-heat}.

The qKZB equations  \Ref{quantum-kzb} 
can be viewed as equations for horizontal sections of a discrete
connection with fiber over $(z_1,...,z_n)$ being the space of functions of $\la$

The problem is to find difference analogs of Conditions I-III for conformal blocks in conformal field theory.
The difference analog of Condition I is indicated in Section 6.5 of \cite{FTV2},
it is a theta function property analogous to I.
It is not difficult to describe a difference analog of Condition II. In Section \ref{weyl-refl}
we consider the case when all $V_i$ are finite dimensional
irreducible $sl_2$ modules. In this case we introduce an action
of the Weyl group on $V[0]$ and show that the qKZB operators commute with this action,
thus allowing us to consider the Weyl anti-symmetric  solutions.

Our first main motivation was to find an analog of vanishing conditions III for
Weyl anti-symmetric solutions of the qKZB equations. It turns out that each coordinate of a
Weyl anti-symmetric hypergeometric solution of the qKZB equations is equal to zero
at  a special set of values of $\la$ depending on the $\h$-weight of the coordinate. 
This set of special values of $\la$ is described in terms of fusion rules of $sl_2$,
see Section \ref{vanishing-cond}.
We show in Section \ref{remarks}
in the simplest nontrivial case that the semiclassical limit of the vanishing
conditions for solutions of the qKZB equations is the vanishing condition III
for conformal blocks.

Now we shall formulate the vanishing conditions in a special case.
For $i=1,...n$, let $V_i$ be the irreducible two dimensional representations of $sl_2$
with the standard weight basis $e[1],e[-1]$. For an even $n$, the zero weight
space $V[0]$ is nontrivial. Any $V[0]$-valued function has the form $u=\sum \, u_{m_1,...,m_n} 
\,e[m_1] \otimes ... \otimes e[m_n]$ where $ m_1+...+m_n=0$ and $u_{m_1,...,m_n}$ are scalar functions.

\begin{thm}\label{thm-1}
Let $u(z_1,...,z_n,\la)$ be a Weyl anti-symmetric $V[0]$-valued hypergeometric
solution of the qKZB equations
and  $u_{m_1,...,m_n}(z_1,...,z_n,\la)$ 
 one of its coordinates. Form a sequence of integers $\Sigma^1,...,\Sigma^n$
where $\Sigma^j=m_1+...+m_j$. Then for a non-negative integer $k$, we have
$u_{m_1,...,m_n}(z_1,...,z_n, 2\eta k)=0$, if at least one of the integers
$-\Sigma^1+k,...,-\Sigma^n+k$ is not positive, and
$u_{m_1,...,m_n}(z_1,...,z_n, -2\eta k)=0$, if at least one of the integers
$\Sigma^1+k,...,\Sigma^n+k$ is not positive.
\end{thm}

Fix $(z_1,...,z_n)\in \C^n$.
Consider the $E_{\tau,\eta}(sl_2)$-module $V_1(z_1)\otimes ... \otimes V_n(z_n)$,  
where all $V_i$ are two dimensional irreducible $sl_2$ modules.
For any complex number $w$, one defines a second order difference operator
$$
T(w)f(\la)\,=\,A(w,\la)f(\la-2\eta)+D(w,\la)f(\la+2\eta)
$$
acting on the space of $V[0]$-valued meromorphic functions of $\la$
and called the transfer matrix, see Section \ref{t-mat}. Here
$A, D : V[0]\to V[0]$ are linear operators meromorphically depending on $w,\la$.
The transfer matrices commute for different values of $w$.
In \cite{FV2}, we introduced an algebraic Bethe ansatz and constructed 
common eigenfunctions of the transfer matrices. The Bethe eigenfunctions are
holomorphic functions of $\la$. The transfer matrices commute with the action
of the Weyl group on $V[0]$. Thus anti-symmetrization of a Bethe eigenfunction
produces a holomorphic Weyl anti-symmetric common eigenfunction of the
transfer matrices. In Section \ref{t-mat} we show that the Weyl anti-symmetric
Bethe eigenfunctions obey the following vanishing conditions.

\begin{thm}\label{thm-2}
Let $u(\la)$ be a Weyl anti-symmetric Bethe eigenfunction
and  $u_{m_1,...,m_n}(\la)$ 
 one of its coordinates. Form a sequence of integers $\tilde \Sigma^1,...,\tilde\Sigma^n$
where $\tilde \Sigma^j=m_j+...+m_n$. Then for a non-negative integer $k$, we have
$u_{m_1,...,m_n}(2\eta k)=0$, if at least one of the integers
$-\tilde \Sigma^1+k,...,-\tilde \Sigma^n+k$ is not positive, and
$u_{m_1,...,m_n}(z_1,...,z_n, -2\eta k)=0$, if at least one of the integers
$\tilde \Sigma^1+k,...,\tilde \Sigma^n+k$ is not positive.
\end{thm}

For $\mu\in \C$, let $C_\mu=\{2\eta(\mu+ j)\,|\,j\in\Z\}$ and let
${\mathcal F}_\mu(V[0])$ be the space of $V[0]$-valued functions on $C_\mu$. 
For generic $\mu$, the transfer matrices act on ${\mathcal F}_\mu(V[0])$
and are well defined. The restriction of Bethe eigenfunctions to $C_\mu$
give common eigenfunctions. The case $\mu=0$ and $2\eta=1/N$ for a natural number $N$
is of special interest. Although the transfer matrices are not well defined
on whole ${\mathcal F}_{\mu=0}(V[0])$, if one considers the finite dimensional subspace
 of $V[0]$-valued function on $C_{\mu=0}^{\{N\}}=
\{2\eta j\,|\,j=1,...,N-1\}$, then the the restriction of the transfer matrices to this subspace
is well defined and the transfer matrices commute for different values of $w$.
The resulting commuting linear operators on this finite dimensional space are called the
row-to-row transfer matrices of the restricted interaction-round-a-face model in statistical mechanics
considered in \cite{ABF}.

Our second main motivation was to construct common eigenvectors of the row-to-row transfer
matrices. In Section \ref{2nd-mot} we show that the restriction of a Weyl anti-symmetric
Bethe eigenfunction to $C_{\mu=0}^{\{N\}}$ is a common eigenvector of the row-to-row transfer
matrices. The vanishing conditions of Theorem \ref{thm-2} play the decisive role in the proof.

Thus, our two main results, the quantization of vanishing conditions for conformal block on tori
and the construction of eigenvectors for restricted models, display another connection between
conformal field theory and statistical mechanics.

The paper is organized as follows. We begin by introducing the notion of R-matrices and
the qKZB equations in Section \ref{sqkzb}.
In Section \ref{fs} we give a geometric construction of R-matrices as transition matrices
between special bases in suitable spaces of functions. The elements of those bases are called
the weight functions. The weight functions are central objects of this paper. In Section
\ref{fs} we prove resonance relations for weight functions, the relations which form the ground for
all variants of vanishing conditions in this paper. 

In Section \ref{section-4} we study
poles of the R-matrices and relations for their matrix coefficients.
In Section \ref{hypergsol} we describe the hypergeometric solutions of the qKZB equations.
Section \ref{Verma} is devoted to resonance relations for solutions of the qKZB equations
with values in a tensor product of Verma modules. 
In Section \ref{Regularity} we show that
the resonance relations are necessary and sufficient conditions for a holomorphic function
of $\la$ remain holomorphic under the action of the qKZB operators and the R-matrices 
permuting tensor factors. 

In Section \ref{remarks} we consider the semiclassical limit of 
the resonance relations in the simplest nontrivial case. 

In Section \ref{weyl-refl}
we introduce an action of the Weyl group and show that the action commutes with the
qKZB operators. 
In Section \ref{vanishing-cond} we prove the vanishing conditions
for Weyl anti-symmetric hypergeometric
solutions of the qKZB equations with values in a tensor product of finite dimensional irreducible
representations and formulate the result in terms of $sl_2$ fusion rules.
Notice that the KZB differential equations \Ref{KZB}, \Ref{heaT} 
have singularities at the hyperplanes $\alpha(\lambda)=r+s \tau$, $r,s\in \Z$.
In \cite{FV5}, it was shown that 
meromorphic solutions of the KZB equations obeying Condition II obey also 
Condition III. In the same spirit  the vanishing conditions are,
for Weyl anti-symmetric functions,
consequences of resonance relations, which are necessary and sufficient
conditions for a holomorphic function of $\lambda$ to remain holomorphic
after acting with the qKZB operators $K_i$ and with the R-matrices permuting
tensor factors.

In Section \ref{sba} we study vanishing conditions for Bethe ansatz eigenfunctions,
the results are formulated in terms of $sl_2$ fusion rules and for special values of
parameters in terms of $U_q(sl_2)$ fusion rules where $q$ is a root of unity.
In Section \ref{RIRF} we construct eigenvectors of restricted interaction-round-a-face models.



\section{$R$-matrices, qKZB equations }\label{sqkzb}

\subsection{$R$-matrices}
The qKZB equations are given in terms of $R$-matrices of elliptic quantum
groups. In the $sl_2$ case, these $R$-matrices have the following properties.
Let $\h=\C h$ be a one-dimensional Lie algebra with generator $h$. For each
$\Lambda\in\C$ consider the $\h$-module $V_\Lambda=\oplus_{j=0}^\infty\C e_j$,
with $he_j=(\Lambda-2j)e_j$. For each pair $\Lambda_1$, $\Lambda_2$ of complex
numbers we have a meromorphic function, called the $R$-matrix,
$R_{\Lambda_1,\Lambda_2}(z,\lambda)$ of two complex variables, with values in
$\End(V_{\Lambda_1}\otimes V_{\Lambda_2})$.

The main properties of the $R$-matrices are
\begin{enumerate}
\item[I.] The zero weight property: for any $\Lambda_i$, $z,\lambda$,
$[R_{\Lambda_1,\Lambda_2}(z,\lambda),h^{(1)}+h^{(2)}]=0$.
\item[II.] For any $\Lambda_1,\Lambda_2,\Lambda_3$, the dynamical Yang--Baxter
equation
\bea
R_{\Lambda_1,\Lambda_2}(z,\lambda-2\eta h^{(3)})^{(12)}
R_{\Lambda_1,\Lambda_3}(z+w,\lambda)^{(13)}
R_{\Lambda_2,\Lambda_3}(w,\lambda-2\eta h^{(1)})^{(23)}
\\
{}=
R_{\Lambda_2,\Lambda_3}(w,\lambda)^{(23)}
R_{\Lambda_1,\Lambda_3}(z+w,\lambda-2\eta h^{(2)})^{(13)}
R_{\Lambda_1,\Lambda_2}(z,\lambda)^{(12)},
\eea
holds in $\End(V_{\Lambda_1}\otimes V_{\Lambda_2}\otimes V_{\Lambda_3})$
for all $z,w,\lambda$.
\item[III.] For all $\Lambda_1$, $\Lambda_2$, $z,\lambda$,
$R_{\Lambda_1,\Lambda_2}(z,\lambda)^{(12)}
R_{\Lambda_2,\Lambda_1}(-z,\lambda)^{(21)}=\Id$.
This property is called the ``unitarity''.
\end{enumerate}

We use the following notation: if $X\in\End(V_i)$, then we denote by
$X^{(i)}\in\End(V_1\otimes\dots\otimes V_n)$ the operator
$\cdots\otimes\Id\otimes X\otimes\Id\otimes\cdots$, acting non-trivially on
the $i$th factor of a tensor product of vector spaces, and if
$X=\sum X_k\otimes Y_k\in\End(V_i\otimes V_j)$, then we set
$X^{(ij)}=\sum X_k^{(i)}Y_k^{(j)}$. If $X(\mu_1,\dots,\mu_n)$ is a function
with values in $\End(V_1\otimes\dots\otimes V_n)$, then
$X(h^{(1)},\dots,h^{(n)})v=X(\mu_1,\dots,\mu_n)v$ if $h^{(i)}v=\mu_iv$,
for all $i=1,\dots,n$.

For each $\tau$ in the upper half plane and generic $\eta\in\C$ (``Planck's
constant'') a system of R-matrices $R_{\Lambda_1,\Lambda_2}(z,\lambda)$
obeying {I\,--\,III} was constructed in \cite{FV1}. They are characterized by
an intertwining property with respect to the action of the elliptic quantum
group $E_{\tau,\eta}(sl_2)$ on tensor products of evaluation Verma modules.
We recall a geometric construction of the R-matrices in Sec. \ref{fs}.

\subsection{The qKZB equations}\label{qkzb}
Let $n>1$ be a natural number.
Fix the parameters $\tau,\eta$. Fix also $n$ complex numbers
$\Lambda_1,\dots,\Lambda_n$ and an additional parameter $p\in\C$.
Let $V_{\vec\La}=V_{\Lambda_1}\otimes\cdots\otimes V_{\Lambda_n}$. The kernel of
$h^{(1)}+\dots+h^{(n)}$ on $V_{\vec\La}$ is called the zero-weight space and is denoted
$V_{\vec\La}[0]$. 

 Let $K_i(z,\tau,p)$ be the qKZB operators acting on
the space ${\mathcal F} (V_{\vec\La}[0])$ of  meromorphic functions
of $\lambda\in\C$ with values in the zero weight space $V_{\vec\Lambda}[0]$.
They have the form:
\begin{eqnarray*}
K_j(z,\tau,p)&=&R_{j,j-1}(z_j-z_{j-1}+p,\tau)\cdots R_{j,1}(z_j-z_1+p,\tau)
\\ & &\Gamma_jR_{j,n}(z_j-z_n,\tau)\cdots R_{j,j+1}(z_j-z_{j+1},\tau).
\end{eqnarray*}
The operators $R_{j,k}(z,\tau)$ are defined by the formula
\[
R_{j,k}(z,\tau)\,v(\lambda)=R_{\Lambda_j,\Lambda_k}(z,\lambda-2\eta
\sum_{{l=1,l\neq j}}^{k-1}h^{(l)},\tau)\,v(\lambda),
\]
and $(\Gamma_jv)(\lambda)=v(\lambda-2\eta\mu)$ if $h^{(j)}v(\lambda)=
\mu v(\lambda)$.

Let $\delta_j$, $j=1,\dots,n$ be the standard basis of $\C^n$.
The qKZB system of difference equations
\begin{equation}\label{qk1}
v(z+p\delta_j)=K_j(z,\tau,p)\,v(z),\qquad j=1,\dots n,
\end{equation}
for a function $v(z)$ on $\C^n$ with values in ${\mathcal F} (V_{\vec\Lambda}[0])$ is
compatible, i.e., we have
\begin{equation}\label{compatib}
K_j(z+p\delta_l,\tau,p)K_l(z,\tau,p)=
K_l(z+p\delta_j,\tau,p)K_j(z,\tau,p),
\end{equation}
for all $j,l$, as a consequence of the dynamical Yang--Baxter equations
satisfied by the $R$-matrices.

 We also consider the {\it mirror} qKZB operators
\begin{eqnarray*}
K_j^\vee(z,\tau,p)
&=&R^\vee_{j,j+1}(z_j-z_{j+1}+p,\tau)\cdots R^\vee_{j,n}(z_j-z_n+p,\tau)
\\ & &\Gamma_j
R^\vee_{j,1}(z_j-z_1,\tau)\cdots R^\vee_{j,j-1}(z_j-z_{j-1},\tau),
\end{eqnarray*}
with
\[
R^\vee_{j,k}(z,\tau)\,v(\lambda)=R_{\Lambda_j,\Lambda_k}(z,\lambda-2\eta
\sum_{{l=k+1,l\neq j}}^{n}h^{(l)},\tau)\,v(\lambda),
\]
The corresponding system of qKZB equations
\begin{equation}\label{qk2}
v(z+p\delta_j)=K^\vee_j(z,\tau,p)\,v(z),\qquad j=1,\dots n,
\end{equation}
is also compatible. In fact, if
we write $x^\vee=(x_n,\dots,x_1)$ for any $x=(x_1,\dots,x_n)\in \C^n$
and let $P:V_{\vec\Lambda}\to V_{\vec\Lambda^\vee}$
be the linear map sending $v_1\otimes\cdots\otimes v_n$ to $v_n\otimes
\cdots\otimes v_1$,
then we have
\[
K^\vee_i(z,\tau,p,\vec\Lambda)=
P^{-1}K_{n+1-i}(z^\vee,\tau,p,\vec\Lambda^\vee)P.
\]

\subsection{Fundamental hypergeometric solutions}\label{sol.qKZB}
In \cite{FTV2} we constructed a {\it fundamental hypergeometric solution} of
the qKZB equations: it is constructed out of
 a {\it universal hypergeometric function}
$u(z,\lambda,\mu,\tau,p)$ taking
values in $V_{\vec\Lambda}[0]\otimes V_{\vec\Lambda}[0]$ and obeying the
equations
\begin{eqnarray}\label{eq1}
u(z+\delta_jp,\tau,p)&=&K_j(z,\tau,p)\otimes D_j\, u(z,\tau,p),\notag\\
u(z+\delta_j\tau,\tau,p)&=&D^\vee_j\otimes K^\vee_j(z,\tau,p)\, u(z,\tau,p),\\
u(z+\delta_j,\tau,p)&=&u(z,\tau,p).
\notag
\end{eqnarray}
Here we view $u$ as taking values in the space of functions of
$\lambda$ and $\mu$ with values in  $V_{\vec\Lambda}[0]\otimes V_{\vec\Lambda}[0]$.
$K_j$ acts on the variable $\lambda$ and $K_j^\vee$ on the variable
$\mu$. The operators $D_j$, $D_j^\vee$ act by multiplication by
diagonal matrices $D_j(\mu)$, $D_j^\vee(\lambda)$, respectively.
For our purpose, the most convenient description of these matrices
is in terms of the function
\[
\alpha(\lambda)=\exp(-{\pi i\lambda^2/4\eta}).
\]
We have, for $j=1,\dots, n$,
\begin{eqnarray*}
D_j(\mu)&=&
\frac
{\alpha(\mu-2\eta(h^{(j+1)}+\cdots+h^{(n)}))}
{\alpha(\mu-2\eta(h^{(j)}+\cdots+h^{(n)}))}
\, e^{\pi i\eta\Lambda_j(\sum_{l=1}^{j-1}\Lambda_l-
                \sum_{l=j+1}^{n}\Lambda_l)},
\\
D^\vee_j(\lambda)&=&
\frac
{\alpha(\lambda-2\eta(h^{(1)}+\cdots+h^{(j-1)}))}
{\alpha(\lambda-2\eta(h^{(1)}+\cdots+h^{(j)}))}
\, e^{-\pi i\eta\Lambda_j(\sum_{l=1}^{j-1}\Lambda_l-
                \sum_{l=j+1}^{n}\Lambda_l)}.
\end{eqnarray*}
The {\it fundamental hypergeometric
solution} is then 
\bean\label{fund}
v=(1\otimes\prod_{j=1}^n D_j^{-z_j/p})\,u.
\eean
It obeys the qKZB equations in the first factor:
\bean\label{first-factor}
v(z+p\delta_j,\tau,p)=K_j(z,\tau,p)\otimes 1\,v(z,\tau,p).
\eean
In other words, for every complex $\mu$ and every linear form on
$V_{\vec\La} [0]$, we have a solution of the qKZB equations. 
We call such solutions {\it the elementary hypergeometric solutions}.
A solution $w(z,\la)$ of the qKZB equations  is called {\it a hypergeometric solution}
if it can be represented in the form
\bean\label{hypergeom}
w(z,\la)\,=\,\sum_j\, a_j(z)\, v_j(z,\la)
\eean
where $v_j$ are elementary hypergeometric solutions and $a_j$ are
scalar functions of $z_1,...,z_n$ periodic with respect to $p$-shifts of
variables.

We recall the construction of the universal
hypergeometric function in Sec. \ref{hypergsol}.

The second system of equations in \Ref{eq1} gives the monodromy of these
solutions, see \cite{FTV2}.

\subsection{Finite-dimensional representations}
If $\Lambda$ is a non-negative integer, $V_\Lambda$ contains the subspace
$SV_\Lambda=\oplus_{j=\Lambda+1}^\infty\C e_j$ with the property that, for any
$\Mu$, $SV_\Lambda\otimes V_\Mu$ and $V_\Mu\otimes SV_\Lambda$ are preserved by
the $R$-matrices $R_{\Lambda,\Mu}(z,\lambda)$ and $R_{\Mu,\Lambda}(z,\lambda)$,
respectively, see \cite{FV1} and Sec. \ref{fs}. Let
$L_\Lambda=V_\Lambda/SV_\Lambda$, $\Lambda\in\Z_{\geq 0}$. Then, in particular,
for any non-negative integers $\Lambda$ and $\Mu$, $R_{\Lambda,\Mu}(z,\lambda)$
induces a map, also denoted by $R_{\Lambda,\Mu}(z,\lambda)$, on the
finite-dimensional space $L_\Lambda\otimes L_\Mu$.

The simplest nontrivial case is $\Lambda=\Mu=1$. Then $R_{1,1}(z,\lambda)$
is defined on a four-dimensional vector space and coincides with
the {\em fundamental} $R$-matrix, the matrix of structure constants of
the elliptic quantum group $E_{\tau,\eta}(sl_2)$.

If $\Lambda_1,\dots,\Lambda_n$ are non-negative integers, we can
consider the qKZB equations \Ref{qk1}, \Ref{qk2} on functions with values in the zero
weight space of $L_{\Lambda_1}\otimes \cdots\otimes L_{\Lambda_n}$.


\section{Function spaces and R-matrices}\label{fs}

In this section we geometrically construct R-matrices of the qKZB equations.


Let us fix complex parameters $\tau$, $\eta$ with $\Im\tau>0$, and complex
numbers $\Lambda_1,\dots, \Lambda_n$. We set $a_i=\eta\Lambda_i$,
$i=1,\dots,n$.

\subsection{A space of symmetric functions, \cite {FTV1}}\label{sssym}
Introduce a space of functions with an action of the symmetric group.
Recall that the Jacobi theta function
\begin{equation}
\theta(t)=\theta(t,\tau)=-\sum_{j\in\Z}
e^{\pi i(j+\half)^2\tau+2\pi i(j+\half)(t+\half)},
\notag
\end{equation}
has multipliers $-1$ and $-\exp(-2\pi it-\pi i\tau)$ as $t\to t+1$ and
$t\to t+\tau$, respectively. It is an odd entire function whose zeros are
simple and lie on the lattice $\Z+\tau \Z$. It has the product formula
\be
\theta(t)\,=\,2e^{\pi i\tau/4}\sin(\pi t)
\prod_{j=1}^\infty(1-q^j)(1-q^je^{2\pi it})(1-q^je^{-2\pi it}),
\qquad q=e^{2\pi i\tau}.
\ee

For complex numbers $a_1,\dots,a_n$, $z_1,\dots,z_n$, $\lambda$,
let $\tilde F^m_{a_1,\dots,a_n}(z_1,\dots,z_n,\lambda)$ be the space of
meromorphic functions $f(t_1,\dots,t_m)$ of $m$ complex variables such that
\begin{enumerate}
\item[(i)] $\prod_{i<j}\theta(t_i-t_j+2\eta)
\prod_{i=1}^m\prod_{k=1}^n\theta(t_i-z_k-a_k)f$
is a holomorphic function on $\C^m$.
\item[(ii)] $f$ is periodic with period 1 in each of its arguments and
\be
f(\cdots,t_j+\tau,\cdots)
=e^{-2\pi i(\lambda+4\eta j-2\eta)}f(\cdots,t_j,\cdots),
\ee
for all $j=1,\dots,m$.
\end{enumerate}

The symmetric group $S_m$ acts on
$\tilde F^m_{a_1,\dots,a_n}(z_1,\dots,z_n,\lambda)$ so that the transposition
of $j$ and $j+1$ acts as
\be
s_jf(t_1,\dots,t_m)=f(t_1,\dots,t_{j+1},t_j,\dots,t_m)
\frac{\theta(t_{j}-t_{j+1}-2\eta)}
{\theta(t_{j}-t_{j+1}+2\eta)}\,.
\ee

For any $m\in\Z_{>0}$, let $F^m_{a_1,\dots,a_n}(z_1,\dots,z_n,\lambda)=
\tilde F^m_{a_1,\dots,a_n} (z_1,\dots,z_n,\lambda)^{S_m}$ be the space of
$S_m$-invariant functions. If $m=0$, then we set
$F^0_{a_1,\dots,a_n}(z_1,\dots,z_n,\lambda)=\C$. We denote by $\Sym$
the symmetrization operator $\Sym=\sum_{s\in S_m}s:\tilde F^m\to F^m$.
Also, we set
\be
F_{a_1,\dots,a_n}(z_1,\dots,z_n,\lambda)
=\oplus_{m=0}^\infty F^m_{a_1,\dots,a_n}(z_1,\dots,z_n,\lambda),
\ee
and define an $\h$-module structure on
$F_{a_1,\dots,a_{n}}(z_1,\dots,z_{n},\lambda)$ by letting $h$ act by
\be
h|_{F^{m}_{a_1,\dots,a_{n}}(z_1,\dots,z_{n},\lambda)}=(\sum_{i=1}^n
\Lambda_i-2m)\Id,\qquad a_i=\eta\Lambda_i.
\ee

Clearly, $F^m_{a_{\sigma(1)},\dots,a_{\sigma(n)}}(z_{\sigma(1)},\dots,
z_{\sigma(n)})=F^{m}_{a_1,\dots,a_{n}}(z_1,\dots,z_{n},\lambda)$
for any permutation $\sigma\in S_n$.


$F^{m}_{a_1,\dots,a_{n}}(z_1,\dots,z_{n},\lambda)$
is a finite-dimensional vector space of dimension
$\left(\begin{matrix}{n+m-1}\\{m}\end{matrix}\right)$.

\begin{example}
Let $n=1$. Then $F^m_a(z,\lambda)$ is a one-dimensional space spanned by
\begin{equation}\label{eb}
\omega_m(t_1,\dots,t_m,\lambda,z)=
\prod_{i<j}\frac{\theta(t_i-t_j)}
{\theta(t_i-t_j+2\eta)}\prod_{j=1}^{m}\frac
{\theta(\lambda+2\eta m+t_j-z-a)}
{\theta(t_j-z-a)}.
\end{equation}
\end{example}

\subsection{Tensor products, \cite{FTV1}}\label{sstp}
Let $n=n'+n''$, $m=m'+m''$ be non-negative integers and $a_1,\dots,a_n$,
$z_1,\dots,z_n$ be complex numbers. The formula
\be
k(t_1,\dots,t_{m})\,=\,
\frac1{m'!m''!}\;\Sym\biggl( f(t_{1},\dots,t_{m'})g(t_{m'+1},\dots,t_{m})
\prod_{
\begin{matrix} \scriptstyle{m'<j\leq m}\\
\scriptstyle{1\leq l\leq n'}\end{matrix}}
\frac{\theta(t_j-z_l+a_l)}{\theta(t_j-z_l-a_l)}
\biggr)\ee
correctly defines a linear map $\Phi:f\otimes g\mapsto k=\Phi(f\otimes g)$,
\bea
\oplus_{m'=0}^{m}
F^{m'}_{a_1,\dots,a_{n'}}(z_1,\dots,z_{n'},\lambda)\otimes
F^{m''}_{a_{n'+1},\dots,a_{n}}(z_{n'+1},\dots,z_{n},\lambda-2\nu)
\\
\to F^{m}_{a_1,\dots,a_{n}}(z_1,\dots,z_{n},\lambda),
\eea
\bea
\oplus_{m'=0}^{m}F^{m'}_{a_1,\dots,a_{n'}}(z_1,\dots,z_{n'},\lambda)\otimes
F^{m''}_{a_{n'+1},\dots,a_{n}}(z_{n'+1},\dots,z_{n},\lambda-2\nu)\to
F^{m}_{a_1,\dots,a_{n}}(z_1,\dots,z_{n},\lambda),
\eea
where
$
\nu=a_{1}+\cdots +a_{n'}-2\eta m'.
$
For generic values of the parameters $z_j$, $\lambda$, the map $\Phi$ is
an isomorphism. Moreover, $\Phi$ is associative in the sense that, for any
three functions $f,g,h$, $\Phi(\Phi(f\otimes g)\otimes h)=
\Phi(f\otimes\Phi(g\otimes h))$, whenever defined.

By iterating this construction, we get for all $n\geq1$ a linear map $\Phi_n$,
defined recursively by $\Phi_1=\Id$, $\Phi_{n}=\Phi(\Phi_{n-1}\otimes\Id)$,
from
\be
\oplus_{m_1+\cdots+m_n=m}\otimes_{i=1}^n
F_{a_i}^{m_i}(z_i,\lambda-2\eta(\mu_{1}+\cdots+\mu_{i-1}))
\ee
to $F^m_{a_1,\dots,a_n}(z_1,\dots,z_n,\lambda)$, with $\mu_j=a_j/\eta-2m_j$,
$j=1,\dots,n$.

Let $V_\Lambda^*=\oplus_{j=0}^\infty\C e_j^*$ be the restricted dual of
the module $V_\Lambda=\oplus_{j=0}^\infty \C e_j$. It is spanned by the basis
$(e_j^*)$ dual to the basis $(e_j)$. We let $\h$ act on $V_\Lambda^*$ by
$he_j^*=(\Lambda-2j)e_j^*$. Then the map that sends $e^*_j$ to $\omega_j$
(see \Ref{eb}) defines an isomorphism of $\h$-modules
\be
\omega(z,\lambda):V_\Lambda^*\to F_a(z,\lambda), \qquad a=\eta\Lambda.
\ee
By composing this with the maps $\Phi$, we obtain
homomorphisms (of $\h$-modules)
\be
\omega(z_1,\dots,z_n,\lambda):
V_{\Lambda_1}^*\otimes\cdots\otimes
V_{\Lambda_n}^*\to F_{a_1,\dots,a_n}(z_1,\dots,z_n,\lambda)
\ee
which are isomorphisms for generic values of $z_1,\dots,z_n,\lambda$.
The restriction of the map $\omega(z_1,\dots,z_n,\lambda)$ to
$ e_{m_1}^*\otimes\cdots\otimes e_{m_n}^*$ is
\be
\Phi_n(\omega(z_1,\lambda) e^*_{m_1}
\otimes
\omega(z_{2},\lambda-2\eta\mu_1) e^*_{m_2}
\otimes
\cdots\otimes\omega(z_n,\lambda-2\eta(\mu_1+\cdots+\mu_{n-1}))
e^*_{m_n}) ,
\ee
where $\mu_j=\Lambda_j-2m_j$, $j=1,\dots,n$. 

For example, if $n=2$,
then $\omega(z_1,z_2,\lambda)$ sends $e^*_j\otimes e^*_k$ to
\be
\frac1{j!k!}\;\Sym\biggl(\omega_j(t_1,\dots,t_{j},\!\lambda,z_1)
\omega_k(t_{j+1},\dots,t_{j+k},\!\lambda-2a_1+4\eta j,z_2)\!
\prod_{i=j+1}^{j+k}\!\frac{\theta(t_i-z_1+a_1)}{\theta(t_i-z_1-a_1)}\biggr),
\ee
where $\{\omega_j(t_1,\dots,t_j,\lambda,z)\}$ is the basis \Ref{eb}
of $F_a(z,\lambda)$.

More generally, we have an explicit formula for the image of
$e_{m_1}^*\otimes\cdots\otimes e_{m_n}^*$, which we discuss next.

\subsection{Bases of $F_{a_1,\dots,a_n}(z_1,\dots,z_n,\la)$, \cite{FTV1}}\label{bases}
The space $V_\Lambda$ comes with a basis $e_j$. Thus we have the natural basis
$e_{m_1}^*\otimes\cdots\otimes e_{m_n}^*$ of the tensor product of
$V^*_{\Lambda_i}$ in terms of the dual bases of the factors. The map
$\omega(z_1,\dots,z_n,\lambda)$ maps, for generic $z_i$, this basis to a basis
of $F_{a_1,\dots, a_n}(z_1,\dots,z_n,\lambda)$, which is an essential part of
our formulae for integral representations for solutions of the qKZB equations.

We give here an explicit formula for the basis vectors.
Let $m\in\Z_{\geq 0}$, $\vec\Lambda=(\Lambda_1,\dots,\Lambda_n)\in\C^n$, and let
$z=(z_1,\dots,z_n)\in\C^n$ be generic. Set $a_i=\eta\Lambda_i$. Let
\be
u(t_1,\dots,t_m)=\prod_{i<j}
\frac{\theta(t_i-t_j+2\eta)}
{\theta(t_i-t_j)}
\ee
Then, for generic $\lambda\in\C$, the functions
\be
\omega_{m_1,\dots,m_n}(t_1,\dots,t_m,\lambda,z)\,=\,\omega(z,\lambda)
\,e_{m_1}^*\otimes\cdots\otimes e_{m_n}^*
\ee
labeled by $\,{m_1,\dots,m_n\in\Z}\,$ with $\,{\sum_km_k=m}\,$ form a basis
of $F^m_a(z,\lambda)$ and are given by the explicit formula
\bean \label{w.f}
 \omega_{m_1,\dots,m_n}(t_1,\dots,t_m,\lambda,z,\tau)\,=\,
u(t_1,\dots,t_m)^{-1}
\sum_{I_1,\dots,I_n}\prod_{l=1}^n\prod_{i\in I_l}\prod_{k=1}^{l-1}
\frac{\theta(t_i-z_k+a_k)}{\theta(t_i-z_k-a_k)}
\\
\times\,\prod_{k<l}\prod_{i\in I_k,j\in I_l}
\frac{\theta(t_i-t_j+2\eta)}{\theta(t_i-t_j)}
\prod_{k=1}^{n}\prod_{j\in I_k}
\frac{\theta(\lambda\!+\!t_j\!-\!z_k\!-\!a_k\!+\!2\eta
m_k\!-\!2\eta\sum_{l=1}^{k-1}(\Lambda_l\!-\!2m_l))}{\theta(t_j-z_k-a_k)}\,.
\notag
\eean
The summation is over all $n$-tuples $I_1,\dots,I_n$ of disjoint subsets of
$\{1,\dots,m\}$ such that $I_j$ has $m_j$ elements, $1\leq j\leq n$.

We shall call the functions
$\omega_{m_1,\dots,m_n}(t_1,\dots,t_m,\lambda,z,\tau)$ {\it the weight functions.}

For any permutation $\sigma \in S_n$ the space $F_{a_{\sigma(1)},...,a_{\sigma(n)}}
(z_{\sigma(1)},...,z_{\sigma(n)},\la)$ coincides with 
$F_{a_1,\dots, a_n}(z_1,\dots,z_n,\lambda)$.
Thus, for any $\sigma$, we have a map
\be
\omega^\sigma (z_1,\dots,z_n,\lambda):
V_{\Lambda_{\sigma (1)}}^*\otimes\cdots\otimes
V_{\Lambda_{\sigma (n)}}^*\to F_{a_1,\dots,a_n}(z_1,\dots,z_n,\lambda)
\ee
and the corresponding basis in $F_{a_1,\dots,a_n}(z_1,\dots,z_n,\lambda)$.
The basis corresponding to the permutation $\sigma$, such that $\sigma(j)=
n-j+1$ for all $j$, is denoted
$\tilde \omega_{m_1,\dots,m_n}(t_1,\dots,t_m,\lambda,z,\tau)$.
We have
\bean \label{m.w.f}
 \tilde \omega_{m_1,\dots,m_n}(t_1,\dots,t_m,\lambda,z,\tau)\,=\,
u(t_1,\dots,t_m)^{-1}
\sum_{I_1,\dots,I_n}\prod_{l=1}^n\prod_{i\in I_l}\prod_{k=l+1}^{n}
\frac{\theta(t_i-z_k+a_k)}{\theta(t_i-z_k-a_k)}
\\
\times\,\prod_{k>l}\prod_{i\in I_k,j\in I_l}
\frac{\theta(t_i-t_j+2\eta)}{\theta(t_i-t_j)}
\prod_{k=1}^{n}\prod_{j\in I_k}
\frac{\theta(\lambda\!+\!t_j\!-\!z_k\!-\!a_k\!+\!2\eta
m_k\!-\!2\eta\sum_{l=k+1}^{n}(\Lambda_l\!-\!2m_l))}{\theta(t_j-z_k-a_k)}\,.
\notag
\eean
We call these functions {\it the mirror weight functions.}

\subsection{Geometric construction of $R$-matrices, \cite {FTV1}}
\label{geom}
Let $a=\eta\Lambda$ and $b=\eta \Mu$ be complex numbers. Since
$F_{ab}(z,w,\lambda)$ coincides with $F_{ba}(w,z,\lambda)$,
 we obtain a family of isomorphisms between
$V^*_\Lambda\otimes V^*_\Mu$ and $V^*_\Mu\otimes V^*_\Lambda $. The composition
of this family with the flip
$P:V^*_\Mu\otimes V^*_\Lambda\to V^*_\Lambda\otimes V^*_\Mu$,
$Pv\otimes w=w\otimes v$ gives a family of automorphisms of
$V^*_\Lambda\otimes V^*_\Mu$:

Let $z,w,\lambda$ be such that
$\omega(z,w,\lambda,\La, \Mu):V^*_\Lambda \otimes V^*_\Mu\to F_{ab}(z,w,\lambda)$ is
invertible. The {\em $R$-matrix}
$R_{\Lambda,\Mu}(z,w,\lambda)\in\End_\h(V_\Lambda \otimes V_\Mu)$ is the dual
map to the composition $R^*_{\Lambda,\Mu}(z,w,\lambda)$:
\be
V^*_\Lambda\otimes V^*_\Mu
\buildrel{P}\over{\longrightarrow}
V^*_\Mu\otimes V^*_\Lambda
\buildrel{\omega(w,z,\lambda,\Mu, \La)}\over{\longrightarrow}
F_{ab}(z,w,\lambda)
\buildrel{\omega(z,w,\lambda,\La,\Mu)^{-1}}\over{\longrightarrow}
V^*_\Lambda \otimes V^*_\Mu,
\ee
where we identify canonically $V_\Lambda ^*\otimes V_\Mu^*$ with
$(V_\Lambda \otimes V_\Mu)^*$.

Alternatively, the $R$-matrix $R_{\Lambda,\Mu}(z,w,\lambda)$ can be thought
of as the transition matrix expressing the basis of mirror weight functions
$\tilde\omega_{ij}=\omega(w,z,\lambda,\Mu,\La)e_j^*\otimes e_i^*$
of the space $F_{ab}(z,w,\lambda)$
in terms of the basis of weight functions
$\omega_{ij}=\omega(z,w,\lambda,\La,\Mu)e_i^*\otimes e_j^*$:
if $R_{\Lambda,\Mu}(z,w,\lambda)e_i\otimes e_j=
\sum_{kl}R_{ij}^{kl}e_k\otimes e_l$, then
\begin{equation} \label{R.c}
\tilde\omega_{kl}=\sum_{ij}R_{ij}^{kl}\omega_{ij}.
\end{equation}

 $R_{\Lambda, \Mu}(z,w,\lambda)$ is a meromorphic function
of $\Lambda,\Mu,z,w,\lambda$.
 If $\Lambda$ is generic, then $R_{\Lambda,\Lambda}(z,w,\lambda)$
is regular at $z=w$ and $\lim_{z\to w}R_{\Lambda,\Lambda}(z,w,\lambda)=P$,
where $P$ is the flip $u\otimes v\mapsto v\otimes u$.
$R_{\Lambda, \Mu}(z,w,\lambda)$ depends only on the difference $z-w$.
Accordingly, we write $R_{\Lambda,\Mu}(z-w,\lambda)$ instead of
$R_{\Lambda,\Mu}(z,w,\lambda)$ in what follows.

The R-matrices  $R_{\Lambda,\Mu}(z,\lambda)$ satisfy the dynamical Yang-Baxter equation, 
they obey I--III of Sec. \ref{sqkzb}.

Consider the case of positive integer weights. In this case
the $R$-matrices have invariant subspaces. If $\Lambda\in\Z_{\geq 0}$
we let $SV_\Lambda$ be the subspace of $V_\Lambda$ spanned by
$e_{\Lambda+1},\,e_{\Lambda+2},\dots$. The $\Lambda+1$-dimensional quotient
$V_\Lambda/SV_\Lambda$ will be denoted  $L_\Lambda$, and will be often
identified with $\oplus_{j=0}^m\C e_j$.

Let $z,\eta,\lambda$ be generic and $\Lambda$, $\Mu\in\C$.
\begin{enumerate}
\item[(i)] If $\Lambda\in\Z_{\geq 0}$, then $R_{\Lambda,\Mu}(z,\lambda)$
preserves $SV_\Lambda\otimes V_\Mu$
\item[(ii)] If $\Mu\in\Z_{\geq 0}$, then $R_{\Lambda,\Mu}(z,\lambda)$
preserves $V_\Lambda\otimes SV_\Mu$
\item[(iii)] If $\Lambda\in\Z_{\geq 0}$ and $\Mu\in\Z_{\geq 0}$,
then $R_{\Lambda,\Mu}(z,\lambda)$ preserves
$SV_\Lambda\otimes V_\Mu+V_\Lambda\otimes SV_\Mu$.
\end{enumerate}

In particular, if $\Lambda$ and/or $\Mu$ are non-negative integers,
then $R_{\Lambda,\Mu}(z,\lambda)$ induces operators, still denoted by
$R_{\Lambda,\Mu}(z,\lambda)$, on the quotients $L_\Lambda\otimes V_\Mu$,
$V_\Lambda\otimes L_\Mu$ and/or $L_\Lambda\otimes L_\Mu$.
They obey the dynamical Yang--Baxter equation.

\subsection{Evaluation Verma modules and their tensor products}
\label{ssevm}
Here we recall the relation between the geometric construction of
tensor products and $R$-matrices and the representation theory of
the elliptic quantum group $E_{\tau,\eta}(sl_2)$ \cite{FV1}.

Recall the definition of a representation of $E_{\tau,\eta}(sl_2)$:
let $\h$ act on $\C^2$ via $h={\mathrm {diag}}(1,-1)$. A representation of
$E_{\tau,\eta}(sl_2)$ is an $\h$-module $W$ with diagonalizable action of $h$
and finite-dimensional eigenspaces, together with an operator
$L(z,\lambda)\in\End(\C^2\otimes W)$ (the ``$L$-operator''),
commuting with $h^{(1)}+h^{(2)}$, and obeying the relations
\bea
R^{(12)}(z-w,\lambda-2\eta h^{(3)}) &&\kern-20pt\, L^{(13)}(z,\lambda)\,
L^{(23)}(w,\lambda-2\eta h^{(1)})
\\[2pt]
{}={}&&\kern-20pt\, L^{(23)}(w,\lambda)\,
L^{(13)}(z,\lambda-2\eta h^{(2)})\,R^{(12)}(z-w,\lambda)
\eea
in $\End(\C^2\otimes\C^2\otimes W)$. The {\em fundamental $R$-matrix}
$R(z,\lambda)\in\End(\C^2\otimes\C^2)$ is the following solution of the
dynamical Yang--Baxter equation: let $e_0$, $e_1$ be the standard basis of
$\C^2$, then with respect to the basis $e_0\otimes e_0$, $e_0\otimes e_1$,
$e_1\otimes e_0$, $e_1\otimes e_1$ of $\C^2\otimes\C^2$,
\be
R(z,\lambda)=
\left(\begin{matrix}
1 & 0 & 0 & 0\\
0 & \alpha(z,\lambda) & \beta(z,\lambda) & 0\\
0 & \beta(z,-\lambda) & \alpha(z,-\lambda) & 0\\
0 & 0 & 0 & 1
\end{matrix}\right).
\ee
where
\be
\alpha(z,\lambda)=
\frac{\theta(\lambda+2\eta)\theta(z)}{\theta(\lambda)\theta(z-2\eta)},\qquad
\beta(z,\lambda)=
-\frac{\theta(\lambda+z)\theta(2\eta)}{\theta(\lambda)\theta(z-2\eta)}.
\ee

\begin{thm}
\label{tfrm} \cite{FTV1}
Let us identify the two-dimensional space $L_1=V_{\Lambda=1}/SV_{\Lambda=1}$
with $\C^2$ via the basis $e_0,e_1$. Then the $R$-matrix
$R_{1,1}(z,\lambda)\in\End(L_1\otimes L_1)$ coincides with
the fundamental $R$-matrix.
\end{thm}

\begin{corollary}
For any $w,\Mu\in\C$, the $\h$-module $V_\Mu$ together with the operator
$L(z,\lambda)=R_{1,\Mu}(z-w,\lambda)\in\End(L_1\otimes V_\Mu)$ defines
a representation of $E_{\tau,\eta}(sl_2)$.
\end{corollary}

This representation is called in \cite{FV1} the evaluation Verma module with
evaluation point $w$ and highest weight $\Mu$. It is denoted  $V_\Mu(w)$.

The tensor product construction of \ref{sstp} is related to
the tensor product of representations of the elliptic quantum group.
Recall that if $W_1$, $W_2$ are representations of the elliptic quantum group
with $L$-operators $L_1(z,\lambda)$, $L_2(z,\lambda)$, then
their tensor product $W=W_1\otimes W_2$ with $L$-operator
\be
L(z,\lambda)=L_1(z,\lambda-2\eta h^{(3)})^{(12)}
L_2(z,\lambda)^{(13)}\in\End(\C^2\otimes W)
\ee
is also a representation of the elliptic quantum group.

\begin{thm}\cite{FTV1}

Let $\Lambda_1,\dots,\Lambda_n\in\C$ and $z_1,\dots,z_n$ be generic complex
numbers. Let $V=V_{\Lambda_1}\otimes\cdots\otimes V_{\Lambda_n}$ and
$L(z,\lambda)\in\End(V_{\Lambda=1}\otimes V)$ be defined by the relation
\be
\omega(z,z_1,\dots,z_n,\lambda)L(z,\lambda)^*=
\omega(z_1,\dots,z_n,z,\lambda)P
\ee
in $\End((V_1\otimes V)^*)=\End(V_1^*\otimes V^*)$, where
$Pv_1\otimes v=v\otimes v_1$, if $v_1\in V_1^*$, $v\in V^*$.
Then $L(z,\lambda)$ is well-defined as an endomorphism of the quotient
$L_1\otimes V=\C^2\otimes V$, and defines a structure of a representation of
$E_{\tau,\eta}(sl_2)$ on $V$. This representation is isomorphic to the tensor
product of evaluation Verma modules
\be
V_{\Lambda_n}(z_n)\otimes\cdots\otimes V_{\Lambda_1}(z_1),
\ee
with the isomorphism
$u_1\otimes\cdots\otimes u_n\mapsto u_n\otimes\cdots\otimes u_1$.
\end{thm}

Finally, the dynamical Yang--Baxter equation in
$L_1\otimes V_\Lambda\otimes V_\Mu$ can be stated as saying
that $R_{\Lambda,\Mu}(z-w,\lambda)P$ is an isomorphism from
$V_\Mu(w)\otimes V_\Lambda(z)$ to $V_\Lambda(z)\otimes V_\Mu(w)$,
see \cite{FV1}.

\subsection{Weight functions}
For a natural number $m$, define
\bean
\Z^2_m =\{ (m_1,m_2)\in \Z^2_{\geq 0}\,|\, m_1+m_2=m\}.
\eean
Define a lexicographical order on $\Z^2_m$: ${}$
$(l_1,l_2)>(m_1,m_2)$ if $l_1>m_1$.

Let $z=(z_1,z_2)\in \C^2$. Define the points $T_M \in \C^m,\, M\in \Z^2_{m},$
by
\bean
T_M=(z_1-\eta\La_1+2\eta (m_1-1), z_1-\eta\La_1+2\eta (m_1-2), ...,
z_1-\eta\La_1,
\\
z_2-\eta\La_2+2\eta (m_2-1), z_2-\eta\La_2+2\eta (m_2-2),
....,z_2-\eta\La_2).
\notag
\eean

\begin{lemma}\label{triang}
Let $\omega_M(t_1,...,t_m,\la,z,\tau),\,M\in \Z^2_m$, be the weight functions associated with
$\vec\La=(\La_1,\La_2)$ and defined in \Ref{w.f}. Then for
all $L,M\in \Z^2_m,\, L>M$, we have
\bean
\omega_M(T_L,\la,z,\tau)=0.
\notag
\eean
\end{lemma}
\begin{proof} According to \Ref{w.f},
\bean\label{w.f.2}
\omega_{m_1,m_2}(t_1,...,t_m,\la,z,\tau)
\,=\,
u(t_1,\dots,t_m)^{-1}
\sum_{I_1,I_2}
\prod_{i\in I_1,j\in I_2}
\frac{\theta(t_i-t_j+2\eta)}{\theta(t_i-t_j)}\,\times
\\
\prod_{i\in I_1}
\frac{\theta(\lambda\!+\!t_i\!-\!z_1\!-\!\eta\La_1\!+\!2\eta
m_1)}{\theta(t_i-z_1-\eta\La_1)}\,
\prod_{j\in I_2}
\frac{\theta(\lambda\!+\!t_j\!-\!z_2\!-\!\eta\La_2\!+\!2\eta
m_2\!-\!2\eta(\Lambda_1\!-\!2m_1))}{\theta(t_j-z_2-\eta\La_2)}\,
\frac{\theta(t_j-z_1+\eta\La_1)}{\theta(t_j-z_1-\eta\La_1)}
\notag
\eean
where the summation is over all pairs $I_1, I_2$ of disjoint subsets of
$\{1,\dots,m\}$ such that $I_j$ has $m_j$ elements, $j=1,2$.

It is convenient to set apart the term of the sum 
in \Ref{w.f.2} corresponding to the distinguished partition
$I_1=\{1,...,m_1\}$, $I_2= \{ m_1+1,...,m_2\}$,
\bean\label{dist}
 \prod_{1\leq i<j \leq m_1}
\frac{\theta(t_i-t_j)}{\theta(t_i-t_j+2\eta)}\,
\, \prod_{m_1+1\leq i<j \leq m}
\frac{\theta(t_i-t_j)}{\theta(t_i-t_j+2\eta)}\,
\prod_{i=1}^{m_1}\,
\frac{\theta(\lambda\!+\!t_i\!-\!z_1\!-\!\eta\La_1\!+\!2\eta m_1)}
{\theta(t_i-z_1-\eta\La_1)}\,
\notag
\\
\times \, \prod_{i=m_1+1}^{m}\,
\frac{\theta(\lambda\!+\!t_i\!-\!z_2\!-\!\eta\La_2\!+\!2\eta m_2
-2\eta(\La_1-2m_1))}
{\theta(t_i-z_2-\eta\La_2)}\,
\frac{\theta(t_i-z_1+\eta\La_1)}{\theta(t_i-z_1-\eta\La_1)}\,.
\eean
This term  will be called {\it  distinguished.}

Let us show that $\omega_M(T_L)=0$ for $L>M$. Consider a term 
in the sum in \Ref{w.f.2} corresponding
to a partition $I_1, I_2$. Since $l_1>m_1$, at least one of the 
numbers $1,2,...,l_1$ belong to $I_2$. If $l_1\in I_2$, then
the factor $\theta(t_{l_1}-z_1+\eta\La_1)$ in \Ref{w.f.2} is zero.
If $l_1$ does not belong to $I_2$, then there is a pair of numbers
$i, i+1 < l_1$ such that $i\in I_2$ and $i+1\in I_1$. In that case
the factor $\theta(t_{i+1}-t_i+2\eta)$ of the first product in \Ref{w.f.2}
equals zero. 
\end{proof}

\begin{lemma}\label{diag}
${}$

\begin{enumerate}
\item[I.]
Let $L,M\in \Z^2_m,\, L\leq M$. According to \Ref{w.f.2}, decompose 
$\omega_M(T_L,\la,z,\tau)$ into the sum of the terms corresponding to
partitions $I_1, I_2$. Then all of the not distinguished terms are equal to zero.
\item[II.]
For any $M\in \Z^2_m$, we have
\bean\label{w.f.3}
\omega_M(T_M,\la,z,\tau)= 
\prod_{l=1}^{m_1}\,[\,
\frac{\theta(2\eta)}{\theta(2\eta l)}\,
\frac{\theta(\lambda\!-2\eta(\La_1 -m_1-l+1))}
{\theta(-2\eta(\La_1 - l + 1))}\,]\,
\\
\times\,\prod_{l=1}^{m_2}\,[\,
\frac{\theta(2\eta)}{\theta(2\eta l)}\,
\frac{\theta (\la - 2\eta (\La_1+\La_2-2m_1 -m_2 - l + 1))}
{\theta( -2\eta (\La_2 - l + 1))}\,
\notag
\\
\times\,\frac{\theta (-z_1+z_2 -\eta\La_1+\eta \La_2 +2\eta (l-1))}
{\theta (-z_1+z_2 -\eta\La_1-\eta \La_2 +2\eta (l-1))}\,]\,.
\notag
\eean
\end{enumerate}
\end{lemma}
\begin{proof} 
The values of all of the not distinguished summands in \Ref{w.f.2} equal zero 
for the same reason as in Lemma \ref{triang}. Namely, if $l_1 \in I_2$,
then the factor $\theta(t_{l_1}-z_1+\eta\La_1)$ in \Ref{w.f.2} is zero.

Now let $l_1\in I_1$. If there is $j<l_1$ such that $j\in I_2$,
then there is a pair of numbers
$i, i+1 < l_1$ such that $i\in I_2$ and $i+1\in I_1$. In that case
the factor $\theta(t_{i+1}-t_i+2\eta)$ of the first product in \Ref{w.f.2}
equals zero. 

So assume that $\{1,...,l_1\}\in I_1$ and the partition $I_1, I_2$ is not
distinguished, then there is a pair of numbers
$l_1<i, i+1 $ such that $i\in I_2$ and $i+1\in I_1$. In that case
the factor $\theta(t_{i+1}-t_i+2\eta)$ of the first product in \Ref{w.f.2}
equals zero. This proves the first part of the Lemma.

The RHS in \Ref{w.f.3} is the value of the distinguished summand
in \Ref{w.f.2}. This proves the second part.
\end{proof}

\begin{lemma}\label{w.f.mir}
Let $\tilde\omega_M(t_1,...,t_m,\la,z,\tau),\,M\in \Z^2_m$, be the 
mirror weight functions associated with
$\vec\La=(\La_1,\La_2)$ and defined in \Ref{w.f}. Then for
all $L,M\in \Z^2_m,\, L<M$, we have
\bean
\tilde\omega_M(T_L,\la,z,\tau)=0.
\notag
\eean
Moreover, for any $M\in \Z^2_m$, we have
\bean\label{w.f.4}
\tilde\omega_M(T_M,\la,z,\tau)= 
\prod_{l=1}^{m_2}\,[\,
\frac{\theta(2\eta)}{\theta(2\eta l)}\,
\frac{\theta(\lambda\!-2\eta(\La_2 -m_2-l+1))}
{\theta(-2\eta(\La_2 - l + 1))}\,]\,
\\
\times\,\prod_{l=1}^{m_1}\,[\,
\frac{\theta(2\eta)}{\theta(2\eta l)}\,
\frac{\theta (\la - 2\eta (\La_1+\La_2-m_1 -2m_2 - l + 1))}
{\theta( -2\eta (\La_1 - l + 1))}\,
\notag
\\
\times\,\frac{\theta (z_1-z_2 +\eta\La_1-\eta \La_2 +2\eta (l-1))}
{\theta (z_1-z_2 -\eta\La_1-\eta \La_2 +2\eta (l-1))}\,]\,.
\notag
\eean
\end{lemma}

The Lemmas have important corollaries. Namely, consider the functional space
$F^m_{a_1,a_2}(z_1,z_2,\la)$ associated with $\vec\La$.
Consider the triangular matrices $A_{m_1,m_2}^{l_1,l_2}=
\om_{m_1,m_2}(T_{l_1,l_2},\la,z_1,z_2)$
and $\tilde A_{m_1,m_2}^{l_1,l_2}=
\tilde\om_{m_1,m_2}(T_{l_1,l_2},\la,z_1,z_2)$ for $(l_1,l_2),
(m_1,m_2)\in \Z^2_m$. Assume that $z_1, z_2, \Lambda_1, \Lambda_2,
\eta $ are such that all elements of these matrices are well defined.
(Notice that these conditions can be easily written explicitly.)

\begin{corollary}\label{basis}
The weight functions $\omega_M(t_1,...,t_m,\la,z,\tau),\,M\in \Z^2_m$ 
form a basis in $F^m_{a_1,a_2}(z_1,z_2,\la)$ unless one of the diagonal
elements \Ref{w.f.3} of the matrix $A$
 equals zero.
Similarly, the mirror weight functions 
$\tilde\omega_M(t_1,...,t_m,\la,z,\tau),\,M\in \Z^2_m$, 
form a basis in $F^m_{a_1,a_2}(z_1,z_2,\la)$ unless one of the diagonal
elements \Ref{w.f.4} of the matrix $\tilde A$ equals zero.
\end{corollary}

\begin{corollary}\label{A-B}
Let $B_{m_1,m_2}^{l_1,l_2}$ be the matrix inverse to $A$, ${}$
$\sum_{l_1,l_2}A_{m_1,m_2}^{l_1,l_2}B_{l_1,l_2}^{k_1,k_2}=\dl_{k_1,m_1}\dl_{k_2,m_2}$.
Consider the elliptic R-matrix $R_{\La_1,\La_2}(z_1-z_2,\la)$
 defined by \Ref{R.c}. Then
\bean\label{rm}
R^{kl}_{ij}\,=\,\sum_{r_1,r_2}\tilde A^{r_1r_2}_{kl}\, B^{ij}_{r_1r_2}.
\eean
In particular, the R-matrix is well defined unless one the diagonal elements
\Ref{w.f.3} of the matrix $A$ equals zero.
The determinant of the R-matrix
$R_{\La_1,\La_2}(z_1-z_2,\la)$ restricted to the weight space of weight
$\La_1+\La_2-2m$ is equal to the product of all of the
 diagonal elements \Ref{w.f.4}
of the matrix $\tilde A$ divided by the product of all of the
diagonal elements
\Ref{w.f.3} of the matrix $A$.
\end{corollary}

\subsection{The Shapovalov operator}\label{shap}
Let $Q^{\Lambda_j}(\la,\tau)$ be the diagonal operator on $V_{\Lambda_j}$
with diagonal matrix elements 
\bean
Q_k^{\Lambda_j}(\la,\tau)=
\left(\frac{\theta'(0,\tau)}{\theta(2\eta,\tau)}\right)^{k}
\prod_{l=1}^k
\frac
{\theta(2\eta(\Lambda_j+1-l),\tau)\theta(2\eta l,\tau)}
{\theta(\la+2\eta(\Lambda_j+1-k-l),\tau)\theta(\la-2\eta l,\tau)}\, .
\notag
\eean
Set
\bean
Q(\la,\tau)=Q^{\Lambda_1}(\la,\tau) 
\otimes 
Q^{\Lambda_2}(\la+2\eta h^{(1)},\tau) 
\otimes\cdots\otimes 
Q^{\Lambda_n}(\la+2\eta\sum_{j=1}^{n-1}h^{(j)},\tau).
\eean
$Q$ defines a diagonal operator on the tensor product $V_{\vec\La}$
called {\it the Shapovalov operator}. The Shapovalov operator plays an
important role in the study of the qKZB equations \cite{FV3}.

Let $\vec\La\in \C^n$ be such that $\La_1+...+\La_n=2m$ for some natural number
$m$. We consider the Shapovalov operator on the zero weight subspace
$V_{\vec\La}[0]$ and in particular its poles as a function of $\la$.

Let $e_M=e_{m_1}\otimes ....\otimes e_{m_n} 
\in V_{\vec\La}[0]$ be a basis vector, where
$M=(m_1,...,m_n), \, m_1+...+m_n=m$, is a vector of non-negative
integers.
The corresponding diagonal matrix element of the operator $Q$ is
\bean
Q_M(\la,\tau)\,=\,\left(\frac{\theta'(0,\tau)}{\theta(2\eta,\tau)}\right)^{m}
\prod_{j=1}^n
\prod_{l=1}^{m_j}
\frac
{\theta(2\eta(\Lambda_j+1-l),\tau)\theta(2\eta l,\tau)}
{\theta(\la+2\eta(m_j-l+1)+2\eta\sum_{k=1}^{j}(\La_k-2m_k),
\tau)}\, 
\notag
\\
\times \,
\frac
{1}
{\theta(\la-2\eta l+2\eta\sum_{k=1}^{j-1}(\La_k-2m_k),\tau)}\, .
\eean
For generic $\vec\La$, this coefficient has $2m$ simple pairwise distinct
poles. It is convenient to parametrize these poles as follows.

For $j=1,...,n-1,\, 0\leq l \leq m_j+m_{j+1},\, l\neq m_j,$ set
\bean\label{l}
\la_{M,j,l}\,=\, -\,2\eta (\La_j - m_j - l\,+\,\sum_{k=1}^{j-1}(\La_k-2m_k)),
\eean
and for $0\leq l \leq m_n+m_1,\, l\neq m_n,$ set
\bean\label{0}
\la_{M,n,l}\,=\, -\,2\eta (m_n - l).
\eean
Denote this set of $2m$ numbers $S_M$.

\begin{lemma}
The set $S_M$ is the set of poles of the coefficient $Q_M$
considered as a function of $\lambda$. $\square$

\end{lemma}

Any $\la_{M,j,l}$ defines  a new index $L$ as follows. If $j<n$, then
$L=(m_1,...,m_{j-1},l, m_j+m_{j+1}-l, m_{j+2},...,m_n)$. If $j=n$, then
$L=(m_1+m_n -l,m_2,...,m_{n-1},l)$. We call this index {\it dual to $M$
with respect to $\la_{M,j,l}$.} 

Let  $L$ be dual to $M$ with respect to $\la_{M,j,l}$. If $j<n$, 
then the set $S_L$ contains the number
$\la_{M,j,l}=\la_{L,j,m_j}$ and moreover
$M$ is dual to $L$ with respect to $\la_{M,j,l}$. 
If  $j=n$, then the set $S_L$ contains $-\la_{M,n,l}=\la_{L,n,m_n}$ and moreover
$M$ is dual to $L$ with respect to $-\la_{M,n,l}$. 

This construction defines a partition of the set of the pairs $(M,\la)$,
$\la \in S_M$, into the dual pairs.

\subsection{The poles of the Shapovalov operator and weight functions}\label{wf}
Define elliptic numbers by 
$$
[k]\,=\,{\theta (2\eta k,\tau)\over \theta (2\eta ,\tau)}
$$
and elliptic factorials by
\bean
[k]! \,=\, [1]\,[2]\, ... \,[k].
\notag
\eean

Let $\vec\La \in \C^n$, $m$ a positive integer.
Consider the weight functions associated with the pair $\vec\La$, $m$ and defined
in \Ref{w.f}.

For $ j < n$, let $m_1,...,m_{j-1},k,m_{j+2},...,m_n$
be non-negative integers such that $m_1+...+m_{j-1}+k+m_{j+2}+...+m_n=m$. Let 
$a,b$ be integers such that $a\neq b,\, 0\leq a,b \leq k$. Let
$M=(m_1,...,m_{j-1},a,k-a,m_{j+2},...,m_n)$ and
$L=(m_1,...,m_{j-1},b,k-b,m_{j+2},...,m_n)$. 
Consider the number $\la_0=2\eta (\La_j - a-b\,+\,\sum_{l=1}^{j-1}(\La_l-2m_l))$.
Notice that $-\la_0$ is a pole of the coefficients $Q_M(\la),\,Q_L(\la)$
of the Shapovalov operator on $V_{\vec\La}$, and the pairs
$M,-\la_0$ and $L,-\la_0$ are dual in the sense of Sec. \ref{shap}.

\begin{thm}\label{w.f.rel}
Under these conditions, for the weight functions defined in \Ref{w.f}, we have
\bean\label{re}
[a]!\,[k-a]!\,
\omega_{M}(t_1,\dots,t_m,2\eta (\La_j - a-b\,+\,\sum_{l=1}^{j-1}(\La_l-2m_l)),z,\tau)\,=\,
\\
\,[b]!\,[k-b]!\,
\omega_{L}(t_1,\dots,t_m,2\eta (\La_j - a-b\,+\,\sum_{l=1}^{j-1}(\La_l-2m_l)),z,\tau)\,.
\notag
\eean
\end{thm}

It is easy to see that formula \Ref{re} for arbitrary $n\geq 2$
and $j<n$ follows from the next Proposition, cf. \Ref{w.f}.

\begin{proposition}\label{w.f.rel.2}
Let $\vec\La=(\La_1, \La_2)$. Let $m, a, b$ be non-negative integers
such that $a\neq b$, $0\leq a,b \leq m$.
Then for the weight functions corresponding to $V_{\vec\La}$ we have
\bean\label{}
[a]!\,[m-a]!\,
\omega_{a,m-a}(t_1,\dots,t_m,2\eta (\La_1 - a-b),z,\tau)\,=\,
\\
\,[b]!\,[m-b]!\,
\omega_{b,m-b}(t_1,\dots,t_m,2\eta (\La_1 - a-b),z,\tau)\,.
\notag
\eean
\end{proposition}

\subsection{Proof of Proposition \ref{w.f.rel.2}}.

According to Lemmas \ref{triang}-\ref{w.f.mir},
it is easy to see that in order to prove the Proposition
it is enough to check that 
for any $N=(c,m-c)\in \Z^2_m$,
\bean\label{rel2}
[a]!\,[m-a]!\,
\omega_{a,m-a}(T_N,2\eta (\La_1 - a-b),z,\tau)\,=\,
\\
\,[b]!\,[m-b]!\,
\omega_{b,m-b}(T_N,2\eta (\La_1 - a-b),z,\tau)\,.
\notag
\eean

Let $a>b$. If $c>a$, then
\bean\label{rel3}
 \omega_{a,m-a}(T_N,2\eta (\La_1 - a-b),z,\tau)\,=\,
 \omega_{b,m-b}(T_N,2\eta (\La_1 - a-b),z,\tau)\,=\,0
\notag
\eean
by Lemma \ref{triang}. If $a\geq c> b$, then 
$\omega_{b,m-b}(T_N,2\eta (\La_1 - a-b),z,\tau)\,=\,0$
by Lemma \ref{triang}.

\begin{lemma}\label{zero}
If $a\geq c > b$, then 
$\omega_{a,m-a}(T_N,2\eta (\La_1 - a-b),z,\tau)\,=\,0$.
\end{lemma}
\begin{proof}
According to \Ref{w.f.2}, $\omega_{a,m-a}(T_N,2\eta (\La_1 - a-b),z,\tau)$
is the sum over all pairs $I_1, I_2$ of disjoint subsets of 
$\{1,\dots,m\}$ such that $I_1$ has $a$ elements, and $I_2$ has $m-a$ elements.
By Lemma \ref{diag}, all terms in this sum but the distinguished one
are equal to zero. The distinguished summand is zero since
the factor 
$\theta(\lambda+t_{c-b}-z_1-\eta\La_1+2\eta a)=
\theta(2\eta(\La_1-a-b)+z_1-\eta\La_1+2\eta b-z_1-\eta\La_1+2\eta a)$
is zero. 
\end{proof}

\begin{lemma}\label{disting}
If $a>b\geq c$, then equation \Ref{rel2} holds.
\end{lemma}
\begin{proof}
By Lemma \ref{diag}, it is enough to show that the distinguished
summands $D_1$ and $D_2$  of respectively
$[a]!\,[m-a]!\,
\omega_{a,m-a}(T_N,2\eta (\La_1 - a-b),z,\tau)$ and
$[b]!\,[m-b]!\,
 \omega_{b,m-b}(T_N,2\eta (\La_1 - a-b),z,\tau)$ are equal.
We have
\bean
D_1\,=\,
\prod_{j=1}^a
{ \theta (2\eta \,j)\over \theta (2\eta)}\,
\prod_{j=1}^{m-a}
{ \theta (2\eta \,j)\over \theta (2\eta)}\,
\prod_{j=1}^c
{ \theta (2\eta )\over \theta (2\eta\,j))}\,
\notag
\\
\times\,\prod_{i=1}^c \prod_{j=c+1}^a
{ \theta (z_1-z_2-\eta\La_1+\eta\La_2+2\eta(c-m+j-i))
\over
 \theta (z_1-z_2-\eta\La_1+\eta\La_2+2\eta(c-m+j-i+1))}\,
\prod_{j=1}^{a-c}
{ \theta (2\eta )\over \theta (2\eta\,j))}\,
\prod_{j=1}^{m-a}
{ \theta (2\eta )\over \theta (2\eta\,j))}\,
\notag
\\
\times\,\prod_{i=1}^{c}
{ \theta (2\eta (c-i-b) )\over \theta (-2\eta\La_1+2\eta (c-i))}\,
\prod_{i=c+1}^a
{ \theta (-z_1+z_2+\eta\La_1-\eta\La_2+2\eta(m-i-b))
\over
\theta (-z_1+z_2-\eta\La_1-\eta\La_2+2\eta(m-i))}\,
\notag
\\
\times\,
\prod_{i=a+1}^m
{ \theta (-2\eta\La_2+2\eta(2m-i-b))
\over
\theta (-2\eta\La_2+2\eta(m-i))}\,
{ \theta (-z_1+z_2+\eta\La_1-\eta\La_2+2\eta(m-i))
\over
\theta (-z_1+z_2-\eta\La_1-\eta\La_2+2\eta(m-i))}\,,
\notag
\eean
\bean
D_2\,=\,
\prod_{j=1}^b
{ \theta (2\eta \,j)\over \theta (2\eta)}\,
\prod_{j=1}^{m-b}
{ \theta (2\eta \,j)\over \theta (2\eta)}\,
\prod_{j=1}^c
{ \theta (2\eta )\over \theta (2\eta\,j))}\,
\notag
\\
\times\,\prod_{i=1}^c \prod_{j=c+1}^b
{ \theta (z_1-z_2-\eta\La_1+\eta\La_2+2\eta(c-m+j-i))
\over
 \theta (z_1-z_2-\eta\La_1+\eta\La_2+2\eta(c-m+j-i+1))}\,
\prod_{j=1}^{b-c}
{ \theta (2\eta )\over \theta (2\eta\,j))}\,
\prod_{j=1}^{m-b}
{ \theta (2\eta )\over \theta (2\eta\,j))}\,
\notag
\\
\times\,\prod_{i=1}^{c}
{ \theta (2\eta (c-i-a) )\over \theta (-2\eta\La_1+2\eta (c-i))}\,
\prod_{i=c+1}^b
{ \theta (-z_1+z_2+\eta\La_1-\eta\La_2+2\eta(m-i-a))
\over
\theta (-z_1+z_2-\eta\La_1-\eta\La_2+2\eta(m-i))}\,
\notag
\\
\times\,
\prod_{i=b+1}^m
{ \theta (-2\eta\La_2+2\eta(2m-i-a))
\over
\theta (-2\eta\La_2+2\eta(m-i))}\,
{ \theta (-z_1+z_2+\eta\La_1-\eta\La_2+2\eta(m-i))
\over
\theta (-z_1+z_2-\eta\La_1-\eta\La_2+2\eta(m-i))}\,.
\notag
\eean
$D_1$ and $D_2$ contain factors of five types:
$\theta (-2\eta\La_1 +2\eta$ $ \text{integer})$\,,
$\theta (-2\eta\La_2 +2\eta $ $\text{integer})$\,,
$\theta (-z_1+z_2+\eta\La_1-\eta\La_1+2\eta$ $ \text{integer})$\,,
$\theta (-z_1+z_2-\eta\La_1-\eta\La_1+2\eta $ $\text{integer})$\,,
$\theta (2\eta $ $\text{integer})$. Comparing the factors of each type
in $D_1$ and $D_2$ we conclude that  $D_1=D_2$. 

Lemma \ref{disting},
Proposition \ref{w.f.rel.2} and Theorem \ref{w.f.rel}  are proved.
\end{proof}

\begin{thm}\label{la-dep}
Let $r,s$ be integers. Then under conditions of Theorem \ref{w.f.rel},
for the weight functions defined in \Ref{w.f}, we have
\bean\label{la-re}
[a]!\,[k-a]!\,
e^{2\pi i\,s\,a\, (z_{j+1}-z_j +\eta\La_{j+1}+\eta\La_{j})}\,
\\
\times\,
\omega_{M}(t_1,\dots,t_m,\,r+s\tau
\,+\,2\eta (\La_j - a-b\,+\,\sum_{l=1}^{j-1}(\La_l-2m_l)),z,\tau)\,=\,
\notag
\\
\,[b]!\,[k-b]!\,
e^{2\pi i\,s\,b\, (z_{j+1}-z_j +\eta\La_{j+1}+\eta\La_{j})}\,
\notag
\\
\times\,
\omega_{L}(t_1,\dots,t_m,\,r+s\tau
\,+\,2\eta (\La_j - a-b\,+\,\sum_{l=1}^{j-1}(\La_l-2m_l)),z,\tau)\,.
\notag
\eean
\end{thm}

The Theorem follows from Theorem \ref{w.f.rel}, the explicit
formula for the weight functions and the transformation properties
of the theta function $\theta (\la,\tau)$ under the shifts
of $\la$ by $1$ and $\tau$.

\section{Properties of R-matrices}\label{section-4}

\subsection{The $\la$-poles of the R-matrix}\label{4.1}
According to Corollary \ref{A-B}, the poles of the R-matrix
\newline
$R_{\La_1,\La_2}(z,\la,\tau)$ considered as a function of $z, \la, \tau$
come from zeros of the functions $\theta(\la  - 2\eta (\La_1 + k))$,
$\theta(\la  - 2\eta (\La_1+ \La_2 +k))$,
$\theta(z - \eta (\La_1+ \La_2 +2k))$,
$\theta(-z  + \eta (- \La_1+ \La_2 +2k))$ where $k$ is an integer.
Thus the poles could be divided into the $\la$-poles and $z$-poles.
The following Theorem describes the $\la$-poles of the R-matrix.

\begin{thm}\label{La-poles}
Let $\La_1, \La_2, \eta, z$ be generic.
Consider the R-matrix $R_{\La_1,\La_2}(z,\la,\tau)$ restricted to the
weight space $(V_{\La_1}\otimes V_{\La_2})_{\La_1+\La_2-2m}$
of weight $\La_1+\La_2 - 2m$. Then
${}$
\begin{enumerate}
\item[I.] All $\la$-poles of the R-matrix have the form
$\la= 2\eta (\La_1 - k) + r + s\tau$ where $k= 1,..., 2m-1$ and $r,s\in \Z$.
\item[II.] All $\la$-poles of the R-matrix are simple.
\item[III.] Let $\la_0= 2\eta (\La_1 - k) + r + s\tau$ be a $\la$-pole
of the R-matrix, where $0<k<2m$ and $r,s\in \Z$. 
Let $K(z,\tau)\in \End ((V_{\La_1}\otimes V_{\La_2})_{\La_1+\La_2-2m})$
be the residue of the R-matrix at $\la_0$. Then the kernel of $K$
is the subspace in $(V_{\La_1}\otimes V_{\La_2})_{\La_1+\La_2-2m}$ consisting of vectors
$\sum_{i+j=m} u_{i,j} e_i\otimes e_j$ satisfying the relations
\bean\label{kernel}
[a]!\,[m-a]!\,
e^{2\pi i\,s\,a\, (-z +\eta\La_{1}+\eta\La_{2})}\, u_{a,m-a}\,=\,
\notag
\\
\,[b]!\,[m-b]!\,
e^{2\pi i\,s\,b\, (-z +\eta\La_{1}+\eta\La_{2})}\,u_{b,m-b}\,
\notag
\eean
where $a,b$ run through the set of all pairs satisfying $a+b=k$.
\end{enumerate}
\end{thm}
The Theorem will be proved in Sections \ref{part-I}-\ref{part-III}.
More information about the $\la$-poles of the R-matrices see in Section \ref{more}

\subsection{Proof of Part I of Theorem \ref{La-poles}}\label{part-I}
For any $M = (m_1,m_2)\in \Z^2_{m}$ and 
any meromorphic function $f(t_1,\dots,t_m)$ in $m$ variables
we let $\res_{M}f$ be the complex number
\bea
\res_{M}f =
\res_{t_1=z_1+\eta \La_1-2\eta(m_1-1)}
\cdots
\res_{t_{m_1-1}=z_1+\eta\La_1-2\eta}
\res_{t_{m_1}=z_1+\eta\La_1}
\\
\res_{t_{m_1+1}=z_2+\eta\La_2-2\eta(m_2-1)}
\cdots
\res_{t_{m_1+m_2-1}=z_2+\eta\La_2-2\eta}
\res_{t_{m_1+m_2}=z_2+a\eta\La_2}.
\eea
According to Proposition 30 in \cite{FTV1}, the matrices $A=\{\res_{L}\omega_M\}$,
$B=\{\res_{L}\tilde\omega_M\}$ are triangular. Namely, $\res_{L}\omega_M=0$ if $L>M$
and $\res_{L}\tilde\omega_M=0$ if $L<M$. Moreover, Proposition 30 gives explicit 
formulas for the diagonal elements of these matrices as alternating products of
theta functions. The formulas show that the diagonal elements do not have
factors of the form $\theta(\la  - 2\eta (\La_1+ \La_2 +k))$. 
Since $R=A^{-1}B$, this shows that
the R-matrix does not have $\la$-poles of the form $\la  = 2\eta (\La_1+ \La_2 +k)+r+s\tau,
k,r,s \in \Z$.

Now it follows from Corollary \ref{A-B} and formulas \Ref{w.f.3} and \Ref{w.f.4} that
the $\la$-poles of the R-matrix have the form $\la  = 2\eta (\La_1- k)$ where
$0< k <2m$ and $r,s \in \Z$.

\subsection{The poles of the Shapovalov operator and the R-matrix}
The following relation of the R-matrix with the Shapovalov operator
is proved in \cite{FV3}.

\begin{lemma}\label{Q-R} \cite{FV3}
For all $j,k,r,s$,  we have
\bean\label{q-r}
\lefteqn{
Q_j^{\Lambda_1}(\la+2\eta(\Lambda_2-2k),\tau)
Q_k^{\Lambda_2}(\la,\tau)
R_{rs}^{jk}(-\la,z,\tau) }
 \\
 & &=
Q_r^{\Lambda_1}(\la,\tau)
Q_s^{\Lambda_2}(\la+2\eta(\Lambda_1-2 r),\tau)
R_{jk}^{rs}(\la+2\eta(\Lambda_1+\Lambda_2-2(r+s)),z\tau).
\notag
\eean
\end{lemma}

\begin{lemma}\label{simple-poles}
If $\Lambda_1,\Lambda_2,\eta$ are generic, then
$R_{\Lambda_1,\Lambda_2}(z,\lambda,\tau)$, as a function of $\lambda$,
has only simple poles.
\end{lemma}

\begin{proof}
According to Part I of Theorem \ref{La-poles}, 
the poles of the R-matrix are of the form $\lambda=2\eta(\Lambda_1+
\mathrm{integer})$ (modulo $\Z+\tau\Z$).  

Consider formula \Ref{q-r}. Modulo $\Z+\tau\Z+2\eta\Z$, the poles of the R-matrix
on the left are at $-2\eta\Lambda_1$, 
whereas the poles of the R-matrix on the right are at
$-2\eta\Lambda_2$, i.e., elsewhere. Thus the poles of the
R-matrix on the right must come from the poles of the $Q$'s on
the left. It is therefore sufficient to show that the product
of $Q$'s on the left-hand side has only simple poles.

Now $Q_j^\Lambda(\la,\tau)$ has simple poles at $\la=2\eta l$, and $\la=
2\eta(-\Lambda+j+l-1)$, $l=1,\dots,j$. Hence the only way
\[
Q_j^{\Lambda_1}(\la+2\eta(\Lambda_2-2k),\tau)
Q_k^{\Lambda_2}(\la,\tau)
\]
can have higher order poles is that a pole of the first factor
at $\la=2\eta(-\Lambda_2+2k+l)$, $l=1,\dots,j$ coincides with 
a pole of the second factor at $\la=2\eta(-\Lambda_2+k+l'-1)$,
$l'=1,\dots,k$. But this never happens since $2k+l>k+l'-1$ for $l>0$ and
$l'\leq k$. So the R-matrix has only simple
poles.
\end{proof}

The $Q-R$-relation in \Ref{q-r} is powerful enough to tell
which matrix elements of the R-matrix could have which poles.

Part II of Theorem \ref{La-poles} is proved.

\subsection{Proof of Part III of Theorem \ref{La-poles}}\label{part-III}
According to Propositions 13 and 14 in \cite{FTV2}, we have
$R_{\Lambda_1,\Lambda_2}(z,\lambda + 1,\tau)=R_{\Lambda_1,\Lambda_2}(z,\lambda,\tau)$
and 
$$
R_{\Lambda_1,\Lambda_2}(z - w,\lambda + \tau,\tau)=
e^{\pi i (h^{(1)}(-z-\eta\La_2)+h^{(2)}(-w+\eta\La_1))}
R_{\Lambda_1,\Lambda_2}(z-w,\lambda,\tau)
e^{\pi i (h^{(1)}(z-\eta\La_2)+h^{(2)}(w+\eta\La_1))}.
$$
These two formulas show that the statement of Part III for arbitrary integers $r, s$
follows from the statement of Part III for $r=s=0$.

\begin{proposition}
Let $\la_0= 2\eta (\La_1 - k)$ be a $\la$-pole
of the R-matrix, where $0<k<2m$.
Let $K(z,\tau)\in \End ((V_{\La_1}\otimes V_{\La_2})_{\La_1+\La_2-2m})$
be the residue of the R-matrix at $\la_0$. Then the kernel of $K$
is the subspace in $(V_{\La_1}\otimes V_{\La_2})_{\La_1+\La_2-2m}$ consisting of vectors
$\sum_{i+j=m} u_{i,j} e_i\otimes e_j$ satisfying the relations
\bean\label{kernel}
[a]!\,[m-a]!\,u_{a,m-a}\, =\,\,[b]!\,[m-b]!\,u_{b,m-b}\,
\notag
\eean
where $a,b$ run through the set of all pairs satisfying $a+b=k$.
\end{proposition}
The Proposition follows from Lemmas \ref{triang} - \ref{w.f.mir},
Corollary \ref{A-B}, Proposition \ref{w.f.rel.2}, and the following Lemma from linear
algebra.
\begin{lemma}
Let $l,n$ be natural numbers, $A(\al)=\{A_{ij}(\al)\}$ be an $n\times n$ matrix depending on a parameter
$\al$. Assume that det $A(\al) = c\al^l+ O(\al^{l+1}),$ $c\neq 0$, as $\al \to 0$, and
$A_{ij}= O(\al)$ for all $j$ and $i=1,...,l$. Then the inverse matrix $A^{-1}(\al)$
has a simple pole at $\al=0$. Moreover, if $K$is the residue of $A^{-1}(\al)$ at $\al=0$,
then $K$ has rank $l$ and $K_{ij}=0$ for all $i$ and $j=l+1,...,n$. $\square$.

\end{lemma}

\subsection{Relations for matrix coefficients of R-matrices}

\begin{thm}\label{thm-coeff}
Let $r,s\in \Z$. Let $z, w, \La_1, \La_2, \eta$ be generic complex numbers.
Then 
\begin{enumerate}
\item[I.] For all $a,a',b, c, c',d$ we have
\bean\label{coeff-1}
e^{2\pi i\,s\,(b(-w+\eta\La_1)+d(z-\eta\La_2))}\,
{[b]!\over [d]!}\,
R_{\Lambda_1,\Lambda_2}(z-w, 2\eta (b'-b)+r+s\tau)^{a,b}_{d,c}\,=
\\
e^{2\pi i\,s\,(b'(-w+\eta\La_1)+d'(z-\eta\La_2))}\,
{[b']!\over [d']!}
R_{\Lambda_1,\Lambda_2}(z-w, 2\eta(b-b')+r+s\tau)^{a,b'}_{d',c}\,.
\notag
\eean
\item[II.]
For all $a,b,b',c,d,d'$ we have
\bean\label{coeff-2}
{}
\\
e^{2\pi i\,s\,(a(-z-\eta\La_2)+c(w+\eta\La_1))}\,
{[a]!\over [c]!}\,
R_{\Lambda_1,\Lambda_2}(z-w, 2\eta(\Lambda_1+\Lambda_2 -2 b-a-a')+r+s\tau)^{a,b}_{d,c}
\,=
\notag
\\
e^{2\pi i\,s\,(a'(-z-\eta\La_2)+c'(w+\eta\La_1))}\,
{[a']!\over [c']!}\,
R_{\Lambda_1,\Lambda_2}(z-w, 2\eta(\Lambda_1+\Lambda_2 -2 b-a-a')+r+s\tau)^{a',b}_{d,c'}
\,.
\notag
\eean
\end{enumerate}
\end{thm}

{\bf Example.} We have $R_{\La_1,\La_2}(z-w,\la)^{0,0}_{d,c}= \dl_{0,c}\dl_{0,d}$.
Then equation \Ref{coeff-1} says that 
$R_{\Lambda_1,\Lambda_2}(z-w, -2\eta k)^{0,k}_{k,0}=1$ and
$R_{\Lambda_1,\Lambda_2}(z-w, -2\eta k)^{0,k}_{i,j}=0$ 
for $(i,j)\neq (k,0)$.
Equation \Ref{coeff-2} says that
$R_{\Lambda_1,\Lambda_2}(z-w, 2\eta(\Lambda_1+\Lambda_2 - k))^{k,0}_{0,k}=1$
and $R_{\Lambda_1,\Lambda_2}(z-w, 2\eta(\Lambda_1+\Lambda_2 - k))^{k,0}_{i,j}=0$
for $(i,j)\neq (0,k)$.

We shall prove Part I of the Theorem, Part II is proved similarly.

According to the transformation properties of the R-matrix, described
in Propositions 13 and 14 in \cite{FTV2}, in order to prove Part I of
the Theorem it is enough to prove Part I for $r=s=0$,
thus it is enough to prove the following Proposition.

\begin{proposition}\label{Coeff}
  for all $a,a',b, c, c',d$ we have
\bean\label{coeff-1'}
{[b]!\over [d]!}\,
R_{\Lambda_1,\Lambda_2}(z-w, 2\eta (b'-b))^{a,b}_{d,c}\,=
\, {[b']!\over [d']!}
R_{\Lambda_1,\Lambda_2}(z-w, 2\eta(b-b'))^{a,b'}_{d',c}\,.
\notag
\eean
\end{proposition}
\begin{proof} For a natural number $m$, define
\bean
\Z^3_m =\{ (m_1,m_2,m_3)\in \Z^3_{\geq 0}\,|\, m_1+m_2+m_3=m\}.
\eean
Let $\omega^{\La_1,\La_2,\La_3}_{m_1,m_2,m_3}(t_1,...,t_m,\la,z_1,z_2,z_3,\tau)$, 
$(m_1,m_2,m_3)\in \Z^3_m$,
be the weight functions associated in Section \ref{bases} 
with the space $F_{a_1,a_2,a_3}^m(z_1,z_2,z_3)$ for $a_j=\eta\La_j, j=1,2,3$.

Define a $V_{\La_1}\otimes V_{\La_2}\otimes V_{\La_3}$-valued function
$\omega^{\La_1,\La_2,\La_3}(t,\la,z_1,z_2,z_3,\tau)$ by
\bean
\omega^{\La_1,\La_2,\La_3}(t_1&,&...,t_m , \la ,z_1,z_2,z_3,\tau)=
\notag
\\
&{}& \sum  _{(m_1,m_2,m_3)\in \Z^3_m} 
\omega^{\La_1,\La_2,\La_3}_{m_1,m_2,m_3}(t_1,...,t_m,\la,z_1,z_2,z_3,\tau)\,
e_{m_1}\otimes e_{m_2}\otimes e_{m_3}\,.
\notag
\eean
We have 
\bean
\omega^{\La_1,\La_3,\La_2}(t,\la,z_1,z_3,z_2,\tau)\,=\,P^{(23)}
R^{(23)}_{\La_2,\La_3}(z_2-z_3,\la-2\eta h^{(1)})\,
\omega^{\La_1,\La_2,\La_3}(t,\la,z_1,z_2,z_3,\tau)\,,
\notag
\eean
i.e.
\bean
\omega^{\La_1,\La_3,\La_2}_{m_1,l_3,l_2}(t,\la,z_1,z_3,z_2,\tau)\,&=&
\notag
\\
\sum_{m_2,m_3}R_{\La_2,\La_3}(z_2-z_3,&\la&-2\eta(\La_1-2m_1))^{l_2,l_3}_{m_2,m_3}
\omega^{\La_1,\La_2,\La_3}_{m_1,m_2,m_3}(t,\la,z_1,z_2,z_3,\tau)\,.
\notag
\eean
The coordinates of the functions $\omega^{\La_1,\La_2,\La_3}(t,\la,z_1,z_2,z_3,\tau)$ 
and $\omega^{\La_1,\La_3,\La_2}(t,\la,z_1,z_3,z_2,\tau)$ satisfy the resonance relations
of Theorem \ref{w.f.rel}. Let us write one of the relation for coordinates of the
function $\omega^{\La_1,\La_3,\La_2}(t,\la,z_1,z_3,z_2,\tau)$,
\bean
[m_1]!\,[l_3]!\,{}\, \omega^{\La_1,\La_3,\La_2}_{m_1,l_3, a}(t,\la_0,z_1,z_3,z_2,\tau)\,=
[\tilde m_1]!\,[\tilde l_3]!\,{}\, \omega^{\La_1,\La_3,\La_2}_{\tilde m_1,\tilde l_3, a}(t,\la_0,z_1,z_3,z_2,\tau)
\notag
\eean
where $\la_0=2\eta(\La_1-m_1-\tilde m_1)$ and $m_1+l_3= \tilde m_1 + \tilde l_3$.
We assume that $m_1 < \tilde m_1$.

We can write this relation as
\bean
[m_1]![l_3]!\, 
\sum_{m_2,m_3} R_{\La_2,\La_3} (z_2-z_3, 2\eta (m_1-\tilde m_1))^{a,l_3}_{m_2,m_3}
\omega^{\La_1,\La_2,\La_3}_{m_1,m_2, m_3}(t,\la_0,z_1,z_2,z_3,\tau)\,=
\notag
\\
 {} [{\tilde m_1}] ! [\tilde l_3]!\, 
\sum_{\tilde m_2,\tilde m_3}R_{\La_2,\La_3}(z_2-z_3,2\eta(\tilde m_1- m_1))
^{a,\tilde l_3}_{\tilde m_2,\tilde m_3}
\omega^{\La_1,\La_2,\La_3}_{\tilde m_1,\tilde m_2, \tilde m_3}(t,\la_0,z_1,z_2,z_3,\tau)\,.
\notag
\eean
Using the resonance relations for the function
$\omega^{\La_1,\La_2,\La_3}(t,\la,z_1,z_2,z_3,\tau)$, 
\bean
[m_1]!\,[m_2]!\,\omega^{\La_1,\La_2,\La_3}_{m_1,m_2, m_3}(t,\la_0,z_1,z_2,z_3,&\tau &)\,=
\notag
\\
{[ {\tilde m_1}]}! \, 
[\tilde m_2]!\, \omega^{\La_1,\La_2,\La_3}_{\tilde m_1,\tilde m_2, m_3}(&t&,\la_0,z_1,z_3,z_2,\tau)\,,
\notag
\eean
we can
rewrite the last equation as
\bean\label{equation}
\sum_{m_3}\,( \,{[l_3]!\over [a+l_3-m_3]!}
\, R_{\La_2,\La_3}(z_2-z_3,2\eta(m_1-\tilde m_1))^{a,l_3}_{a+l_3-m_3,m_3}
\, -
\\
{[ \tilde l_3]!\over [a+ \tilde l_3-m_3]!}
R_{\La_2,\La_3}(z_2-z_3,2\eta(\tilde m_1- m_1))
^{a,\tilde l_3}_{a+\tilde l_3-m_3,m_3}\,)\,
\notag
\\
\times \, 
\omega^{\La_1,\La_2,\La_3}_{m_1, a + l_3 - m_3, m_3}(t,\la_0,z_1,z_2,z_3,\tau)\,=\,0.
\notag
\eean
In this equation the R-matrix elements with negative indices are considered to be zero and
$1/[j]!$ is zero for negative $j$.

\begin{lemma} Let $\la_0=2\eta(\La_1-m_1-\tilde m_1)$ where $m_1, \tilde m_1$ are non-zero integers and
$m_1 < \tilde m_1$.
Consider the functions 
$\omega^{\La_1,\La_2,\La_3}_{m_1, a+ l_3 - m_3, m_3}(t,\la_0,z_1,z_2,z_3,\tau)$,
$m_3=0,..., a+l_3$, as functions of $t_1,...,t_m$. Then these functions 
are linearly independent over $\C$.
\end{lemma}
Indeed, for any $(k_1,k_2,k_3)\in \Z^3_m$ define
a point $T_{k_1,k_2,k_3}\in \C^m$ by
\bean
T_{k_1,k_2,k_3}=(z_1-\eta\La_1+2\eta (k_1-1), z_1-\eta\La_1+2\eta (k_1-2), ...,
z_1-\eta\La_1,
\notag
\\
z_2-\eta\La_2+2\eta (k_2-1), z_2-\eta\La_2+2\eta (k_2-2),
...., z_2-\eta\La_2,
\notag
\\
 z_3-\eta\La_3+2\eta (k_3-1), z_3-\eta\La_3+2\eta (k_3-2),
....,z_3-\eta\La_3).
\notag
\eean
It is easy to see that the matrix 
$\{A_{m_3,k_3}=
 \omega^{\La_1,\La_2,\La_3}_{m_1, a+ l_3 - m_3, m_3}(T_{m_1, a+l_3-k_3,k_3},\la_0,z_1,z_2,z_3,\tau)\}$,
with  $0\leq m_3, k_3 \leq a+l_3$, is triangular with nonzero diagonal elements, cf. Lemmas \ref{triang}, 
\ref{diag}.
This proves the Lemma.

By the Lemma, the coefficients of the weight functions in equation \Ref{equation}
must be zero. Thus the Proposition is proved.
\end{proof}


\section{Hypergeometric solutions}\label{hypergsol}

\subsection{ Phase functions}
We assume that $p$ has a positive imaginary part, and set $r=e^{2\pi ip}$,
$q=e^{2\pi i\tau}$. Then the convergent infinite product
\begin{equation}\label{phase1}
\Omega_a(t):=\Omega_a(t,\tau,p)=\prod_{j=0}^{\infty}\prod_{k=0}^{\infty}
\frac{(1-r^jq^ke^{2\pi i(t-a)})(1-r^{j+1}q^{k+1}e^{-2\pi i(t+a)})}
{(1-r^jq^ke^{2\pi i(t+a)})(1-r^{j+1}q^{k+1}e^{-2\pi i(t-a)})}\,,
\end{equation}
is called a (one variable) phase function.
It obeys the identities
\[
\Omega_a(z+p,\tau,p)
=e^{2\pi i a}\frac{\theta(z+a;\tau)}
{\theta(z-a;\tau)}
\Omega_a(z,\tau,p),
\]
\[
\Omega_a(z+1,\tau,p)=\Omega_a(z,\tau,p),
\]
\[
\Omega_a(z,\tau,p)
=\Omega_a(z,p,\tau),
\]
see \cite{FV3}. Notice that the phase function obeys many
other remarkable identities that lead to an $SL(3, \Z)$ symmetry,
see \cite{FV3, FV4}.

Fix $a=(a_1,...,a_n)\in \C^n$ and a natural number $m$. Define
the $m$-variable phase function by
\begin{equation}\label{phasenew}
\Om_a(t_1,\dots,t_m,z_1,\dots,z_n,\tau,p)=\prod_{j=1}^m\prod_{l=1}^n\Om_{a_l}(t_j-z_l)
\prod_{1\leq i<j\leq m}\Om_{-2\eta}(t_i-t_j).
\end{equation}

\subsection{The universal hypergeometric function, \cite{FTV2}}
Let  $\vec\La=(\Lambda_1,\ldots,\Lambda_n)$ be 
such that $\Lambda_1+\ldots+\Lambda_n=2m$ for 
some positive integer $m$, and set $a_i=\eta\Lambda_i$. 

Introduce a meromorphic $V_{\vec\La}[0]\otimes V_{\vec\La}[0]$-valued function $u$ of variables
$z, \vec \La \in \C^n,\, \la,\mu \in \C$ by
\bean\label{sol}
u(z,\la,\mu,\tau,p, \vec \La)\,=\,\sum_{I,J,\,|I|=|J|=m}\,
u_{IJ}(z,\la, \mu, \tau,p,\vec \La)\, e_I\otimes e_J
\eean
where
\bean\label{coord}
u_{IJ}(z,\la,\mu,\tau,p,\vec \La)\,=&
\\
e^{-\pi i{\mu\la\over2\eta}}\,
\int_{T^m}\,&\Om(t,z,a,\tau,p,\eta)\,\om_I(t,\la,z,a,\tau,\eta)\,
\tilde \om_J(t,\mu,z,a,p,\eta)\,dt\,.
\notag
\eean
Here, for every $M=(m_1,...,m_n)$, $e_M$ denotes $e_{m_1}\otimes ... \otimes e_{m_n}$,
${}$  $\Om$ is the $m$-variable phase function \Ref{phasenew}, ${}$ $\om_I(t,\la,z,a,\tau,\eta)$ and 
$\tilde \om_J(t,\mu,z,a,p,\eta)$ are the weight functions \Ref{w.f},\Ref{m.w.f}, ${}$
$dt=dt_1\wedge ... \wedge dt_m$. 
The integrand is $1$-periodic with respect to shifts of variables $t=(t_1,...,t_m)$ and
therefore defines a meromorphic function of $t$ on $\C^m/\Z^m$. $T^m=T^m(z,\tau, p,\vec\La)$ 
is an $m$-dimensional
cycle in $\C^m/\Z^m$ which is a suitable deformation of the torus
$\R^m/\Z^m \subset \C^m/\Z^m$ determined by $z, \tau, p,\vec\La $, see \cite {FTV2}.
Namely, the integral is defined by analytic continuation from the region where
Re $(\La_i)<0$, $z_i\in \R$ and $\Im (\eta) < 0$. In this region, the integration 
cycle is just the torus $\R^m/\Z^m $.

The function $u$ will be called {\it the universal hypergeometric
function associated with $V_{\vec\La}$}.

For a fixed generic $\vec\La$, the function $u$
is a meromorphic function of the remaining parameters and 
satisfies equations \Ref{eq1}, see explicit conditions on $\vec\La$ in
\cite{FTV2} and \cite {MV}.

{\bf Remark.} In \cite{FTV2} we used only weight functions
and no mirror weight functions. Then the qKZB equations \Ref{eq1}
only involve qKZB operators and no mirror qKZB operators. The choices
of this paper and \cite{FV3, FV4} make the qKZB heat equation of
\cite{FV3, FV4} more transparent. The proof that $u$ obeys the relations
\Ref{eq1} is the same as the proof of Theorem 31 in \cite{FTV2}.
Note however that the conventions in the definition of $D_j$ are different
there.

\subsection{Finite dimensional representations, \cite{MV}}\label{finite}
Assume that $\La^0_1,...,\La^0_n$ are natural numbers, and for some natural number $m$,
$\La^0_1+...+\La_n^0=2m$. Set $\vec\La^0=(\La_1^0,...,\La_n^0)$.

Let $M=(m_1,...,m_n)$ be a vector of non-negative integers. Say that $M$ is {$\vec\La^0$-
admissible}, if $m_j\leq \La_j$ for all $j$.

\begin{thm}\label{f.dim}\cite{MV}
The coordinates $u_{IJ}$ 
of the universal hypergeometric function defined by \Ref{sol}
are meromorphic functions of its variables.
If at least one of the indices $I, \,J$ is $\vec\La^0$-admissible,
then $u_{IJ}|_{\vec\La=\vec\La^0}$,defined by analytic continuation,
 is a well defined meromorphic function of the remaining
variables. Thus, for every fixed $\vec\La^0$-admissible index $M$, 
\bean
\sum_{I,\,|I|=m}u_{IM}|_{\vec\La=\vec\La^0}\,e_I\otimes e_M
\qquad
\text{and}
\qquad
\sum_{I,\,|I|=m}u_{MI}|_{\vec\La=\vec\La^0}\,e_M\otimes e_I,
\eean
are well defined meromorphic functions. Moreover, the first of the functions
satisfies the first 
and the third system of equations in \Ref{eq1}, and the second of the functions
satisfies the second and the third system of equations in \Ref{eq1}.
\end{thm}

Consider the tensor product $L_{\vec\La^0}=L_{\vec\La_1^0}\otimes ...
\otimes L_{\vec\La^0_n}$ of the corresponding finite dimensional vector spaces and the
$V_{\vec\La^0}[0]\otimes L_{\vec\La^0}[0]$-valued function
\bean
\tilde u(z,\la,\mu,\tau,p, \vec \La^0)\,=\,
\sum_{ I,\,{}|I|=m,   \atop \text{adm} \,J,\,{}|J|=m}\,
u_{IJ}(z,\la, \mu, \tau,p,\vec \La^0)\, e_I\otimes e_J.
\notag
\eean
The function 
$$
\tilde v=(1\otimes\prod_{j=1}^n D_j^{-z_j/p})\,\tilde u
$$
obeys the qKZB
equations in the first factor \Ref{first-factor}.
It means that for every complex $\mu$ and every linear form on
$L_{\vec\La} [0]$, we have a solution of the qKZB equations. 
We call such solutions {\it the elementary hypergeometric solutions}
with values in $V_{\vec\La^0}[0]$.
A solution $w(z,\la)$ of the qKZB equations with values in $V_{\vec\La^0}[0]$ 
is called {\it a hypergeometric solution}
if it can be represented in the form
\Ref{hypergeom}. 

Define {\it the universal hypergeometric function associated with } $L_{\vec\La^0}$
by
\bean\label{admsol}
u(z,\la,\mu,\tau,p, \vec \La^0)\,=\,\sum_{\text{adm} \,I,J,\,{}|I|=|J|=m}\,
u_{IJ}(z,\la, \mu, \tau,p,\vec \La^0)\, e_I\otimes e_J
\eean
where the sum is over all $\vec\La^0$-admissible indices $I,J$.
This is a $L_{\vec\La^0}[0]\otimes L_{\vec\La^0}[0]$-valued function.
Let $\pi: V_{\vec\La^0}[0]\to L_{\vec\La^0}[0]$ be the canonical projection.
Then $u\,=\, \pi \otimes 1\,{}\,\tilde u$.

It follows from the properties of R-matrices formulated in Sec. \ref{geom}
that the universal hypergeometric function 
satisfies equations \Ref{eq1}. Define {\it the fundamental hypergeometric
solution $v$ associated with} $L_{\vec\La^0}$ by formula \Ref{fund}.
The $L_{\vec\La^0}[0]\otimes L_{\vec\La^0}[0]$-valued function $v$ obeys the qKZB
equations in the first factor \Ref{first-factor}. Thus
for every complex $\mu$ and every linear form on
$L_{\vec\La^0} [0]$, we have a solution of the qKZB equations. 
We call such solutions {\it the elementary hypergeometric solutions} with values
in $L_{\vec\La^0}[0]$. 
A solution $w(z,\la)$ of the qKZB equations with values in $L_{\vec\La^0}[0]$ 
is called {\it a hypergeometric solution}
if it can be represented in the form
\Ref{hypergeom}. 

The second system of equations in \Ref{eq1} gives the
monodromy of these solutions, see  \cite{MV}.

\section {Resonance relations for solutions with values in Verma modules}\label{Verma}

\subsection{Resonance relations}
Let $u=\sum u_{IJ}e_I\otimes e_J$ be the universal hypergeometric function
associated with $V_{\vec\La}$, $\La_1+...+\La_n=2m$.

\begin{thm}\label{v} Let $r,s$ be integers.
${}$

\begin{enumerate}
\item[I.]
Let $n>1$. For $ j < n$, let $m_1,...,m_{j-1},k,m_{j+2},...,m_n$
be non-negative integers such that $m_1+...+m_{j-1}+k+m_{j+2}+...+m_n=m$. Let 
$a,b$ be integers such that $a\neq b,\, 0\leq a,b \leq k$. Let
$M=(m_1,...,m_{j-1},a,k-a,m_{j+2},...,m_n)$ and
$L=(m_1,...,m_{j-1},b,k-b,m_{j+2},...,m_n)$. Then for any $J$
we have
\bean\label{v1}
[a]!\,[k-a]!\,
e^{2\pi i\,s\,a\, (z_{j+1}-z_j +\eta\La_{j+1}+\eta\La_{j})}\,
\\
\times\,
u_{MJ}(z,\,r+s\tau\,+\,2\eta (\La_j - a-b\,+\,\sum_{l=1}^{j-1}(\La_l-2m_l)),\mu,\tau,p,\vec\La)=
\notag
\\
\,[b]!\,[k-b]!\,
e^{2\pi i\,s\,b\, (z_{j+1}-z_j +\eta\La_{j+1}+\eta\La_{j})}\,
\notag
\\
\times\,
u_{LJ}(z,\,r+s\tau\,+\,
2\eta (\La_j - a-b\,+\,\sum_{l=1}^{j-1}(\La_l-2m_l)),\mu,\tau,p,\vec\La).
\notag
\eean
\item[II.] 
Let $n>1$. Let  $m_2,...,m_{n-1},k$
be non-negative integers such that $m_2+...+m_{n-1}+k=m$. Let 
$a,b$ be such that $a\neq b,\, 0\leq a,b \leq k$. Let
$M=(k-a,m_2,...,m_{n-1},a)$ and
$L=(k-b,m_2,...,m_{n-1},b)$. Then for any $J$
we have
\bean\label{v2}
[a]!\,[k-a]!\,
e^{2\pi i\,s\,a\, (z_{1}-z_n +\eta\La_{1}+\eta\La_n-p)}\,
\\
\times\,
u_{MJ}(z,\,r+s\tau\,+\,2\eta (a-b),\mu,\tau,p,\vec\La)\,=\,
\notag
\\
\,[b]!\,[k-b]!\,
e^{2\pi i\,s\,b\, (z_{1}-z_n +\eta\La_{1}+\eta\La_n-p)}\,
\notag
\\
\times\,
u_{LJ}(z,\,r+s\tau\,+\,2\eta (b-a),\mu,\tau,p,\vec\La).
\notag
\eean
\item[III.] Let $n=1$, $\vec \Lambda = 2m$. In this case
$u=u_{mm} e_m\otimes e_m$ and $u$ does not depend on $z$. We claim that
\bean\label{n1}
{}\,{}\, e^{2\pi i\,s\,a\, (4\eta m-p)}\,
u_{mm}(\,r+s\tau\,+\,2\eta a,\mu,\tau,p,\vec\La)\,=\,
u_{mm}(\,r+s\tau\,-\,2\eta a,\mu,\tau,p,\vec\La)\,
\eean
for $a = 1,...,m$.
\end{enumerate}
\end{thm}

\begin{corollary}\label{part-1}
Any hypergeometric solution $w(z,\la)=\sum_J \,w_J(z,\la)\, e_J$ 
defined in \Ref{hypergeom} obeys
the resonance relations
\bean\label{res-hyp-1}
[a]!\,[k-a]!\,
e^{2\pi i\,s\,a\, (z_{j+1}-z_j +\eta\La_{j+1}+\eta\La_{j})}\,
\\
\times\,
w_{M}(z,\,r+s\tau\,+\,2\eta (\La_j - a-b\,+\,\sum_{l=1}^{j-1}(\La_l-2m_l)))=
\notag
\\
\,[b]!\,[k-b]!\,
e^{2\pi i\,s\,b\, (z_{j+1}-z_j +\eta\La_{j+1}+\eta\La_{j})}\,
\notag
\\
\times\,
w_{L}(z,\,r+s\tau\,+\,2\eta (\La_j - a-b\,+\,\sum_{l=1}^{j-1}(\La_l-2m_l)))
\notag
\eean
for  $j, a, b, k, L, M$ defined in Part I of Theorem \ref{v} and the resonance
relations
\bean\label{}
[a]!\,[k-a]!\,
e^{2\pi i\,s\,a\, (z_{1}-z_n +\eta\La_{1}+\eta\La_n-p)}\,
w_{M}(z,\,r+s\tau\,+\,2\eta (a-b))\,=\,
\\
{}\,[b]!\,[k-b]!\,
e^{2\pi i\,s\,b\, (z_{1}-z_n +\eta\La_{1}+\eta\La_n-p)}\,
w_{L}(z,\,r+s\tau\,+\,2\eta (b-a))
\notag
\eean
for $a, b, k, L, M$ defined in Part II of Theorem \ref{v}.
\end{corollary}

Formula \Ref{v1} of Theorem \ref{v} follows from Theorem \ref{la-dep}
and formula \Ref{coord}. 
Hence any hypergeometric solution $w=\sum_J\, w_J\, e_J$ 
defined in \Ref{hypergeom} obeys relations \Ref{res-hyp-1} for
$j, a, b, k, L, M$ defined in Part I of Theorem \ref{v}.

We  deduce \Ref{v2} from \Ref{res-hyp-1} in Sec. \ref{proof}.
Part III is proved in Sec. \ref{proof-iii}.

\subsection{ The monodromy with respect to permutations of variables}
In order to prove Part II of Theorem \ref{v} we recall some facts about monodromy
properties of the qKZB equations.

\begin{thm}\label{Perm}\cite{FTV2}
For $\vec\La\in \C^n$ such that $\La_1+...+\La_n =2m$, 
let $v (z_1,\ldots,z_n,\la)$ be a solution of 
the qKZB equations with values
in $V_{\vec\La}[0]=
(V_{\La_1}\otimes \ldots\otimes V_{\La_n})[0]$, step $p$ and modulus $\tau$.
Then for any $j=1,\ldots,n-1$, the function
\bean\label{perm}
P^{(j,j+1)}\,R_{\La_j, \La_{j+1}}^{(j,j+1)}
(z_{j+1}-z_j,\la-2\eta \sum_{l=1}^{j-1} h^{(l)},\tau)\,
v(z_1,\ldots,z_{j+1},z_j,\ldots ,z_n,\la)
\notag
\eean
is a solution of the qKZB equations with values in 
$(V_{\La_1}\otimes\ldots
V_{\La_{j+1}}\otimes V_{\La_j}\ldots\otimes V_{\La_n})[0]$, step $p$ and
modulus $\tau$. Here $P^{(j,j+1)}$ is the permutation of the $j$-th and
$j+1$-th factors, and $R_{\La_j,\La_{j+1}}(z,\la,\tau)\in
\End(V_{\La_{j}}\otimes V_{\La_{j+1}})$ is the elliptic $R$-matrix
with modulus $\tau$.
\end{thm}

Let $v$ denote the fundamental hypergeometric solution associated
with $\vec\La$. Let $\vec\La^j$ denote the vector
$(\La_1,...,\La_{j+1},\La_j,...,\La_n)$
and $v^j$ denote the fundamental hypergeometric solution associated
with $\vec\La^j$. According to Theorem \ref{Perm}, the $V_{\vec\La^j}[0]\otimes
V_{\vec\La}[0]$-valued function
$$
u^j\,=\,P^{(j,j+1)}\,R_{\La_j, \La_{j+1}}^{(j,j+1)}
(z_{j+1}-z_j,\la-2\eta \sum_{l=1}^{j-1} h^{(l)},\tau)\otimes 1 \,{}\,{}
\text{}\,
v(z_1,\ldots,z_{j+1},z_j,\ldots ,z_n,\la,\mu,\tau,p)
\notag
$$
and the $V_{\vec\La^j}[0]\otimes V_{\vec\La^j}[0]$-valued function
$v^j(z,\la,\mu,\tau,p)$ satisfy the same qKZB equations in the first factor.

The next Theorem describes a relation between the two solutions and can be
considered as a description of the monodromy of the hypergeometric solutions
constructed in Sec. \ref{sol.qKZB} with respect to permutation of variables.

Introduce a new $R$-matrix $\tilde R_{A,B}(z,\mu,p)\in\End(V_A\otimes V_B)$ by
\bean\label{new-r}
\tilde R_{A,B} (z,\mu,p)\,={}
\\
e^{2\pi i ABz/p}\,\Bigl(\,\frac{\al(\mu)}{\al(\mu-2\eta h^{(2)})}\Bigr)^{z/p}\,
R_{A,B} (z,\mu, p)\,\Bigl(\,\frac{\al(\mu-2\eta(h^{(1)}+h^{(2)}))}
{\al(\mu -2\eta h^{(1)})}\Bigr)^{z/p}.
\notag
\eean

\begin{thm}\label{permut}\cite{FTV2}
\bean\label{perm-1}
u^j(z,\la,\mu,\tau,p)\,={}
\\
1\otimes P^{(j,j+1)}
\,(\tilde R_{\La_{j+1}, \La_{j}}^{(j,j+1)}
( z_{j}- z_{j+1}, \mu - 2\eta \sum_{l=j+1}^{n} h^{(l)},p)
\,)^{-1}\,{}\,v^j(z,\la,\mu,\tau,p)\,.
\notag
\eean
\end{thm}

{\bf Remark.} According to Proposition 12 in \cite{FTV2}, the matrix
$\tilde R_{A,B}(z,\mu,p)$ is $p$-periodic,
$$
\tilde R_{A,B}(z+p,\mu,p)\,=\,\tilde R_{A,B}(z,\mu,p)\,.
$$
Hence, formula \Ref{perm-1} expresses the solution $u^j$ as a linear combination
of solutions $v^j$ with $p$-periodic coefficients.

{\bf Remark.} Although we used in \cite{FTV2} only weight
functions and no mirror weight functions, the proof of Theorem \ref{permut}
is the same as the proof on Theorem 36 in \cite{FTV2}.

\subsection{Proof of Part II of Theorem \ref{v}}\label{proof}

For $\vec\La\in \C^n$ such that $\La_1+...+\La_n =2m$, 
let $v (z_1,\ldots,z_n,\la)$ be the fundamental hypergeometric solution
associated with $V_{\vec\La}$. Let $\vec\La^\vee=(\La_1^\vee,...,\La_n^\vee)$
be defined by $(\La_1^\vee,...,\La_n^\vee)=(\La_n,\La_1,\La_2,...
\La_{n-1})$. Let $\Dl:V_{\vec\La}\to V_{\vec\La^\vee}$ be the linear operator
defined by
\bean\label{delta}
\Dl\,=\, \Gamma_1 P^{(1,2)} P^{(2,3)} ... P^{(n-1,n)}\, .
\eean
Consider the $V_{\vec\La^\vee}[0]\otimes V_{\vec\La}[0]$-valued function
\bean
w(z_1,...,z_n,\la,\mu,\tau,p)\,=\,\Dl\,\otimes \,1\,{}\,
v(z_2,z_3,...,z_n,z_1,\la,\mu,\tau, p)\,.
\notag
\eean
Write $w$ in coordinate form $w=\sum \,w_{IJ}\,e_I\otimes e_J$
where $\{e_I\}$ is the standard basis in $V_{\vec\La^\vee}[0]$
and  $\{e_J\}$ is the standard basis in $V_{\vec\La}[0]$.

\begin{thm}\label{Transf}
Let $n>1$. Let $M=(a,k-a,m_3,...,m_n)$ and $L=(b,k-b,m_3,...,m_n)$.
Then
\bean\label{}
{}\,
\\
{}\,
[a]!\,[k-a]!\,
e^{2\pi i\,s\,a\, (z_2-z_1 +\eta\La_{1}+\eta\La_n-p)}\,
w_{MJ}(z,\,r +s\tau\,+\,
2\eta (\La_n - a-b),\mu,\tau,p)=
\notag
\\
\,[b]!\,[k-b]!\,
e^{2\pi i\,s\,b\, (z_2-z_1 +\eta\La_{1}+\eta\La_n-p)}\,
w_{LJ}(z,\,r +s\tau\,+\,2\eta (\La_n - a-b)
,\mu,\tau,p).
\notag
\eean
\end{thm}
\begin{corollary}\label{part-2}
For any $a, b, k, L, M$ defined in Part II of Theorem \ref{v}, the coordinates
of the hypergeometric solution $v\,=\,\sum\,v_{IJ}\,e_I\otimes e_J$
obey the resonance relation in \Ref{v2}.
\end{corollary}

This completes the proof of Part II of Theorem \ref{v}
since $u_{LJ}$ is obtained from $v_{LJ}$ by multiplying with a nonzero
factor independent of $L$ or $\lambda$, see \Ref{fund}.

{\bf Proof of the Corollary.} We have
\bean
w_{(a,k-a,m_3,...,m_n)\,J}(z_1,...,z_n,\,r+s\tau\,+\,2\eta(\La_n-a-b),\mu,\tau,p)=
\notag
\\
v_{(k-a,m_3,...,m_n,a)\,J}(z_2,...,z_n,z_1,\,r+s\tau\,+\,2\eta(a-b),\mu,\tau,p)
\notag
\eean
and similarly for $w_{LJ}(z,\,r+s\tau\,+\,2\eta (\La_n - a-b),\mu,\tau,p)$.
$\square$

{\bf Proof of Theorem \ref{Transf}.}
Apply to the function $v$ the transformation of Theorem \ref{Perm} for $j=n-1$,
then apply to the result
the transformation of Theorem \ref{Perm} for $j=n-2$, then repeatedly apply
the transformations for $j=n-3,...,1$.
The resulting $V_{\vec\La^\vee}[0]\otimes V_{\vec\La}[0]$-valued function
$w'(z,\la,\mu,\tau,p)$ has the form
\bean
w'(z_1,...,z_n,\la,\mu,\tau,p)\,=\,R^{(2,1)}_{\La_1,\La_n}(z_2-z_1,\la,\tau)\,
R^{(3,1)}_{\La_2,\La_n}(z_3-z_1,\la-2\eta h^{(2)},\tau)\,
...
\notag
\\
R^{(n,1)}_{\La_{n-1},\La_n}(z_n-z_1,\la-2\eta \sum_{l=2}^{n-2}h^{(l)},\tau)\,
P^{(1,2)} P^{(2,3)} ... P^{(n-1,n)}\, 
\,\otimes\, 1\,{}\, v(z_2,...,z_n,z_1, \la,\mu,\tau,p)\,.
\notag
\eean
By Theorem \ref{Perm} the function $w'$ satisfies the qKZB equations in the first
factor. In particular, $w'(z_1+p,z_2,...,z_n)= K_1(z)\otimes 1\,
w'(z_1,z_2,...,z_n)$. But
\bean
K_1(z)\otimes 1\,{}\,w'(z_1,z_2,...,z_n)=
\Dl\,\otimes \,1\,{}\,
v(z_2,z_3,...,z_n, z_1,\la,\mu,\tau, p)\,.
\notag
\eean
Write the function $w'$ in coordinate form, $w'\,=\,\sum\,w'_{IJ}\,e_I\otimes e_J$.
By Theorem \ref{permut}, for any $J$, the 
$V_{\vec\La^\vee}[0]$-valued function $\sum_J\,w'_{IJ}\,e_I$ is a hypergeometric
solution of the qKZB equations. Hence this function obeys the resonance relations
\Ref{res-hyp-1} of Corollary \ref{part-1}.
$\square$

\subsection{Proof of Part III of Theorem \ref{v}}\label{proof-iii}
The function $u_{m,m}$ has the form
\bean
u_{m,m}(\lambda,\mu,\tau,p)=e^{-\pi i{\mu\la\over2\eta}}\,
\int_{T^m}\,&\Om(t,\tau,p,\eta)\,\om_m(t,\la,\tau)\, \om_m(t,\mu,p)\,dt\,.
\notag
\eean
Here
\begin{equation}\label{}
\Om(t,\tau,p)=\prod_{j=1}^m\Om_{2\eta m}(t_j,\tau,p)
\prod_{1\leq i<j\leq m}\Om_{-2\eta}(t_i-t_j,\tau,p)
\notag
\end{equation}
and
\begin{equation}\label{}
\omega_m(t,\lambda,\tau)=
\prod_{1\leq i<j\leq m}\frac{\theta(t_i-t_j,\tau)}
{\theta(t_i-t_j+2\eta,\tau)}\prod_{j=1}^{m}\frac
{\theta(\lambda+t_j,\tau)}
{\theta(t_j-2\eta m,\tau)}.
\notag
\end{equation}
Notice that $\omega_m(t,\la,\tau)$ generates 
the one-dimensional space  $F^m_{2\eta m}(0,\lambda)$,
see Sec. \ref{sssym}.

For any function $f(t_1,...,t_m)$, let
 Sym $f(t_1,...,t_m)=\sum_{s\in S_m}[f(t_1,...,t_m)]_s$
be the symmetrization with respect to the action of the symmetric group
$S_m$ defined in Sec. \ref{sssym}.
For any $a=0,1,...,m$ introduce functions $\Sigma_a =$ Sym $\sigma_a$, $\Sigma_a'=$
Sym $\sigma'_a$ where
\bean
\sigma_a(t,\lambda,\tau)\,=
\,
\prod_{1\leq i<j\leq m-a}\frac{\theta(t_i-t_j)}
{\theta(t_i-t_j+2\eta)}
\prod_{j=1}^{m-a}\frac
{\theta(\lambda+2\eta (m-a)+t_j-2\eta m)}
{\theta(t_j-2\eta m)}
\notag
\\
\times
\prod_{m-a+1\leq i<j\leq m}\frac{\theta(t_i-t_j)}
{\theta(t_i-t_j+2\eta)}
\prod_{j=m-a+1}^{m}\frac
{\theta(\lambda+2\eta a +t_j-2\eta m + 4\eta (m-a))}
{\theta(t_j-2\eta m)}\,{}\,
\notag
\eean
and
\bean
\sigma_a'(t,\lambda,\tau)\,=
\,
\prod_{m-a+1\leq i<j\leq m}\frac{\theta(t_i-t_j)}
{\theta(t_i-t_j+2\eta)}
\prod_{j=m-a+1}^{m}\frac
{\theta(\lambda+2\eta a+t_j-2\eta m)}
{\theta(t_j-2\eta m)}
\notag
\\
\times
\prod_{1\leq i<j\leq m-a}\frac{\theta(t_i-t_j)}
{\theta(t_i-t_j+2\eta)}
\prod_{j=1}^{m-a}\frac
{\theta(\lambda+2\eta (m-a)+t_j-2\eta m + 4\eta a)}
{\theta(t_j-2\eta m)}\,
\notag
\\
\times
\prod_{i=1}^{m-a} \prod_{j=m-a+1}^{m}
\frac{\theta(t_i-t_j-2\eta)}
{\theta(t_i-t_j+2\eta)}\,{}\,.
\notag
\eean

It is easy to see that the functions $\Sigma_a, \Sigma_a'$ are elements
of the space $F^m_{2\eta m}(0,\lambda)$ and hence are proportional.

\begin{lemma}\label{} 
For any $a = 0,1,...,m$ we have
\bean\label{equality}
[m]!\,\omega_m(t,\lambda,\tau)\,=\,[a]![m-a]!\,\Sigma_a(t,\lambda,\tau)
\,=\,[a]![m-a]!\,\Sigma_a'(t,\lambda,\tau)
\eean
where $[l]!$ is the elliptic factorial.
\end{lemma}
\begin{proof} To prove the formula it is enough to compare the residues
of all functions at the point $(t_1,...,t_m)=(2\eta, 4\eta,...,2m\eta)$.
\end{proof}

To prove Part II of Theorem \ref{v} it is enough to notice that for $a = 1,...,m$
and any $s\in S_m$ we have
\bean
e^{2\pi i\,s\,a\, (4\eta m-p)}\,
e^{-\pi i { \mu\over 2\eta}(r+s\tau + 2\eta a)}\,
\int_{T^m}\,\Om(t,\tau,p,\eta)\,[\,\sigma_a(t,r+s\tau +
2\eta a,\tau)\,]_s\, \om_m(t,\mu,p)\,dt\,=
\notag
\\
e^{-\pi i { \mu\over 2\eta}(r+s\tau - 2\eta a)}\,
\int_{T^m}\,\Om(t,\tau,p,\eta)\,[\,\sigma_a'(t,r+s\tau -
2\eta a,\tau)\,]_s\, \om_m(t,\mu,p)\,dt\,.
\notag
\eean
For example, for $s =$ id, set $
F_a(t,\la,\mu,\tau,p)=
e^{-\pi i{\mu\la\over2\eta}}\Om(t,\tau,p,\eta)\sigma_a(t,\la,\tau) \om_m(t,\mu,p)
$ and $
F_a'(t,\la,\mu,\tau,p)=
e^{-\pi i{\mu\la\over2\eta}}\Om(t,\tau,p,\eta)\sigma_a'(t,\la,\tau) \om_m(t,\mu,p)\,.
$ Then we have
\bean
e^{2\pi i\,s\,a\, (4\eta m-p)}\,F_a(t_1,...,t_m,r+s\tau+2\eta a,\mu,\tau,p)=
\notag
\\
F_a'(t_1,...,t_{m-a}, t_{m-a+1}+p,...,t_m+p,r+s\tau-2\eta a,\mu,\tau,p,\eta).
\notag
\eean
This formula easily follows from the explicit formulas for $\sigma_a(t,\lambda,\tau),
\sigma_a'(t,\lambda,\tau)$ and transformation properties of $ \Omega(t,\tau,p)$
and $\omega(t,\mu,p)$.

For an arbitrary permutation $s\in S_m$ the proof is similar.

\subsection{Resonance relations for reduced coefficients}\label{transf}
For $\vec\La \in \C^n$, consider the standard
basis $e_J=e_{j_1}\otimes ... \otimes e_{j_n}$ in $V_{\vec\La}$.
Introduce a new {\it reduced} basis by
$$
E_J\,=\, {1\over [j_1]!\,[j_2]!\,...\,[j_n]!}\,e_J
$$
where $[j]!$ is the elliptic factorial of $j$.

Consider a hypergeometric solution $w(z,\la)=\sum_J \,w_J(z,\la)\, e_J$ 
of the qKZB equations with values in $V_{\vec\La}[0]$. Using the reduced
basis, the solution can be written as
$w(z,\la)=\sum_J \,W_J(z,\la)\, E_J$ where its {\it reduced} coefficients
are defined by
$$
W_J\,=\,  [j_1]!\,[j_2]!\,...\,[j_n]!\,w_J.
$$
Corollary \ref{part-1} can be reformulated in terms of the reduced coefficients.


>From now on we will be interested only in the special case of Corollary \ref{part-1}
when the integers $r, s$ are equal to zero. We have

\begin{corollary}\label{part-red}
The reduced coefficients of a hypergeometric solution obey
the resonance relations
\bean\label{res-hyp-red}
W_{M}(z, 2\eta (\La_j - a-b\,+\,\sum_{l=1}^{j-1}(\La_l-2m_l)))=
\\
W_{L}(z,2\eta (\La_j - a-b\,+\,\sum_{l=1}^{j-1}(\La_l-2m_l)))
\notag
\eean
for  $j, a, b, k, L, M$ defined in Part I of Theorem \ref{v} and the resonance
relations
\bean\label{res-hyp-2}
W_{M}(z, 2\eta (a-b))\,=\,
W_{L}(z, 2\eta (b-a))
\eean
for $ a, b, k, L, M$ defined in Part II of Theorem \ref{v}.
\end{corollary}

Consider a relation of Corollary \ref{part-red} involving a coefficient
$W_{M}(z,2\eta k )$ where $k$ is either $\La_j - a-b\,+\,\sum_{l=1}^{j-1}(\La_l-2m_l)$
or $a-b$.
Assume that the pair of indices $M, L$ in a relation of Corollary \ref{part-red}
is ordered from $M$ to $L$, then a relation will be called {\it a transformation
at level $k$ 
from $W_{M}$ to $W_{L}$}. The transformation in \Ref{res-hyp-red} will be denoted $T_j(k),\,
j=1,...,n-1,$ and the transformation in \Ref{res-hyp-2} will be denoted $T_n(k)$.

If $M=(m_1,...,m_n)$ and $j < n$, then the transformation $T_j(k)$ can be applied to
$W_M(z, 2\eta k)$ if there exists an integer $a$ such that
\bean\label{trans-j-1}
a \,\in \,[0,m_j+m_{j+1}]\, -\{m_j\}\,,
\\
\label{trans-j-2}
\sum_{l=1}^{j-1}(\La_l-2m_l) + \La_j - m_j - a = k\,.
\eean
Thus,
\bean\label{a}
a\, = \,\sum_{l=1}^{j-1}(\La_l-2m_l) + \La_j - m_j - k\,.
\eean
The result of the transformation $T_j(k)$ is the index
\bean\label{result-j-1}
{}
\\
L\,=\,(m_1,...,m_{j-1}, \,\sum_{l=1}^{j-1}(\La_l-2m_l) + \La_j - m_j - k,\,
m_{j+1} + k - \sum_{l=1}^{j}(\La_l-2m_l),\,m_{j+2},...,m_n)\,.
\notag
\eean
and the value $W_L(z,2\eta k)$.

If $M=(m_1,...,m_n)$, then the transformation $T_n(k)$ can be applied to
$W_M(z, 2\eta k)$ if there exists an integer $a$ such that
\bean\label{trans-n-1}
a \,\in \,[0,m_n+m_{1}]\, -\{m_n\}\,,
\\
\label{trans-n-2}
m_n - a \,= \,k\,.
\eean
Note that equation \Ref{trans-n-2} can be written as
$\sum_{l=1}^{n-1}(\La_l-2m_l) + \La_n - m_n - a = k$ similarly to 
\Ref{trans-j-2}.

The result of the transformation $T_n(k)$ is the index
\bean\label{result-j-1}
L\,=\,(m_1 + k , m_2, ...,m_{n-1}, m_n - k)
\eean
and the value $W_L(z,\,-2\eta k)$.


\section{Regularity and Resonance Relations}\label{Regularity}

\subsection{Statement of results}

Hypergeometric solutions are functions holomorphic with respect to variable $\la$.
At the same time they satisfy the qKZB equations, $v(z_1,...,z_j+p,...,z_n)$
$=K_j(z_1,...,z_n,\tau,p)v(z_1,...,z_n)$,
where the qKZB operators $K_j$ have poles with respect
to  $\la$.
Also if $v(z_1,...,z_n)$ is a hypergeometric solution with values in 
$V_{\La_1}\otimes ... \otimes V_{\La_n}\,[0]$, then
for any $j=1,\ldots,n-1$, the function
$P^{(j,j+1)}$ $R_{\La_j, \La_{j+1}}^{(j,j+1)}
(z_{j+1}-z_j,\la-2\eta \sum_{l=1}^{j-1} h^{(l)},\tau)
v(z_1,\ldots,z_{j+1},z_j,\ldots ,z_n,\la)$
is a hypergeometric solution of the qKZB equations with values in 
$(V_{\La_1}\otimes\ldots
V_{\La_{j+1}}\otimes V_{\La_j}\ldots\otimes V_{\La_n})[0]$. The R-matrix has $\la$-poles
while the new solution remains holomorphic with respect to $\la$.
To preserve the $\la$-holomorphy under the action of these transformations
a function has to possess special properties. 
In this Section we show that these  properties are  the resonance relations discussed 
in Section \ref{Verma}.

Let $n>1$ be a natural number.
Let $\vec\La=(\Lambda_1,\dots,\Lambda_n) \in \C^n$ be such that
$\La_1+...+\La_n=2m$ for some natural number $m$. 
Let $u(\la)=\sum_{m_1+...+m_n=m}u_{m_1,...,m_n}(\la)\,e_{m_1}\otimes ... 
\otimes e_{m_n}$ be a $V_{\vec\La}[0]$-valued holomorphic function of $\la$.
Let $j$ be a natural number, $1\leq j\leq n-1$.
Let $(z_1,...,z_n) \in \C^n$. We say that the 
function $u(\la)$ obeys  {\it the resonance relation $C_j$ 
with respect to the vector $\vec\La$ at the point} $(z_1,...,z_n)$ if the following holds.
Let  $m_1,...,m_{j-1},k,m_{j+2},...,m_n$
be any non-negative integers such that $m_1+...+m_{j-1}+k+m_{j+2}+...+m_n=m$. Let 
$a,b$ be any integers such that $a\neq b,\, 0\leq a,b \leq k$. Let
$M=(m_1,...,m_{j-1},a,k-a,m_{j+2},...,m_n)$ and
$L=(m_1,...,m_{j-1},b,k-b,m_{j+2},...,m_n)$. Then 
we require that
\bean\label{C-j}
{}
\\
{}\,{[a]!\,[k-a]!}\,e^{2\pi i\,s\,a\, (z_{j+1}-z_j +\eta\La_{j+1}+\eta\La_{j})}\,
u_{M}(r+s\tau\,+\,2\eta (\La_j - a-b\,+\,\sum_{l=1}^{j-1}(\La_l-2m_l)))\,=
\notag
\\
\,[b]!\,[k-b]!\,
e^{2\pi i\,s\,b\, (z_{j+1}-z_j +\eta\La_{j+1}+\eta\La_{j})}\,
u_{L}(r+s\tau\,+\,
2\eta (\La_j - a-b\,+\,\sum_{l=1}^{j-1}(\La_l-2m_l)))\,.
\notag   
\eean
We say that the function $u(\la)$ obeys  {\it the resonance relation $C_n$
with respect to the vector $\vec\La$ at the point} $(z_1,...,z_n)$ 
if the following holds. Let  $m_2,...,m_{n-1},k$
be any non-negative integers such that $m_2+...+m_{n-1}+k=m$. Let 
$a,b$ be any non-negative integers that $a\neq b,\, 0\leq a,b \leq k$. Let
$M=(k-a,m_2,...,m_{n-1},a)$ and
$L=(k-b,m_2,...,m_{n-1},b)$. Then we require
\bean\label{C-n}
{}
\\
{[a]!}\,[k-a]!\,
e^{2\pi i\,s\,a\, (z_{1}-z_n +\eta\La_{1}+\eta\La_n-p)}\,
u_{M}(r+s\tau\,+\,2\eta (a-b))\,=\,
\notag
\\
\,[b]!\,[k-b]!\,
e^{2\pi i\,s\,b\, (z_{1}-z_n +\eta\La_{1}+\eta\La_n-p)}\,
u_{L}(r+s\tau\,+\,2\eta (b-a))\,.
\notag
\eean
For $z\in \C$ and $1\leq j\leq n-1$, we introduce an operator
$s_j(z) : {\mathcal F} (V_{\La_1,...,\La_n})$ $ \to$ $
{\mathcal F} (V_{\La_1,...,\La_{j-1},\La_{j+1},\La_j,\La_{j+2},...,\La_n})$
by
\[
s_j(z)u(\lambda)
=P^{(j,j+1)}R_{\La_j,\La_{j+1}}(z,\lambda - 2\eta\sum_{l=1}^{j-1}
h^{(l)})^{(j,j+1)}u(\lambda)\,.
\]
\begin{thm}\label{reg-thm}
Let $n>1$ be a natural number.
Let $\vec\La=(\Lambda_1,\dots,\Lambda_n)$ be a generic vector in $\C^n$ with property
$\La_1+...+\La_n=2m$ for some natural number $m$.
Let $(z_1,...,z_n)\in \C^n$ be a generic point.
Suppose that a $V_{\vec\La}[0]$-valued function $u(\lambda)$ is holomorphic
and   obeys the  conditions
$C_l, l=1,...,n$,  with respect to the vector $(\Lambda_1,\dots,\Lambda_n)$
at the point $(z_1,...,z_n)$.
Then
\begin{enumerate}
\item[I.] For $1\leq j \leq n-1$, the function $s_j(z_j-z_{j+1})u(\la)$ is holomorphic
and obeys the  conditions $C_l,\, l=1,...,n$, with respect to the vector 
$(\La_1,...,\La_{j+1},\La_j,...,\La_n)$ at the point 
$(z_1,...,z_{j+1},z_j,...,z_n)$.
\item[II.] 
For $1\leq j \leq n$, the function $K_j(z_1,...,z_n,\tau,p)u(\la)$ is holomorphic
and obeys the conditions
$C_l,\,l=1,...,n,$ with respect to the vector $(\La_1,...,\La_n)$
at the point $(z_1,...,z_{j-1},,z_j+p,z_{j+1},...,z_n)$.
\end{enumerate}
Conversely, suppose that for some generic vector $(\Lambda_1,...,\Lambda_n)$
with $\La_1+...+\La_n=2m$, generic point $(z_1,...,z_n)\in \C^n$ and a holomorphic
$V_{\vec\La}[0]$-valued function  $u(\la)$ we have the following properties:
\begin{enumerate}
\item[a)] For any $1\leq j \leq n-1$, the function $s_j(z_j-z_{j+1})u(\la)$ of $\la$
is a holomorphic function of $\la$. 
\item[b)] The function $K_{2}(z_1,...,z_n) u(\la)$ is
a holomorphic function of $\la$.
\end{enumerate}
Then $u(\la)$ 
obeys the  conditions $C_l,\,l=1,...,n,$ with respect to the vector $(\Lambda_1,...,\Lambda_n)$ 
at the point $(z_1,...,z_n)$.

\end{thm}
{\bf Remark.} The vector $(\Lambda_1,...,\Lambda_n)$ is generic if all $\lambda$-poles of the R-matrices
appearing in the proof below are simple.
The point $(z_1,...,z_n)$ is generic  if all R-matrices $R_{\La_i,\La_j}(z_i-z_j \,(+p),\la)$
appearing in the proof below are well defined meromorphic functions of $\la$.

\subsection{Proof of Theorem \ref{reg-thm}  }\label{more}
Let $\vec\La^\vee=(\La_n,\La_1,\La_2,...,
\La_{n-1})$. Let $\Dl:V_{\vec\La}\to V_{\vec\La^\vee}$ be the linear operator
defined by $\Dl\,=\, \Gamma_1 P^{(1,2)} P^{(2,3)} ... P^{(n-1,n)}$.

We have $s_j(z)s_j(-z)={\mathrm{Id}}$, by the ``unitarity'' of the $R$-matrix.
\begin{lemma}\label{l01} 
Let $z=(z_1,\dots,z_n)\in\C^n$.
Then the $j$-th qKZB operator can be written as:
\begin{eqnarray*}
K_j(z,\tau,p)&=&s_{j-1}(z_j-z_{j-1}+p)
\cdots s_2(z_j-z_2+p)s_1(z_j-z_1+p)
\\ & &\Delta\,
s_{n-1}(z_j-z_n)\cdots s_j(z_j-z_{j+1}).
\end{eqnarray*}
\end{lemma}

\begin{proof}
Recall that if we set $R_{jk}=R(z_j-z_k(+p),\lambda-2\eta\sum_{j<k,\, j\neq l}
h^{(l)})$,  where we add $p$ only if $k<j$, 
and $\hat R_{jk}=P^{(j,k)}R_{jk}$, then
\begin{eqnarray*}
K_j(z,\tau,p)&=&R_{jj-1}\cdots R_{j1}\Gamma_jR_{jn}\cdots R_{jj+1}\\
      &=&P^{(j,j-1)}\hat R_{jj-1}\cdots P^{(j,1)}\hat R_{j1}
         \Gamma_jP^{(j,n)}\hat R_{jn}\cdots P^{(j,j+1)}\hat R_{jj+1}\\
      &=&s_{j-1}(z_j-z_{j-1}+p)\cdots s_1(z_j-z_1+p)\\
      & &P^{(j,j-1)}\cdots P^{(j,1)}\Gamma_jP^{(j,n)}\cdots P^{(j,j+1)}
           s_{n-1}(z_j-z_{n})\cdots s_j(z_j-z_{j+1}).
\end{eqnarray*}
We have moved the $P^{(j,k)}$ to the middle using the rule
$P^{(i,k)}X^{(k)}=X^{(i)}P^{(i,k)}$. Using the same rule to move
$P^{(j,j-1)}$ to the right, we obtain
\[
P^{(j,j-1)}\cdots P^{(j,1)}\Gamma_jP^{(j,n)}\cdots P^{(j,j+1)}
=P^{(j-1,j-2)}\cdots P^{(j-1,1)}\Gamma_{j-1}P^{(j-1,n)}\cdots P^{(j-1,j)}.
\]
This expression is therefore independent of $j$, and can thus be replaced
by its value at $j=1$, which is, by definition, $\Delta$.
\end{proof}

\begin{lemma}\label{l02} 
Let $(\La_1,...,\La_n)$ be  a generic vector.
Let $(z_1,...,z_n)\in \C^n$ be a generic point.
Let $u(\la)$ be a   holomorphic function  with values in $V_{\vec\La}$.
Let  $1\leq j\leq n-1$.
Then,  $s_j(z_j-z_{j+1})u(\lambda)$ is holomorphic
if and only if $u(\la)$
obeys  the condition $C_j$ with respect to the vector 
$(\La_1,...,\La_n)$ at the point $(z_1,...,z_n)$.
\end{lemma}
\begin{proof}
For generic $\La_j, \La_{j+1}$, the R-matrix $R_{\La_j,\La_{j+1}}(z_j-z_{j+1},\la 
- 2\eta\sum_{l=1}^{j-1} h^{(l)})$
has only simple $\la$-poles. The kernel of the residue of this matrix at a
$\la$-pole is described in Theorem \ref{La-poles}. 
It follows from Theorem \ref{La-poles} that a
function $u(\la)$ obeys the condition $C_j$ with respect to the vector
$(\La_1,...,\La_n)$ at a point $(z_1,...,z_n)$ if and only if for any
$\la$-pole $\la_0$ of  $R_{\La_j,\La_{j+1}}(z_j-z_{j+1},\la- 2\eta\sum_{l=1}^{j-1} h^{(l)})$ 
the vector $u(\la_0)$ belongs to the kernel of the residue of
$R_{\La_j,\La_{j+1}}(z_j-z_{j+1},\la - 2\eta\sum_{l=1}^{j-1} h^{(l)})$ at $\la_0$. Thus the function
$R_{\La_j,\La_{j+1}}(z_j-z_{j+1},\la- 2\eta\sum_{l=1}^{j-1} h^{(l)})u(\la)$ is holomorphic
if and only if $u(\la)$
obeys the condition $C_j$. The Lemma is proved since $s_j=PR$ and the permutation $P$ is regular.
\end{proof}

\begin{lemma}\label{l03}
Let $1\leq k\leq n-1$.
If $u(\la)$ is a holomorphic function with values in $V_{\vec\La}$
and obeys the condition
$C_k$ with respect to the vector $(\La_1,...,\La_n)$ 
at a generic point $(z_1,...,z_n)$,  then for any $j$
the function $s_j(z_j-z_{j+1})u(\la)$ obeys the condition $C_k$ with respect to the vector 
$(\La_1,...,\La_{j+1},\La_j,...,\La_n)$
at the point $(z_1,...,z_{j+1},z_j,...,z_n)$.
\end{lemma}
\begin{proof}
The only non-trivial cases are when $|j-k|\leq 1$. 

Let $j=k$. By Lemma \ref{l02}, $s_j(z_j-z_{j+1})u(\la)$ is holomorphic
and $s_j(z_{j+1}-z_{j})s_j(z_j-z_{j+1})u(\la)=u(\la)$
is holomorphic. Hence by Lemma \ref{l02}, $s_j(z_j-z_{j+1})u(\la)$ obeys $C_j$.

Let $j=k+1$. The fact, that $s_j(z_j-z_{j+1})u(\la)$ satisfies the condition $C_{j+1}$
if $u(\la)$ satisfies the condition $C_{j+1}$, follows from identities \Ref{coeff-1}
and is similar to the proof of the identities \Ref{coeff-1}.

The case $j=k-1$ follows from identities \Ref{coeff-2}.

\end{proof}

\begin{lemma}\label{l05}
If $u(\la)$ is a holomorphic function with values in $V_{\vec\La}$
and obeys the condition
$C_n$ with respect to the vector $(\La_1,...,\La_n)$
at a generic point $(z_1,...,z_n)$,  then for any $j$
the function $s_j(z_j-z_{j+1})u(\la)$ obeys  the condition $C_n$ with respect to the vector 
$(\La_1,...,\La_{j+1},\La_j,...,\La_n)$
at the point $(z_1,...,z_{j+1},z_j,...,z_n)$.
\end{lemma}
\begin{proof} The statement is trivial if $2\leq j\leq n-2$.
Suppose that $j=1$, and let
\[
u(\lambda)=\sum_{m_1,m_2,m_n}e_{m_1}\otimes e_{m_2}\otimes u_{m_1,m_2,m_n}(\lambda)\otimes e_{m_n},
\]
with $u_{m_1,m_2,m_n}(\lambda)\in V_{\La_3}\otimes ... \otimes V_{\La_{n-1}}$.
Let
\[
v(\lambda)=s_1(z_1-z_2)u(\la)=
\sum_{l_2,l_1,m_n}e_{l_2}\otimes e_{l_1}\otimes v_{l_2,l_1,m_n}(\lambda)\otimes e_{m_n},
\]
where
\bean\label{nuzhnoe}
v_{l_2,l_1,m_n}(\lambda)=\sum_{m_1,m_2}
R_{\La_1,\La_2}(z_1-z_2,\la)^{l_1,l_2}_{m_1,m_2}u_{m_1,m_2,m_n}(\lambda).
\eean
Our goal is to prove the identity
\bean\label{conj-identity}
[l_2]![m_n]!e^{2\pi i s m_n(z_2-z_n+\eta\La_2+\eta\La_n-p)}
v_{l_2,l_1,m_n}(r+s\tau+2\eta(m_n-\tilde m_n))=
\\
{[\tilde l_2]!}[\tilde m_n]!e^{2\pi i s \tilde m_n(z_2-z_n+\eta\La_2+\eta\La_n-p)}
v_{\tilde l_2,l_1,\tilde m_n}(r+s\tau+2\eta(\tilde m_n-m_n))\,
\notag
\eean
for any $l_1,l_2, \tilde l_2, m_n \tilde m_n$. Using \Ref{nuzhnoe},
we can rewrite \Ref{conj-identity} as
\bean\label{second}
[l_2]![m_n]!e^{2\pi i s m_n(z_2-z_n+\eta\La_2+\eta\La_n-p)}\times &{}&
\\
\sum_{m_2}R_{\La_1,\La_2}(z_1-z_2, r+s\tau+ &{}&2\eta(\tilde l_2-l_2))^{l_1,l_2}
_{l_1+l_2-m_2,m_2}\times \notag \\
u_{l_1+l_2-m_2,m_2 ,m_n}&{}&(r+s\tau+2\eta(m_n-\tilde m_n))=
\notag
\\
{[\tilde l_2]!}[\tilde m_n]!e^{2\pi i s \tilde m_n(z_2-z_n+\eta\La_2+\eta\La_n-p)}\times &{}&
\notag
\\
\sum_{m_2}R_{\La_1,\La_2}(z_1-z_2, r+s\tau+ &{}&2\eta( l_2-\tilde l_2))^{l_1,\tilde l_2}
_{l_1+\tilde l_2-m_2,m_2}\times \notag \\
u_{l_1+\tilde l_2-m_2,m_2 ,\tilde m_n}&{}&(r+s\tau+2\eta(\tilde m_n- m_n)).
\notag
\eean
Now \Ref{second} follows from identities \Ref{coeff-1} and  the condition $C_n$ for
the function $u(\la)$.

A similar consideration shows that the Lemma holds also for $j=n-1$.
\end{proof}

\begin{lemma}\label{posled}
Let $u(\la)$ be  a meromorphic function with values in $V_{\vec\La}$.
Then for any $j=1,...,n$, the function $\Dl u(\la)$ is holomorphic and
obeys $C_j$ with respect to the vector
$(\La_n,\La_1,\La_2,...,\La_{n-1})$ at the point $(z_n+p,z_1,...,z_{n-1})$ if and only if
the function $u(\la)$ is holomorphic
and 
obeys $C_{j-1}$ with respect to the vector
$(\La_1,\La_2,...,\La_{n})$ at the point $(z_1,z_2,...,z_n)$ (where $C_0=C_n$).
\end{lemma}
The Lemma is proved by direct verification.

The proof of Theorem \ref{reg-thm} can be concluded.
It follows from Lemma \ref{l02}
that if $u(\la)$  is holomorphic and
obeys the conditions $C_1,...,C_n$, then
 the functions $s_j(z_j-z_{j+1})u(\la)$ are holomorphic and obey $C_1,...,C_n$. 
Using the fact that
$K_j(z_1,\dots,z_n)$ can be expressed through $s_k(z)$ and 
$\Delta$ (Lemma \ref{l01}), we deduce
that the functions $K_j(z_1,\dots,z_n,\tau,p)u(\la)$ are holomorphic and obey
$C_1,...,C_n$.

Now we prove the  "conversely" part. By Lemma \ref{l02}, the function
$u(\la)$ obey conditions $C_1,...,C_{n-1}$ with respect to the vector
$(\La_1,...,\La_n)$ at the point $(z_1,...,z_n)$. We have
$K_2(z_1,...,z_n,\tau,p)=s_{1}(z_2-z_{1}+p)\Delta
s_{n-1}(z_2-z_n)... s_2(z_2-z_{3})$.
For $j=2,...,n-1,$ set $v_j(\la)=s_{j}(z_2-z_{j+1})... s_2(z_2-z_{3})u(\la)$.
By Lemmas \ref{l02} and \ref{l03}, the functions $v_2(\la),...,v_{n-1}(\la)$
are holomorphic and the function $v_{n-1}(\la)$ obeys $C_1,...,C_{n-1}$
with respect to the vector $(\La_1,\La_3,...,\La_n,\La_2)$ at the point
$(z_1,z_3,...,z_n,z_2)$. By Lemma \ref{posled}, the function $\Dl v_{n-1}(\la)$
is holomorphic and obeys $C_2,...,C_n$ with respect to the vector
$(\La_2,\La_1,\La_3,...,\La_n)$ at the point $(z_2+p,z_1,z_3,...,z_n)$.
The function $s_{1}(z_2-z_{1}+p)\Delta v_{n-1}(\la)$ is holomorphic.
By Lemma \ref{l02}, the function $\Dl v_{n-1}(\la)$ obeys $C_1$
with respect to the vector $(\La_2,\La_1,\La_3,...,\La_n)$ at 
the point $(z_2+p,z_1,z_3,...,z_n)$. By Lemma \ref{posled},
the function $v_{n-1}(\la)$ obeys $C_n$ with respect to the vector
$(\La_1,\La_3,...,\La_n,\La_2)$ at the point $(z_1,z_3,...,z_n,z_2)$.
By Lemma \ref{posled}, the function $v_{n-2}(\la)=s_{n-1}(z_n-z_2)v_{n-1}(\la)$ 
obeys $C_n$ with respect to the vector
$(\La_1,\La_3,...,\La_{n-1},\La_2,\La_n)$ at the point $(z_1,z_3,...,z_{n-1},z_2,z_n)$.
Repeating this procedure we conclude that $u(\la)$ obeys $C_n$ with respect to the vector
$(\La_1,...,\La_n)$ at the point $(z_1,...,z_n)$. The "conversely" part 
is proved.

\section{Remarks on Semiclassical Limit of Resonance Relations}\label{remarks}

We consider the semiclassical limit of the resonance relations in the simplest nontrivial case.
 Let $n=1$, $\vec \Lambda = 2m$. 
The universal hypergeometric function $u$ associated with $\vec\La$ is a scalar function.
The resonance relations are $ u(2\eta a,\mu,\tau,p,\eta,\vec\La)$ $=$
$u(-2\eta a,\mu,\tau,p,\eta,\vec\La)$ for $a = 1,...,m$.
\begin{lemma}\label{HEAT}
Let $v(\la,\eta)$ be a smooth function such that
for small $\eta$ we have $ v(2\eta a,\eta)=v(-2\eta a,\eta)$ for $a = 1,...,m$.
Then 
\bean\label{semiclass}
 {\partial v\over \partial \la}  (0,0)= {\partial^3 v\over \partial \la^3}  (0,0)=...
={\partial^{2m-1} v\over \partial \la^{2m-1}}  (0,0)=0.
\notag
\eean
Moreover, for the anti-symmetric function $Av(\la)=v(\la)-v(-\la)$ we have
${\partial^{a} Av\over \partial \la^{a}}  (0,0)=0$ for $a = 0,1,...,2m$.
\end{lemma} 
\begin{proof} 
\bean
2{\partial^{2m-1} v\over \partial \la^{2m-1}}  (0,0)=\text{lim}_{\eta\to 0}
{1\over (2\eta)^{2m-1}}(\sum_{k=0}^{k=2m-1}(-1)^k &{}& {2m-1\choose k}
v(2\eta(m-k),\eta)
+
\notag
\\
 \sum_{k=0}^{k=2m-1}&{}&(-1)^k  {2m-1\choose k} v(2\eta(m-k-1),\eta)\,)\,.
\notag
\eean
It is easy to see that the expression in the parentheses can be written
as $\sum_{a=1}^m c_a (v(2\eta a,\eta)-v(-2\eta a,\eta))$ for suitable coefficients
$c_a$. Thus the expression is zero.
The other equalities to zero are proved similarly.
\end{proof}
The semiclassical limit of the qKZB equations is the limit $\eta, p \to 0$ with
$p/\eta=-2\kappa$  fixed. In this limit the qKZB equations after
a suitable normalization turn into the KZB differential equations
\cite{FW, FV5}. In the classical case, in addition to the KZB equations,
that are associated to the variation of the marked points on the elliptic curve,
one also has an equation associated to the variation of the modulus $\tau$ of
the elliptic curve. This additional differential equation, compatible with
the KZB equations, is called the KZB {\it heat equation} (its difference
version is suggested in \cite{FV3}). For example,
if $n=1$ and $\vec\La=2m$, then the heat differential equation takes the form
\bean\label{heat}
2\pi i \kappa {\partial u\over \partial \tau}={\partial^2 u\over \partial \la^2}
-m(m+1)\,\wp (\la,\tau)u
\eean
where $\wp$ is the Weierstrass function. The Weierstrass function is an even function of $\la$
with second order pole at $\la=0$. The solutions of the
heat equation in the neighborhood of $\la=0$ are meromorpchic functions.
The Weyl reflection $\la \to -\la$ preserves the equation and acts on the space of  solutions.
If $m$ is even, then a $\la$-odd solution  has a $\la$-pole of order $m$ at the origin
and a $\la$-even solution has zero of order $m+1$.
If $m$ is odd, then a $\la$-even solution has a $\la$-pole of order $m$ and a $\la$-odd solution
has zero of order $m+1$. Change the variable in \Ref{heat}, $w\,=\,\theta(\la,\tau)^m\,u$. Then
the Weyl reflection still acts on the space of solutions of the transformed heat equation.
Any $\la$-meromorphic solution of the transformed heat equation is now regular at the origin.
Its Taylor expansion has the form $w(\la,\tau)= \sum_{k=0}^{m-1}w_{2k}(\tau)\la^{2k}
+\sum_{k>2m}w_{k}(\tau)\la^{k}$, and the Taylor expansion of a $\la$-odd solution of the transformed
equation has zero of order $2m+1$.

It is expected that the semiclassical limit of the
universal hypergeometric function $ u(\la,\mu,\tau,p,\eta,\vec\La)$  after multiplication by a suitable
function of $\tau$ will give solutions of the transformed heat equation (depending
on $\mu$ as a parameter). Lemma \ref{HEAT} shows that the semiclassical limit of resonance relations
turn into the regularity properties of the solutions of the heat equation.

Equation \Ref{heat} can be considered as the equation for horizontal sections of a connection
over the upper half plane whose fiber is the space of functions of $\la$. 
If $\kappa$ is an integer, then the connection has an invariant finite dimensional
subbundle of conformal blocks coming from conformal field theory. The subbundle consists of
Weyl anti-symmetric theta functions  satisfying the vanishing conditions that we find in the 
semiclassical limit of our resonance relations, see \cite{FG}, \cite{FW}, \cite{EFK}.

\section{The Weyl reflection}\label{weyl-refl}

\subsection{The Weyl reflection, R-matrices and qKZB operators}
For a positive integer $\Lambda$, let $L_\La(z)$ be the $\La + 1$ 
dimensional $E_{\tau,\eta}(sl_2)$ evaluation module with the standard basis
$e_j,\, j=1,...,\Lambda$. Introduce a new basis $E_j$ by $E_j=e_j/[j]!$
where $[j]!$ is the elliptic factorial. Define an involution
$$
s_\La\,:\,L_\La\,\to L_\La,\, \qquad E_j\,\mapsto \, E_{\La-j},\qquad
\text{for} \qquad j=0,...,\La\,.
$$
\begin{thm}\label{Rm-weyl}
For natural numbers $\La_1, \La_2$, the R-matrix $R_{\La_1,\La_2}(z_1-z_2,\la)
\in \End \,(L_{\La_1}(z_1)\otimes L_{\La_1}(z_2))$
obeys
\bean\label{rm-weyl}
R_{\La_1,\La_2}(z_1-z_2,-\la)\,{}\,s_{\La_1}\otimes s_{\La_2}\,=\,
s_{\La_1}\otimes s_{\La_2}\,{}\,
R_{\La_1,\La_2}(z_1-z_2,\la)\,.
\eean
\end{thm}

Let $\vec\La=(\La_1,...,\La_n)$ be a vector of natural numbers.
Consider the qKZB equations with values in $L_{\vec\La}[0]$ and the
corresponding qKZB operators $K_i(z,\tau,p)$ acting on the space 
${\mathcal F} (L_{\vec\La}[0])$ of  meromorphic functions
of $\lambda\in\C$ with values in the zero weight space $L_{\vec\Lambda}[0]$.

Define {\it the Weyl reflection} $S\, :\, {\mathcal F} (L_{\vec\La}[0])\,\to\,
{\mathcal F} (L_{\vec\La}[0])$ by
\bean\label{weyl}
 (Sf)(\la)\,=\, s_{\La_1}\otimes ... \otimes s_{\La_n}\,{} f(-\la).
\eean 
\begin{thm}\label{Qkzb-weyl} For all $i$, the qKZB operators obey
the relation
\bean\label{qkzb-weyl}
S\,K_i(z,\tau,p)\,=\,K_i(z,\tau,p)\,S\,.
\eean
\end{thm}
The Theorem is a direct corollary of \Ref{rm-weyl} and the obvious
relation $S\,\Gamma_i\,=\,\Gamma_i\,S$.

\begin{corollary}
If $u$ is a solution of the qKZB equations, then $Su$ is.
\end{corollary}

\subsection{Proof of Theorem \ref{Rm-weyl}} The Theorem is proved
by induction on highest weights $\La_1, \La_2$. First we prove
some properties of R-matrices.

Let $V, W$ be representations of $E_{\tau,\eta}(sl_2)$  with
L-operators $L_V\in \End\,(\C^2\otimes V)$ and $L_W\in \End\,(\C^2\otimes W)$.
An R-matrix $R_{V,W}(\la) \in \End\,(V\otimes W)$ is an operator
depending on $\la$ and such that $R_{V,W}(\la)P_{V,W}$ is an 
 $E_{\tau,\eta}(sl_2)$  intertwiner from 
$W\otimes V$ to $V\otimes W$, thus
\bean
R_{VW}(\lambda)^{(23)}P_{VW}L_W(z,\lambda-2\eta h^{(3)})^{(12)}
L_V(z,\lambda)^{(13)}=
\notag
\\
L_V(z,\lambda-2\eta h^{(3)})^{(12)}
L_W(z,\lambda)^{(13)}R_{VW}(\lambda-2\eta h^{(1)})^{(23)}P_{VW}.
\notag
\eean
So the relation determining $R$ is
\bean
R_{VW}(\lambda)^{(23)}L_W(z,\lambda-2\eta h^{(2)})^{(13)}
L_V(z,\lambda)^{(12)}=
\notag
\\
L_V(z,\lambda-2\eta h^{(3)})^{(12)}
L_W(z,\lambda)^{(13)}R_{VW}(\lambda-2\eta h^{(1)})^{(23)}.
\notag
\eean

\begin{lemma}\label{intertw} 
Let $V_1, V_2, V_3$ be representations of 
$E_{\tau,\eta}(sl_2)$. Then
$$
R_{V_1,V_2\otimes V_3}(\lambda)=
R_{V_1,V_2}(\lambda-2\eta h^{(3)})^{(12)}
R_{V_1,V_3}(\lambda)^{(13)}
$$ 
is an $R$-matrix for representations $V_1,\, V_2\otimes V_3$
and
$$
R_{V_1\otimes V_2, V_3}(\lambda)=
R_{V_2,V_3}(\lambda)^{(23)} 
R_{V_1,V_3}(\lambda-2\eta h^{(2)})^{(13)}
$$ 
is an $R$-matrix for representations $V_1\otimes V_2,\, V_3$.
\end{lemma}
\begin{proof} We prove the first relation.
\bean
L_{V_1}(z,\lambda-2\eta(h^{(3)}+h^{(4)}))^{(12)}L_{V_2\otimes V_3}(z,\lambda)^{(1,34)}
R_{V_1,V_2\otimes V_3}(\lambda-2\eta h^{(1)})^{(2,34)}
\notag
\\
=
L_{V_1}(z,\lambda-2\eta(h^{(3)}+h^{(4)}))^{(12)}
L_{V_2}(z,\lambda-2\eta h^{(4)})^{(13)}
L_{V_3}(z,\lambda)^{(14)}
\notag
\\
\times\,R_{V_1,V_2}(\lambda-2\eta(h^{(1)}+h^{(4)})^{(23)}
R_{V_1,V_3}(\lambda-2\eta h^{(1)})^{(24)}
\notag
\\
=
L_{V_1}(z,\lambda-2\eta(h^{(3)}+h^{(4)}))^{(12)}
L_{V_2}(z,\lambda-2\eta h^{(4)})^{(13)}
\notag
\\
\times
R_{V_1,V_2}(\lambda-2\eta(h^{(1)}+h^{(4)})^{(23)}
L_{V_3}(z,\lambda)^{(14)}
R_{V_1,V_3}(\lambda-2\eta h^{(1)})^{(24)}
\notag
\eean
At this point it is easy to move the $R$-matrices to the left,
with the desired  result.
\end{proof}

For a positive integer $\Lambda$, let $L_\La(z)$ be the $\La + 1$ 
dimensional $E_{\tau,\eta}(sl_2)$ evaluation module with the basis
$E_j,\, j=1,...,\Lambda$.  Then  Theorem 10 in \cite{FV1} says that
the map
\bean
j_\Lambda\,:\,L_{\Lambda+1}(z)\to L_{1}(z-\eta\Lambda)\otimes L_\Lambda(z+\eta)
, \qquad E_k \,\mapsto\,\sum_l E_l\otimes E_{k-l}
\notag
\eean
define  an embedding of $E_{\tau,\eta}(sl_2)$ modules.

\begin{lemma}\label{recurs}
For natural numbers $\Lambda_1, \La_2$, let
$R_{\Lambda_1,\Lambda_2}(z,\lambda)\in \End\,(L_{\La_1}\otimes L_{\La_2})$ 
be the R-matrix constructed
in Sec. \ref{fs}. 
Then the  R-matrix $R_{\Lambda_1,\Lambda_2}(z,\lambda)$ obeys the
recursion relations
\bean
j_{\Lambda_1}\otimes 1 \,R_{\Lambda_1+1,\Lambda_2}(z,\lambda)=
R_{\Lambda_1,\Lambda_2}(z+\eta,\lambda)^{(23)}
R_{1,\Lambda_2}(z-\eta\Lambda_1,\lambda-2\eta h^{(2)})^{(13)}
j_{\Lambda_1}\otimes 1,
\notag
\\
1\otimes j_{\Lambda_2}\,R_{\Lambda_1,\Lambda_2+1}(z,\lambda)=
R_{\Lambda_1,1}(z+\eta\Lambda_2,\lambda-2\eta h^{(3)})^{(12)}
R_{\Lambda_1,\Lambda_2}(z-\eta,\lambda)^{(13)}
1\otimes j_{\Lambda_2}.
\notag
\eean
\end{lemma}
\begin{proof} The Lemma follows from Lemma \ref{intertw} and the normalization
property $R_{\La, \Mu}(z,\la)\, E_0\otimes E_0\,=\,E_0\otimes E_0$.
\end{proof}

Now Theorem \ref{Rm-weyl} is proved
by induction in $\Lambda_1$, $\Lambda_2$. For $\Lambda_1=
\Lambda_2=1$ it is checked by direct calculation. The general case follows
using 
$$
(s_{1}\otimes s_{\Lambda})\circ j_\Lambda= j_{\Lambda+1}\circ s_{\Lambda+1}\,.
$$

\section{Vanishing conditions for hypergeometric solutions with
values in finite dimensional modules}\label{vanishing-cond}

\subsection{Fusion rules for $sl_2$ }\label{fus}
For a non-negative integer $a$, let $L_a$ denote the finite dimensional
$sl_2$-module with highest weight $a$. Let $a, b, c$ be non-negative integers.
We say that $a, b, c$ obey {\it the fusion rules} for $sl_2$
if the $sl_2$-module $L_a\otimes L_b$ contains $L_c$ as a submodule.
It means that $a-b \in \{-c, -c+2,...,c-2, c\}$ and $c\leq a+b$.
If the triple $a, b, c$ obeys the fusion rules, then each of its its permutation does.

Let $\vec\La=(\La_1,...,\La_n)$ be an $n$-tuple of positive integers.
Let $a=(a_1,...,a_n)$ be an $n$-tuple of integers. We say that $a$ obeys
{\it the fusion rules for $sl_2$
with respect to $\vec\La$}, if $a_1,...,a_n$
are non-negative
and for every $j=1,...,n$, the triple $a_{j-1}, a_{j}, \La_j$ obeys
the fusion rules. Here we assume that $a_{0}=a_n$.

Define {\it the weight vector},  $w(a)=(w_1,...,w_n)$, of an $n$-tuple
$a$ by $w_j=a_j-a_{j-1}$.
The $n$-tuple $a$ obeys the fusion rules if and only if for every $j$,
we have
$w_j\in \{-\La_j, -\La_j+2,...,\La_j-2, \La_j\}$ and $\La_j \leq a_j+a_{j-1}$.

If  $(a_1,...,a_n)$ obeys the fusion rules with respect to $\vec\La$,
then $(a_1+1,...,a_n+1)$ also obeys them. 

Let $w=(w_1,...,w_n)$ be an $n$-tuple of integers such that 
$w_1+...+w_n=0$ and  $w_j\in \{-\La_j, -\La_j+2,..., \La_j-2,
\La_j\}$
for all $j=1,...,n$. We call $w$  {\it a
weight vector}. The vector $-w$ will be called {\it the weight vector
dual to} $w$.

For a weight vector $w$, define an $n$-tuple $\Si(w)=(\Si^1,...,\Si^n)$
by the rule $\Si^j=\sum_{l=1}^jw_j$. Note that $\Si^n=0$.

It is clear that there exists a (non-negative) integer $k(w)$ 
such that for every integer $k$ the $n$-tuple
$(\Si^1+k,...,\Si^n+k)$ obeys the fusion rules with respect to
$\vec\La$ if and only if $k\geq k(w)$.
The integer $k(w)$ will be called {\it the shift number} associated with
a weight vector $w$. 
Every $n$-tuple obeying the fusion rules has the form
$\Si(w)+(k,...,k)$ for a suitable weight vector $w$ and a non-negative
integer $k$.

\subsection{Vanishing conditions and fusion rules}
Let $\vec\La=(\La_1,...,\La_n)$ be a vector of natural numbers,
$\La_1+...+\La_n=2m$, and let
$v(z,\la)$ be a hypergeometric solution of the qKZB equations with values
in $V_{\vec\La}[0]$ in the sense of Sec. \ref{finite}. Write the solution in coordinates,
$$
v(z,\la)\,=\, \sum_{M,\,{}|M|=m}\, v_M(z,\la)\,E_M
$$
where $\{E_M\}$ is the reduced basis in $V_{\vec\La}[0]$.
The following Theorem describes the resonance relations for admissible coordinates
of the hypergeometric solution.

Let $M=(m_1,...,m_n)$ be an admissible index,
$m_j \leq \La_j$ for all $j$.  Let $w_M= (\La_1-2m_1,...,
\La_n-2m_n)$ be the $\h$-weight vector of the basis vector
$E_M$ and $-w_M$ the weight vector dual to $w_M$. 
The weight vector $-w_M$ is the $\h$-weight of the
basis vector $E_{s(M)}$, $s(M)=(\La_1 - m_1,...,\La_n - m_n)$.
Let $k(w_M)$ and $k(-w_M)$ be the shift numbers defined in
Sec. \ref{fus}.

Consider the resonance relations of Corollary \ref{part-red} involving
the coefficient $v_M(z,\la)$. Since $\vec\La$ is a vector of natural
numbers, in any such a  relation we have $\la=2\eta k$ for a
suitable integer $k$.

\begin{thm}\label{rel-f.d}
${}$

Let $M=(m_1,...,m_n)$ be an admissible index.
\begin{enumerate}
\item[I.] If $k\,>\,k(-w_M)$ or $k\,<\,-k(w_M)$, then there is no
resonance relation of Corollary \ref{part-red} involving
$v_M(\,z,\,2\eta k\,)$.
\item[II.] If $k\,\in\, [-k(w_M), \,k(-w_M)]$, then the resonance relations
of Corollary \ref{part-red} imply the relation
\bean
v_M(\,z,\,2\eta k\,)\,=\,v_{s(M)}(\,z,\,-\,2\eta k\,).
\notag
\eean
\end{enumerate}
\end{thm}

Consider a hypergeometric solution $\tilde v (z,\la)$ of the qKZB equations
with values in $L_{\vec\La}[0]$ in the sense of Sec. \ref{finite}.
Any such a solution has the form $\tilde v = \pi \, v$, where $v$
is a hypergeometric solution of the qKZB equations with values in
$V_{\vec\La}[0]$ and  
$\pi: V_{\vec\La}[0]\to L_{\vec\La}[0]$ is the canonical projection.
Consider the anti-symmetrization,
$$
A\tilde v(z,\la)\,=\, \tilde v(z,\la)\,- S\tilde v(z,\la)\,=\, 
\tilde v(z,\la)\,- \,s_{\La_1}\otimes ...\otimes s_{\La_n}\,{}\,
\tilde v(z,-\la),
$$
of $\tilde v$ with respect to the Weyl reflection. By Theorem \ref{Rm-weyl},
$A\tilde v$ is a solution of the qKZB equations.
Write $A\tilde v$ in coordinates, 
$$
A\tilde v(z,\la)\,=\,\sum_{\text{adm}\,M,\,{}|M|=m}\,
A\tilde v_M(z,\la)\, E_M.
$$
\begin {corollary} {\bf (Vanishing Conditions)}\label{VC}
For any (admissible) index $M$ and any integer $k\in\, [-k(w_M), \,k(-w_M)]$, 
we have
$$
A\tilde v_M(\,z,\,2\eta k\,)\,=\,0\,.
$$
\end{corollary}

In other words, let $M=(m_1,...,m_n)$ be an admissible index.
Let $w=(w_1,...,w_n)$, $w_j=\La_j-2m_j$, be the $\h$-weight of the
basis vector $E_M$. Form the sums $\Sigma^j=\sum_{i=1}^j w_j$.
Let $k$ be a non-negative integer. Corollary \ref{VC} says
that $A\tilde v_M(\,z,\,2\eta k\,)\,=\,0\,$ unless the vector
$$
(-\Sigma^1 +k-1, -\Sigma^2 +k-1,...,-\Sigma^n +k-1)
$$
satisfies the $sl_2$ fusion rules with respect to $(\La_1,...,
\La_n)$ and
$A\tilde v_M(\,z,\,-\,2\eta k\,)\,=\,0\,$ unless the vector
$$
(\Sigma^1 +k-1, \Sigma^2 +k-1,...,\Sigma^n +k-1)
$$
satisfies the $sl_2$ fusion rules with respect to $(\La_1,...,
\La_n)$. In particular, $A\tilde v_M(\,z,\,0\,)\,=\,0\,$ since
$\Sigma^n=0$.

{\bf Remark.} Corollary \ref{VC} has an analog for the case when
$n=1$, $\vec \La = 2m$. Namely, let $u=u_{mm}(\la,\mu,\tau,p) e_m\otimes e_m
\in V_{\vec \La}[0]\otimes V_{\vec \La}[0]$ 
be the corresponding universal hypergeometric function.
Define the Weyl anti-symmetric function
$
Au(\la,\mu,\tau,p)=(u_{mm}(\la,\mu,\tau,p) - u_{mm}(-\la,\mu,\tau,p))
e_m\otimes e_m.
$
Then for all $k=-m,-m+1,...,m-1,m$ we have
\bean
Au(2\eta k,\mu,\tau,p)\,=\,0,
\notag
\eean
see Part III of Theorem \ref{v}.

\subsection{Proof of Theorem \ref{rel-f.d}}
\begin{lemma}
If $k\,>\,k(-w_M)$ or $k\,<\,-k(w_M)$, then there is no
resonance relation of Corollary \ref{part-red} involving
$v_M(\,z,\,2\eta k\,)$.
\end{lemma}
\begin{proof}
For a given integer $k$, a relation involving $v_M(\,z,\,2\eta k\,)$ exists,
if for some $j\in [1,n]$, there is $a\,\in\, [0,m_j+m_{j+1}]\,-\,\{m_j\}$
such that
\bean\label{cond}
\sum_{l=1}^{j-1}(\La_l-2m_l)\,+\,\La_j-m_j-a\,=\,k\,.
\eean

Let $k\,>\,k(-w_M)$. It means that for all $j=1,...,n$ we have
$$
(-\sum_{l=1}^{j-1}(\La_l-2m_l)+k)\,+\,
(-\sum_{l=1}^{j}(\La_l-2m_l)+k)\,>\,\La_j.
$$
Here $m_{n+1}=m_1$. This is equivalent to
\bean\label{below}
k\,>\,\sum_{l=1}^{j-1}(\La_l-2m_l)\,+\,\La_j - m_j\,.
\notag
\eean
This inequality contradicts \Ref{cond}.

Let $k\,<\,-\,k(w_M)$. It means that for all $j=1,...,n$ we have
$$
(\sum_{l=1}^{j-1}(\La_l-2m_l)-k)\,+\,
(\sum_{l=1}^{j}(\La_l-2m_l)-k)\,>\,\La_j.
$$
This is equivalent to
\bean\label{}
k\,<\,\sum_{l=1}^{j-2}(\La_l-2m_l)\,+\,\La_{j-1} - m_{j-1}-(m_{j-1}+m_{j})\,.
\notag
\eean
This inequality contradicts \Ref{cond}.
\end{proof}
\begin{lemma}\label{exist}
If $k\,\in\,[-k(w_M), k(-w_M)]$, then either there is a
resonance relation of Corollary \ref{part-red} involving
$v_M(\,z,\,2\eta k\,)$ or $k=0$, all $\La_j,\,j=1,...,n,$ are even, $M=(\La_1/2,...
,\La_n/2)$, $s(M)=M$ and we have the statement of Theorem \ref{rel-f.d},
$v_M(\,z,\,0\,)\,=\,v_{s(M)}(\,z,\,0\,)$.
\end{lemma}
\begin{proof} There are only four cases. We prove  the Lemma for each of them. 
\begin{enumerate}
\item[1.] For all $j \in [1,n]$,
\bean\label{1}
\sum_{l=1}^j (\La_l-2m_l)\,<\, k.
\eean
\item[2.] For all $j \in [1,n]$,
\bean\label{2}
\sum_{l=1}^j (\La_l-2m_l)\,>\, k.
\eean
\item[3.] There is $j \in [1,n]$, such that
\bean
\sum_{l=1}^{j-1} (\La_l-2m_l)\,= \, k.
\notag
\eean
\item[4.] There is $j \in [1,n]$, such that
\bean
\sum_{l=1}^{j-1} (\La_l-2m_l)\,> \, k\,>\,\sum_{l=1}^{j} (\La_l-2m_l)\,.
\notag
\eean
\end{enumerate}
In case 1, we have $0< k \leq k(-w_M)$. This means that
there is $j\in [1,n]$ such that
$$
(-\sum_{l=1}^{j-1}(\La_l-2m_l)+k)\,+\,
(-\sum_{l=1}^{j}(\La_l-2m_l)+k)\,\leq \,\La_j,
$$
which can be rewritten as
\bean\label{}
k\,\leq \,\sum_{l=1}^{j-1}(\La_l-2m_l)\,+\,\La_j - m_j\,.
\notag
\eean
Combining with \Ref{1} we get
$$
\sum_{l=1}^{j-1}(\La_l-2m_l)\,+\,\La_j - 2m_j\,<\,
k\,\leq \,\sum_{l=1}^{j-1}(\La_l-2m_l)\,+\,\La_j - m_j\,.
$$
Hence there is $a\in [0,m_j)$ such that
$$
\sum_{l=1}^{j-1}(\La_l-2m_l)\,+\,\La_j - m_j-a\,=\,k.
$$
Then there is a relation involving $v_M(\,z,\,2\eta k\,)$ corresponding to this equation.

In case 2, we have $-k(w_M)\leq k < 0$. This means that
there is $j\in [1,n]$ such that
$$
(\sum_{l=1}^{j-1}(\La_l-2m_l)-k)\,+\,
(\sum_{l=1}^{j}(\La_l-2m_l)-k)\,\leq \,\La_j,
$$
which can be rewritten as
\bean\label{}
k\,\geq \,\sum_{l=1}^{j-2}(\La_l-2m_l)\,+\,\La_{j-1} - m_{j-1}-(m_{j-1}+m_j)
\,.
\notag
\eean
Combining with \Ref{2} we get
$$
\sum_{l=1}^{j-2}(\La_l-2m_l)\,+\,\La_{j-1} - 2m_{j-1}\,>\,
k\,\geq \,\sum_{l=1}^{j-2}(\La_l-2m_l)\,+\,\La_{j-1} -  m_{j-1}-(m_{j-1}+m_j).
$$
Hence there is $a\in (m_{j-1},m_{j-1}+m_j]$ such that
$$
\sum_{l=1}^{j-2}(\La_l-2m_l)\,+\,\La_{j-1} - m_{j-1}-a\,=\,k.
$$
Then there is a relation involving $v_M(\,z,\,2\eta k\,)$ corresponding to this equation.

In case 3, either $k=0$ and $M=(\La_1/2,...,\La_n/2)$ or there is
$j\in [1,n]$ such that
$$
\sum_{l=1}^{j}(\La_l-2m_l)\,=\,k\,>\,\sum_{l=1}^{j+1}(\La_l-2m_l)\,.
$$
In the last case, we have $0\leq \La_{j+1}-m_{j+1}< m_{j+1}$ and a relation
of Corollary \ref{part-red} corresponding to the equation
$$
\sum_{l=1}^{j}(\La_l-2m_l)\,+\,\La_{j+1}-m_{j+1}-(\La_{j+1}-m_{j+1})\,=\,k\,.
$$

In case 4, there is $a\in [0,m_j)$ such that
$$
\sum_{l=1}^{j-2}(\La_l-2m_l)\,+\,\La_{j-1} - m_{j-1}-a\,=\,k\,,
$$
and then there is a relation
of Corollary \ref{part-red} corresponding to this equation.
\end{proof}

\begin{lemma}\label{right}
Let $k\,\in\,[-k(w_M), k(-w_M)]$ and for all $j \in [1,n]$,
\bean\label{i}
\sum_{l=1}^j (\La_l-2m_l)\,<\, k.
\eean
Then there is a sequence of relations of Corollary \ref{part-red}
implying $v_M(\,z,\,2\eta k\,)\,=\,v_{s(M)}(\,z,\,-\,2\eta k\,)$.
\end{lemma}
\begin{proof} By Lemma \ref{exist}, there are $j \in [1,n]$ and 
$a \in [0, m_j+m_{j+1}]-\{m_j\}$ such that
\bean\label{begin}
\sum_{l=1}^{j-1}(\La_l-2m_l)\,+\,\La_j - m_j-a\,=\,k.
\eean
Here $a\,=\,m_j\,+\,\sum_{l=1}^{j}(\La_l-2m_l)\,-\,k$.
This means that the transformation $T_j(k)$ can be applied to $v_M(z,2\eta k)$,
see Sec. \ref{transf}. We prove that the product of transformations
$T_{j-1}(-k)...T_2(-k)$ $ T_1(-k)$ $ T_n(k)$ $...$
$T_{j+1}(k)T_j(k)$ can be applied to $v_M(z,2\eta k)$
and gives $v_{s(M)}(z,\,-2\eta k)$.

We assume that $1\leq j <n$, the proof for $j=n$ is similar.

For all $i$, define $\Si^i$ by $\Si^i=\sum_{l=1}^i (\La_l-2m_l)$.
Equation \Ref{begin} implies a relation
\bean\label{1st}
v_M(\,z,\,2\eta k\,)\,=\,v_{(m_1,...,m_{j-1}, \, m_j\,+\,\Si^{j}\,-\,k, \, 
m_{j+1} \,-\,\Si^{j}\,+\,k, \, 
m_{j+2},...,m_n)}(\,z,\,\,2\eta k\,)\,.
\notag
\eean
We have an equation
\bean\label{ii}
{}
\\
\Si^{j-1}+\La_j - 2(m_j+\Si^{j}\,-k) 
+ \La_{j+1} - (m_{j+1}-\Si^{j}+k) - (\La_{j+1}-m_{j+1})\,=\,k\,,
\notag
\eean
and assumption \Ref{i} implies
\bean\label{i1}
(\La_{j+1}-m_{j+1})\,<\,m_{j+1}-\Si^{j}+k\,.
\eean
Formulas \Ref{ii} and \Ref{i1} imply a relation
\bean\label{2nd}
v_{(m_1,...,m_{j-1}, 
m_j\,+\,\Si^{j}\,-\,k, \, 
m_{j+1} \,-\,\Si^{j}\,+\,k, \, 
m_{j+2},...,m_n)}(\,z,\,\,2\eta k\,)\,=
\notag
\\
v_{(m_1,...,m_{j-1}, 
m_j+\Si^{j}-k, \, \La_{j+1}-m_{j+1},\,m_{j+2}-\Si^{j+1}+k, \, 
m_{j+3},...,m_n)}(\,z,\,\,2\eta k\,)\,.
\notag
\eean
Now we continue  the same way and get
\bean\label{3d}
v_M(\,z,\,\,2\eta k\,)=
v_{(m_1,...,m_{j-1}, 
m_j+\Si^{j}-k, \, \La_{j+1}-m_{j+1},...,
\La_{n-1}-m_{n-1},\,
m_n-\Si^{n-1}+k)}(\,z,\,\,2\eta k\,)\,.
\notag
\eean
We have an equation
\bean\label{3}
(m_n-\Si^{n-1}+k)-(\La_n-m_n)\,=\,k\, ,
\eean
and assumption \Ref{i} implies
\bean\label{i2}
(m_n-\Si^{n-1}+k)-(\La_n-m_n)\,>\,0\,.
\eean
Formulas \Ref{3} and \Ref{i2} imply a relation
\bean\label{4th}
v_M(\,z,\,\,2\eta k\,)=
v_{(m_1+k,\,m_2,...,m_{j-1}, 
m_j+\Si^{j}-k, \, \La_{j+1}-m_{j+1},...,
\,\La_{n}-m_{n})} (\,z,\,\,-\,2\eta k\,)\,.
\notag
\eean
We have an equation
\bean\label{4}
\La_1-(m_1+k)-(\La_1-m_1)\,=\,-\,k\,,
\eean
and assumption \Ref{i} implies
\bean\label{i3}
(m_1+k)-(\La_1-m_1)\,>\,0\,.
\eean
Formulas \Ref{4} and \Ref{i3} imply a relation
\bean\label{5th}
v_M(\,z,\,\,2\eta k\,)=
v_{(\La_1-m_1,\, m_2+k-\Si^1,\,m_3,...,m_{j-1}, 
m_j+\Si^{j}-k, \, \La_{j+1}-m_{j+1},...,
\,\La_{n}-m_{n})} (\,z,\,\,-\,2\eta k\,)\,.
\notag
\eean
Continuing this way, we get
\bean\label{6th}
v_M(\,z,\,\,2\eta k\,)=
v_{(\La_1-m_1,...,\La_{j-2}-m_{j-2},\,
m_{j-1}+k-\Si^{j-2},\,
m_j+\Si^{j}-k, \, \La_{j+1}-m_{j+1},...,
\,\La_{n}-m_{n})} (\,z,\,\,-\,2\eta k\,)\,,
\notag
\eean
and finally
\bean\label{6th}
v_M(\,z,\,\,2\eta k\,)\,=v_{s(M)}(\,z,\,\,-\,2\eta k\,)\,.
\notag
\eean
Thus, starting from a relation, corresponding to a reconstruction of
the pair of  $j$-th and $j+1$-st indices, then moving the total circle to
the right and  applying the relations, corresponding to reconstruction of pairs of indices
$(j+1, j+2)$,..., $(n,1)$,...,$(j-1,j)$, we proved the Lemma.
\end{proof}

\begin{lemma}\label{left}
Let $k\,\in\,[-k(w_M), k(-w_M)]$ and for all $j \in [1,n]$,
\bean\label{1r}
\sum_{l=1}^j (\La_l-2m_l)\,>\, k.
\eean
Then there is a sequence of relations of Corollary \ref{part-red}
implying $v_M(\,z,\,2\eta k\,)\,=\,v_{s(M)}(\,z,\,-\,2\eta k\,)$.
\end{lemma}
\begin{proof} The proof of this Lemma is similar to the proof of
Lemma \ref{right}, but now we move over the total circle to the left.

Namely, by Lemma \ref{exist}, there are $j \in [1,n]$ and 
$a \in [0, m_j+m_{j+1}]-\{m_j\}$ such that equation \Ref{begin}
is satisfied. 
This means that the transformation $T_j(k)$ can be applied to $v_M(z,2\eta k)$.

 We prove that the product of transformations
$T_{j+1}(-k)...T_{n-1}(-k)
$ $T_n (k) T_1(k)T_2(k)$ $...$ 
$ T_{j-1}(k)T_j(k)$ can be  applied to $v_M(z,2\eta k)$
and gives $v_{s(M)}(z,\,-2\eta k)$.

We assume that $1 < j \leq n$, the proof for $j=1$ is similar.

Equation \Ref{begin} implies a relation
\bean\label{}
v_M(\,z,\,2\eta k\,)\,=\,v_{(m_1,...,m_{j-1}, \, m_j\,+\,\Si^{j}\,-\,k, \, 
m_{j+1} \,-\,\Si^{j}\,+\,k, \, 
m_{j+2},...,m_n)}(\,z,\,\,2\eta k\,)\,.
\notag
\eean
We also have a relation
\bean\label{2nd}
v_{(m_1,...,m_{j-1}, 
m_j\,+\,\Si^{j}\,-\,k, \, 
m_{j+1} \,-\,\Si^{j}\,+\,k, \, 
m_{j+2},...,m_n)}(\,z,\,\,2\eta k\,)\,=
\\
v_{(m_1,...,m_{j-2}, m_{j-1}+\Si^{j-1}-k, \,\La_j-m_j,
m_{j+1} \,-\,\Si^{j}\,+\,k, \, 
m_{j+2},...,m_n)}(\,z,\,\,2\eta k\,)\,,
\notag
\eean
since we have equations
\bean
\Si^{j-2}+\La_{j-1}-m_{j-1}-(m_{j-1}+\Si^{j-1}-k)\,=\,k\,,
\notag
\\
m_{j-1} + (m_{j}+\Si^{j}-k)\,=\,
(m_{j-1}+\Si^{j-1}-k)+(\La_j - m_j)\,,
\notag
\eean
and an inequality
\bean
m_{j}+\Si^{j}-k\,>\,
\La_j - m_j\,
\notag
\eean
implied by assumption \Ref{1r}.

We continue  the same way and get
\bean\label{3d}
v_M(\,z,\,\,2\eta k\,)=
v_{(m_1+\Si^1-k,\La_2-m_2,...,\La_j-m_j,\,
m_{j+1} \,-\,\Si^{j}\,+\,k, \, 
m_{j+2},...,m_n)}(\,z,\,\,2\eta k\,)\,.
\notag
\eean
Similarly we have a relation
\bean\label{3d}
v_M(\,z,\,\,2\eta k\,)=
v_{(\La_1-m_1,...,\La_j-m_j,\,
m_{j+1} \,-\,\Si^{j}\,+\,k, \, 
m_{j+2},...,m_{n-1}, m_n-k)}(\,z,\,\,-\,2\eta k\,)\,,
\notag
\eean
then relations
\bean\label{3d}
v_M(\,z,\,\,2\eta k\,)=
v_{(\La_1-m_1,...,\La_j-m_j,\,
m_{j+1} \,-\,\Si^{j}\,+\,k, \, 
m_{j+2},...,m_{n-2}, m_{n-1}+\Si^{n-1}-k,\La_n-m_n)}(\,z,\,\,-\,2\eta k\,)\,,
\notag
\eean
\bean\label{3d}
v_M(\,z,\,\,2\eta k\,)=
v_{(\La_1-m_1,...,\La_j-m_j,\,
m_{j+1} \,-\,\Si^{j}\,+\,k, \, 
m_{j+2}+\Si^{j+2}-k,\La_{j+3}-m_{j+3},...,\La_n-m_n)}(\,z,\,\,-\,2\eta k\,)\,,
\notag
\eean
and finally
\bean\label{3d}
v_M(\,z,\,\,2\eta k\,)=
v_{s(M)}(\,z,\,\,-\,2\eta k\,)\,.
\notag
\eean
The Lemma is proved.
\end{proof}

In order to finish the proof of Theorem \ref{rel-f.d} we assume that
$$
\text{min}\{\Si^1,...,\Si^n\}\,\leq \,k \, \leq\,
\text{max}\{\Si^1,...,\Si^n\}
$$
and prove that \bean
v_M(\,z,\,2\eta k\,)\,=\,v_{s(M)}(\,z,\,-\,2\eta k\,).
\notag
\eean

Considering the resonance relations it is convenient to think
of the factors of the tensor product $V_{\vec\La}$ as positioned
in the cyclic order so that the factor $V_{\La_1}$ follows the factor
$V_{\La_n}$. Therefore, considering the index of a factor we shall consider
it modulo $n$ so that the $n+l$-th factor is the same as the $l$-th factor.
The sums $\Si^j$ introduced above correspond to this convention since
$\Si^n=0$.

We call $j\in [1,n]$ {\it a distinguished vertex} if $\Si^j=k$.

We define  top and bottom intervals.

{\it A  top interval} is a sequence of vertices
$r, r+1, ... , s$ such that $\Si^r, \Si^{r+1},..., \Si^s \,> \, k$ and
$\Si^{r-1}, \Si^{s+1} \leq k$.
Here $r$ can be equal to $s$.

{\it A  bottom interval} is a sequence of vertices
$r, r+1, ... , s$ such that $\Si^r, \Si^{r+1},..., \Si^s \,< \, k$ and
$\Si^{r-1}, \Si^{s+1} \geq k$.
Here $r$ can be equal to $s$.

In each of the two cases, the vertices $r, r+1,...,s$ are called the vertices
of the corresponding  interval.

Each not distinguished vertex is a vertex of a unique top or bottom
interval.

We define notions of the boundary and proper  intervals.
A top or bottom interval $r, r+1,...,s$
is called {\it the boundary interval} if $ r\leq n \leq s$.
All other top or bottom intervals are called {\it the proper  intervals}.

There is no  boundary interval if and only if $k=0$.

Our goal is to transform the coordinates of the index
$M=(m_1,...,m_n)$ into the coordinates of the index
$s(M)=(\La_1-m_1,...,\La_n-m_n)$. We construct a specially ordered product
of transformations $T_i$, corresponding to not distinguished
vertices, which can be applied to $v_M(z, 2\eta k)$ and makes the transformation.

First we associate a product $T^I$ of transformations $T_i$ to each 
top and bottom interval $I=\{r, r+1,...,s\}$.

Namely, let $I$ be a proper top interval,
then $T^I= T_r(k)T_{r+1}(k)...T_s(k)$.

Let $I$ be the boundary top interval, then
$T^I= T_r(-k)T_{r+1}(-k)...  T_{n-1}(-k)T_{n}(k)...T_s(k)$.

Let $I$ be a proper  bottom interval,
then $T^I= T_{s}(k)...T_{r+1}(k)T_{r}(k)$.

Let $I$ be the boundary bottom interval, then
$T^I= T_{s}(-k)... T_{n+1}(-k)T_n(k)...T_{r+1}(k)T_{r}(k)$.

We order all  top or bottom intervals as follows. If $I$ is boundary and
$J$ is proper, then $I<J$. If $I$ and $J$ are proper, then $I<J$ if there exist
$1\leq i < j \leq n$ such that $i$ is a vertex of $I$ and
$j$ is a vertex of $J$.

Let $\{I_1, I_2,..., I_r\}$ be the set of all  top or bottom intervals
and $I_1< I_2<...< I_r$. Define the product of transformations $T_i$ corresponding
to $v_M(z,2\eta k)$ by
$$
P\,=\, T^{I_1}\,T^{I_2}...T^{I_r}.
$$
\begin{lemma} 
The product $P$ of transformations $T_i$ can be applied to
$v_M(z,2\eta k)$ and transforms $v_M(z,2\eta k)$ to
$v_{s(M)}(z,\,-2\eta k)$.
\end{lemma}
The Lemma finishes the proof of Theorem \ref{rel-f.d}.

The proof of the Lemma is by direct verification similar to the 
proofs of Lemmas \ref{right} and \ref{left}.


\section{Resonance relations for Bethe ansatz eigenfunctions}\label{sba}

\subsection{Bethe ansatz eigenfunctions of commuting difference operators}
In this section we fix 
$\vec\La=(\Lambda_1,\dots,\Lambda_n)\in \C^n$ such that $\sum \Lambda_l
=2 m$ for some positive integer $m$.
 Let $K_i(z,\tau,p)$ be the qKZB operators acting on
the space ${\mathcal F} (V_{\vec\La}[0])$ of  meromorphic functions
of $\lambda\in\C$ with values in the zero weight space $V_{\vec\Lambda}[0]$.
The qKZB operators satisfy the compatibility conditions \Ref{compatib}.

There is a closely related set of commuting difference operators $H_j(z)$
acting on ${\mathcal F} (V_{\vec\La}[0])$ and the corresponding eigenvalue problem
\be
H_j(z)\psi=\epsilon_j\psi, \qquad j=1,\dots,n,\qquad \psi\in
{\mathcal F} (V_{\vec\La}[0])\,.
\ee
Here 
\begin{equation}\label{Go}
H_j(z)=R_{j,j-1}(z_j\!-\!z_{j-1})\cdots
R_{j,1}(z_j\!-\!z_{1})
\Gamma_j
R_{j,n}(z_j\!-\!z_n)\cdots,R_{j,j+1}(z_j\!-\!z_{j+1}),
\end{equation}
 $z=(z_1,\dots,z_n)$ is a fixed generic point in $\C^n$, $\psi$
is in ${\mathcal F} (V_{\vec\La}[0])$ and the operators $R_{j,k}$
are defined in Sec. \ref{qkzb}.  The fact that the operators $H_j(z)$ commute with each
other follows from the compatibility of the qKZB operators as $p\to 0$.

In \cite{FTV1} we gave a formula for common quasiperiodic eigenfunctions of the operators
$H_j$, i.e., functions $\psi$ such that $H_j\psi=\epsilon_j\psi$,
$j=1,\dots, n$, for some $\epsilon_j\in\C$. The quasiperiodicity
assumption means that we require that $\psi(\lambda+1)=
\mu\psi(\lambda)$ for some multiplier $\mu\in\C^\times$.

The formula is given in terms  of the
weight functions $\om_{j_1,...,j_n}(t_1,...,t_m,z,\la)$ associated with the tensor product
$V_{\vec\La}$ and defined in Sec. \ref{bases}.

\begin{thm}\label{tba}
 \cite{FTV1}
Let $c\in\C$.
Suppose that $t_1,\dots,t_m$ obey the system of
``Bethe ansatz'' equations
\begin{equation}\label{bae}
\prod_{l=1}^{n}
\frac{\theta(t_j-z_l+\eta\La_l)}
{\theta(t_j-z_l-\eta\La_l)}
\prod_{k:k\neq j}
\frac{\theta(t_j-t_k-2\eta)}{\theta(t_j-t_k+2\eta)}
=e^{-4\eta c},\qquad j=1,\dots,m.
\end{equation}
Then
\begin{equation}\label{psi}
\psi(\lambda)=\sum_{J,\,|J|=m}
e^{c\lambda}
\omega_{J}(t_1,\dots,t_m,z,\lambda)\,
e_J
\end{equation}
is a common eigenfunction of the commuting operators
$H_j$, $j=1,\dots,n$, with eigenvalues
\be
\epsilon_j=e^{-2c\eta\La_j}\prod_{k=1}^m
\frac{\theta(t_k-z_j-\eta\La_j)}
{\theta(t_k-z_j+\eta\La_j)}
\ee
and multiplier $\mu=(-1)^me^c$. Moreover, if $\,t_1,\dots,t_m$ are a
solution of
\Ref{bae}, and $\sigma\in S_m$ is any permutation,
then $t_{\sigma(1)},\dots,t_{\sigma(m)}$ are also a solution.
The eigenfunctions $\psi$ corresponding to these two solutions
are proportional.
\end{thm}

\subsection{Resonance relations }\label{8.2}
Let $\{E_J\}$ be the reduced basis in $V_{\vec\La}[0]$.
Let $\psi(\la)$ be an eigenfunction defined in Theorem \ref{tba}
and corresponding to a solution $(t_1,\dots,t_m)$ of the Bethe ansatz equations. Write
the eigenfunction in the reduced bases,
\begin{equation}\label{psi-reduced}
\psi(\lambda)\,=\,\sum_{J}\,
\psi_{J}(t_1,\dots,t_m,z,\lambda)\,
E_J\,,
\notag
\end{equation}
where the reduced coefficients $\psi_{J}(t_1,\dots,t_m,\lambda)$ are defined
by
$$
\psi_{J}(t_1,\dots,t_m,z,\lambda)\,=\,[j_1]!...[j_n]!\,
e^{c\lambda} \, \omega_{J}(t_1,\dots,t_m,z,\lambda)\,.
$$

\begin{thm}\label{relat}
${}$

\begin{enumerate}
\item[I.] Let $n>1$. Then
the reduced coefficients of an eigenfunction $\psi(\la)$ obey
the resonance relations
\bean\label{res-eig-1}
\psi_{M}(t_1,\dots,t_m,z,2\eta (\La_j - a-b\,+\,\sum_{l=1}^{j-1}(\La_l-2m_l)))=
\\
\psi_{L}(t_1,\dots,t_m,z,2\eta (\La_j - a-b\,+\,\sum_{l=1}^{j-1}(\La_l-2m_l)))
\notag
\eean
for  $j, a, b, k, L, M$ defined in Part I of Theorem \ref{v} and the resonance
relations
\bean\label{res-eig-2}
\psi_{M}(t_1,\dots,t_m,z,2\eta (a-b))\,=\,
\psi_{L}(t_1,\dots,t_m,z,2\eta (b-a))
\eean
for $a, b, k, L, M$ defined in Part II of Theorem \ref{v}.
\item[II.] 
Let $n=1$, $\vec \Lambda = 2m$. In this case
$\psi=\psi_{m} E_m$ and $\psi$ does not depend on $z$. We claim that
\bean\label{}
\psi_{m}(t_1,...,t_m, 2\eta a,\tau)\,=\,
\psi_{m}(t_1,...,t_m,-2\eta a,\tau)\,
\eean
for $a=1,...,m$.
\end{enumerate}
\end{thm}

\subsection{Proof of Theorem \ref{relat}}
The proof of Part II is similar to the proof of Part III of Theorem \ref{v}.

The resonance relations in formula \Ref{res-eig-1}
follow from Theorem \ref{w.f.rel}. The proof of the resonance relations
in formula \Ref{res-eig-2} is similar to the proof of Part II of Theorem
\ref{v}.

Namely, let $S_n$ act on $\C^n$ by
permutations of the coordinates, and let $s_j$ be the
transposition of the $j$-th and $(j+1)$-st coordinates.

The equations \Ref{bae} do
not depend on the ordering of the parameters $z_i,\Lambda_i$. Let
us fix a solution $t^*$ of \Ref{bae} with parameters
${\vec\Lambda}=(\Lambda_1,\dots,\Lambda_n)$
and ${z}=(z_1,\dots,z_n)$. Then
for each permutation $\sigma\in S_n$, the point $t^*$ is still a
solution of the equations \Ref{bae} with parameters $\sigma{\vec\Lambda}$,
$\sigma{z}$ and we have a corresponding function
$\psi_\sigma(\lambda)\in {\mathcal F} (V_{\sigma \vec\La}[0])$
given by the formula \Ref{psi} with
parameters $\sigma{\vec\Lambda}$, $\sigma{z}$. 

It follows from the definition of R-matrices that
\be
\psi_{s_j}(\lambda)=P^{(j,j+1)}
R_{\Lambda_j\Lambda_{j+1}}
^{(j,j+1)}(z_j-z_{j+1},\lambda-2\eta(h^{(1)}+\cdots+h^{(j-1)}))\psi(\lambda).
\ee
Introduce a permutation $s^\vee\in S_n$ by $s^\vee =s_1 s_2 ... s_{n-1}$.
Then $s^\vee\vec \La=\vec\La^\vee=(\La_n,\La_1,...,\La_{n-1})$.
Introduce  an operator $\Dl : V_{\vec \La}\to V_{\vec \La^\vee}$ by formula
\Ref{delta}.

\begin{lemma}
$$
\Dl\,\psi(\lambda)\,=\,\epsilon_1\,
\psi_{s^\vee}(\lambda).
$$
\end{lemma}
The proof of the Lemma is similar to the argument in the proof of Theorem
\ref{Transf}.

The Lemma finishes the proof of Theorem \ref{relat} since the resonance relations
for the eigenfunction $\psi(\la)$ in formula \Ref{res-eig-2}  
are the resonance relations
for the eigenfunction $\psi_{s^\vee}(\la)$ in formula \Ref{res-eig-1}
for $j=1$.

\subsection{ The fusion rules and the
resonance relations for the Bethe ansatz eigenfunctions
with values in finite dimensional modules }\label{8.4}

Let $\vec\La=(\La_1,...,\La_n)$ be a vector of natural numbers,
$\La_1+...+\La_n=2m$, and let
$\psi(\la)$ be a Bethe ansatz eigenfunction with values in 
$V_{\vec\Lambda}[0]$ and corresponding to a solution 
of the Bethe ansatz equations. Write the eigenfunction in coordinates,
$$
\psi(\la)\,=\, \sum_{M,\,{}|M|=m}\, \psi_M(\la)\,E_M
$$
where $\{E_M\}$ is the reduced basis in $V_{\vec\La}[0]$.
The following Theorem describes the resonance relations for admissible coordinates
of the eigenfunction.     

Let $M=(m_1,...,m_n)$ be an admissible index,
$m_j \leq \La_j$ for all $j$.  Let $w_M= (\La_1-2m_1,...,
\La_n-2m_n)$ be the $\h$-weight vector of the basis vector
$E_M$ and $-w_M$ the weight vector dual to $w_M$. 
Let $k(w_M)$ and $k(-w_M)$ be the shift numbers defined in
Sec. \ref{fus}.

Consider the resonance relations of Theorem \ref{relat} involving
the coefficient $\psi_M(\la)$. Since $\vec\La$ is a vector of natural
numbers, in any such a  relation we have $\la=2\eta k$ for a
suitable integer $k$.

\begin{thm}\label{bethe-f.d}
${}$

Let $M=(m_1,...,m_n)$ be an admissible index.
\begin{enumerate}
\item[I.] If $k\,>\,k(-w_M)$ or $k\,<\,-k(w_M)$, then there is no
resonance relation of Theorem \ref{relat} involving
$\psi_M(2\eta k)$.
\item[II.] If $k\,\in\, [-k(w_M), \,k(-w_M)]$, then the resonance relations
of Theorem \ref{relat} imply the relation
\bean
\psi_M(2\eta k)\,=\,\psi_{s(M)}(\,-\,2\eta k).
\notag
\eean
\end{enumerate}
\end{thm}

The proof of this Theorem coincides with the proof of Theorem \ref{rel-f.d}.

Let $L_{\vec\La}=L_{\La_1}\otimes ... \otimes L_{\La_n}$ be the corresponding 
product of the finite dimensional spaces and $\pi : V_{\vec\La} \to
L_{\vec\La}$ the canonical projection. The operators $H_j(z)$ induce on
$L_{\vec\La}[0]$ a system of commuting difference operators which also will be 
denoted $H_j(z)$. 

Let $\psi(\la)$ be  the Bethe ansatz eigenfunction  with values in 
$V_{\vec\Lambda}[0]$ and corresponding to a solution 
of the Bethe ansatz equations. 
Then the function
$$
\tilde\psi(\la)\,=\, \pi \psi (\la)\,=\,
\sum_{\text{adm}\,M,\,{}|M|=m}\, \psi_M(\la)\,E_M
$$
is an eigenfunction of the commuting operators $H_j(z)$ acting on
${\mathcal F} (L_{\vec\La}[0])$, $H_j(z)\tilde\psi(\la)= \epsilon_j\tilde
\psi(\la)$, where the eigenvalues $\epsilon_j$ are defined in Theorem 
\ref{tba}.
The values of $\tilde\psi(\la)$ obey the resonance relations
of Theorem \ref{bethe-f.d}.

Let $S\, :\, {\mathcal F} (L_{\vec\La}[0])\,\to\,
{\mathcal F} (L_{\vec\La}[0])$
be the Weyl reflection.
The Weyl reflection commutes with the operators $H_j(z)$,
$S\,H_j(z)\,=\,H_j(z)\,S$. Hence $S\tilde \psi(\la)$ is also an eigenfunction
of $H_j(z)$ with the same eigenvalues.

Consider the anti-symmetrization,
$$
A\tilde \psi(\la)\,=\, \tilde \psi (\la)\,- S\tilde \psi (\la)\,=\, 
\tilde \psi(\la)\,- \,s_{\La_1}\otimes ...\otimes s_{\La_n}\,{}\,
\tilde \psi (-\la),
$$
of $\tilde \psi(\la)$ with respect to the Weyl reflection. 
Write $A\tilde \psi (\la)$ in coordinates, 
$$
A\tilde \psi(\la)\,=\,\sum_{\text{adm}\,M,\,{}|M|=m}\,
A\tilde \psi_M(\la)\, E_M.
$$
\begin {corollary}\label{vc} {\bf (Vanishing Conditions)}
For any (admissible) index $M$ and any integer $k\in\, [-k(w_M), \,k(-w_M)]$, 
we have
$$
A\tilde \psi_M(2\eta k)\,=\,0\,.
$$
\end{corollary}

In other words, let $M=(m_1,...,m_n)$ be an admissible index and
 $w=(w_1,...,w_n)$, $w_j=\La_j-2m_j$,  the $\h$-weight of the
basis vector $E_M$. Form the sums $\Sigma^j=\sum_{i=1}^j w_j$.
Let $k$ be a non-negative integer. Corollary \ref{vc} says
that $A\tilde \psi_M(2\eta k\,)\,=\,0\,$ unless the vector
$$
(-\Sigma^1 +k-1, -\Sigma^2 +k-1,...,-\Sigma^n +k-1)
$$
satisfies the $sl_2$ fusion rules with respect to $(\La_1,...,
\La_n)$ and
$A\tilde \psi_M(-\,2\eta k\,)\,=\,0\,$ unless the vector
$$
(\Sigma^1 +k-1, \Sigma^2 +k-1,...,\Sigma^n +k-1)
$$
satisfies the $sl_2$ fusion rules with respect to $(\La_1,...,
\La_n)$. In particular, $A\tilde \psi_M(\,0\,)\,=\,0\,$ since
$\Sigma^n=0$.

{\bf Remark.} The Corollary \ref{vc} has an analog for the case when
$n=1$, $\vec \La = 2m$. Namely, let $\psi(\la)=\psi_{m}(\la) \,e_m$ 
be  the Bethe ansatz eigenfunction  with values in 
$V_{\vec\Lambda}[0]$ and corresponding to a solution 
of the Bethe ansatz equations. Consider the Weyl anti-symmetric eigenfunction
$A\psi(\la)=(\psi_m(\la)-\psi_m(-\la))\, e_m$. Then for $k=-m, -m+1,...,m-1, m$
we have 
$$
A \psi(2\eta k)\,=\,0\,.
$$

\subsection{ The case $\eta=1/2N$ and $e^{2c}=1$; 
the resonance relations for the Bethe ansatz eigenfunctions
and the fusion rules for the quantum group
$U_{e^{2\pi i/N}}(sl_2)$}\label{8.5}

In this section we assume that the parameter $\eta$ has the form
$\eta\,=\, {1\over 2N}, $
where $N$ is a natural number, and assume that
the parameter $c$ in the Bethe ansatz equations
obeys the relation $e^{2c}\,=\, 1,$ 
see Theorem \ref{tba}. 

Let $\vec\La=(\La_1,...,\La_n)$ be a vector of natural numbers,
$\La_1+...+\La_n=2m$, and let
$$
\tilde\psi(\la)\,=\,\sum_{\text{adm} \, M,\,{}|M|=m}\, \tilde\psi_M(\la)\,E_M
$$
be the Bethe ansatz eigenfunction with values in 
$L_{\vec\Lambda}[0]$ corresponding to a solution 
of the Bethe ansatz equations with parameter $c\in \C$.
Let $\tilde\psi_M(2\eta k)\,=\,\tilde\psi_{s(M)}(\,-\,2\eta k)$ for
$k\in [-k(w_M), k(-w_M)]$ be a resonance relation
of Theorem \ref{bethe-f.d} for the coefficients of the eigenfunction.
Since the function $\tilde\psi(\la)$ is quasiperiodic,
$\tilde\psi(\la+1)=(-1)^me^c\tilde\psi(\la)$, and since $2\eta = 1/N$,
the resonance relation implies a relation
$e^c\tilde\psi_M(2\eta (k+N))\,=\, e^{-c}\tilde\psi_{s(M)}(\,-\,2\eta (k+N))$.
If $e^{2c}=1$, then  a resonance relation of Theorem \ref{bethe-f.d},
 $\tilde\psi_M(2\eta k)\,=\,\tilde\psi_{s(M)}(\,-\,2\eta k)$, implies a series of 
new resonance relations,
\bean
\tilde\psi_M(2\eta (k+lN))\,=\,\tilde\psi_{s(M)}(\,-\,2\eta (k+lN)),\qquad \text{for}
\qquad l\in \Z\,.
\notag
\eean
Therefore, under the above assumptions on $\eta$ and $c$, the resonance
relations become  $N$-periodic with respect to $k$ and to describe 
all the resonance relations it is enough to describe the relations for $k\in [0,N)$.

\begin{corollary}\label{1/2N}
Under assumptions $\eta=1/2N$ and $e^{2c}=1$, for all $k\in [0,N-1]$, we have
$$
\tilde\psi_M(2\eta k)\,=\,\tilde\psi_{s(M)}(\,-\,2\eta k)
\qquad \text{and} \qquad
A\tilde\psi_M(2\eta k)\,=\,0
$$
unless
\bean\label{rel-quant} 
k(-w_M)\, <\,k \,<\, N\,-\, k(w_M).
\eean

\end{corollary}

The  relation \Ref{rel-quant} 
can be interpreted in terms of the fusion rules for the quantum group $U_q(sl_2)$
with $q=e^{2\pi i/N}$.

First we remind the fusion rules for $U_q(sl_2)$ with $q=e^{2\pi i / N}.$

For an  integer  $a\in [0,N-2]$,
let $L_a$ denote the $a+1-$dimensional
$U_q(sl_2)$-module with highest weight $q^{a}$. 
We say that integers
$a, b, c$ obey {\it the fusion rules for $U_q(sl_2),\, q^{2\pi i/N},$}
if $a, b, c \in [0, N-2]$ and 
the $U_{q}(sl_2)$-module $L_a\otimes L_b$ contains $L_c$ as a submodule.

It means that $a-b \in \{-c, -c+2,...,c-2, c\}$ and $c\leq  a+b\,\leq \,
 2N - c -4$.
If the triple $a, b, c$ obeys the fusion rules, then each of its its permutation does.

Let $\vec\La=(\La_1,...,\La_n)$,  $a=(a_1,...,a_n)$ be  $n$-tuples
 of integers. We say that $a$ obeys
{\it the fusion rules for $U_{q}(sl_2)$
with respect to $\vec\La$}, if 
 for every $j=1,...,n$, the triple $a_{j-1}, a_{j}, \La_j$ obeys
the fusion rules for $U_{q}(sl_2)$. 
This means that for $j=1,...,n$ we have
$a_{j-1}, a_j, \La_j \in [0, N-2]$ and
$\La_j\,\leq \, a_{j-1} + a_j \,\leq 2N - \La_j - 4$.
Here we assume that $a_{0}=a_n$.

\begin{lemma}

An integer $k$ satisfies \Ref{rel-quant}, if and only if
the vector
\bean\label{vect}
(-\Si^1 + k -1, -\Si^2+ k -1, ... , -\Si^n + k -1)
\eean
obeys the fusion rules for $U_q(sl_2)$ with respect to $\vec\La =
(\La_1, ... , \La_n).$
\end{lemma}
According to the Lemma, if $k$ satisfies \Ref{rel-quant}, then
$-\Si^j + k - 1, \La_j$ belong to  $[0, N-2]$ for all $j$. The Theorem also
implies that if there is
$j$ such that $\La_j > N-2$, then we have
$$
\tilde\psi_M(2\eta k)\,=\,\tilde\psi_{s(M)}(\,-\,2\eta k)
\qquad \text{and} \qquad
A\tilde\psi_M(2\eta k)\,=\,0
$$
for all $M$ and $k$.

\begin{proof} First we prove that \Ref{rel-quant} implies the $U_q(sl_2)$ fusion rules.
Indeed, the inequality $k(-w_M)<k$ implies
$
\La_j<-\Si^{j-1}+k - \Si^j +k
$
for all $j$.
Similarly, the inequality $k<N-k(w_M)$ implies
$
-\Si^{j-1}+k - \Si^j +k < 2N-\La_j
$
for all $j$. Hence,
$$
\La_j<-\Si^{j-1}+k - \Si^j +k < 2N - \La_j.
$$
Notice that the number
$ -\Si^{j-1}+k - \Si^j +k $ has the same parity as $\La_j$. 
Therefore, $\La_j \in [0, N-2]$ and
$$
\La_j \leq \Si^{j-1}+k-1 - \Si^j +k-1 \leq 2N - \La_j -4.
$$

The vector
$(-\Si^1 + k(-w_M), -\Si^2+ k(-w_M), ... , -\Si^n + k(-w_M) )$ obeys 
the fusion rules for $sl_2$. Hence, the inequality $k(-w_M)<k$ implies
$-\Si^j+k > 0$ for all $j$.

The vector
$(\Si^1 + k(w_M), \Si^2+ k(w_M), ... , \Si^n + k(w_M) )$ obeys 
the fusion rules for $sl_2$. Hence, the inequality $k(w_M)<N-k$ implies
$N>-\Si^j+k$ for all $j$. Thus $-\Si^j+k-1\in[0, N-2]$. Thus the vector
\Ref{vect} obeys the $U_q(sl_2)$ fusion rules. The converse implication is 
proved similarly.

\end{proof}

\subsection{Resonance Relations for Eigenfunctions of Transfer Matrices}\label{t-mat}
In this section we discuss the resonance relations for the eigenfunctions of the
transfer matrix of highest weight representations of the elliptic quantum group
$E_{\tau,\eta}(sl_2)$, see \cite{FV2, FTV1}.

Let $W$ be a representation of $E_{\tau,\eta}(sl_2)$ with L-operator
 $L(z,\la)\in \End\,(\C^2\otimes W)$.
Let $e_0, e_1$ be the standard basis of $\C^2$ and
 $E_{ij}$  the linear operators on $\C^2$ defined by $E_{ij}e_k= \dl_{jk} e_i$.
Writing $L= E_{00}\otimes a +  E_{01}\otimes b + E_{10}\otimes c +
 E_{11}\otimes d$
we get four operators, $a(z,\lambda)$, $b(z,\lambda)$, $c(z,\lambda)$,
$d(z,\lambda)$, the matrix elements of the L-operator, acting
on $W$, and obeying the various relations of $E_{\tau,\eta}(sl_2)$.
The transfer matrix $T_W(z)\in\End \,({\mathcal F} (W[0]))$ acts on functions
by
\be
T_W(z)f(\lambda)=a(z,\lambda)f(\lambda-2\eta)
+d(z,\lambda)f(\lambda+2\eta).
\ee
The relations imply that $T_W(z)T_W(w)=T_W(w)T_W(z)$ for all $z,w\in\C$.

In \cite{FV2}, common eigenfunctions of $T_W(z)$, $z\in\C$,
are constructed in the form
\be
b(t_1)\cdots b(t_m)v_c
\ee
where
\be
v_c(\lambda)=e^{c\lambda}\prod_{j=1}^m\frac
{\theta(\lambda-2\eta j)}
{\theta(2\eta)}v^0,
\ee
$v^0$ is the highest weight vector of $W$ 
and $b(t)$ is the difference operator
$(b(t)f)(\lambda)=b(t,\lambda)f(\lambda+2\eta)$, $f\in {\mathcal F} (W)$
(both the transfer matrix and the difference operators $b(t)$
are part of the {\em operator algebra} of the elliptic
quantum group, see \cite{FV1, FV2}).
The variables $t_1,\dots,t_m$ obey a set of Bethe ansatz equations,
which are up to a shift the same as the ones in \Ref{bae}.

The eigenfunctions of the transfer matrix in a tensor product of Verma modules
are described by the following two Theorems.

\begin{thm}\label{previous}
\cite{FTV1}
Let $ V_{\vec\La}=V_{\Lambda_1}(z_1)\otimes\cdots\otimes V_{\Lambda_n}(z_n)$ be
a tensor product of $E_{\tau,\eta}(sl_2)$
evaluation Verma modules with generic evaluation
points $z_1,\dots, z_n$ and let
\be
v_c(\lambda)=e^{c\lambda}\prod_{j=1}^m
\frac{\theta(\lambda-2\eta j)}{\theta(2\eta)}\,e_0\otimes ...\otimes e_0\,\in\, 
{\mathcal F} (V_{\vec\La}).
\ee
Then
\bea
\prod_{j=1}^mb(t_j+\eta)\,v_c
&=&
e^{c(\lambda+2\eta m)}
(-1)^m\prod_{i<j}\frac
{\theta(t_i-t_j+2\eta)}
{\theta(t_i-t_j)}\\
& &\times
\sum_{j_1+\cdots+j_n=m}\tilde\omega_{j_1,\dots,j_n}(t_1,\dots,t_n,\lambda,z,\tau)
\,e_{j_1}\otimes\cdots\otimes e_{j_n},
\eea
where $\tilde\omega_{j_1,\dots,j_n}(t_1,\dots,t_n,\lambda,z, \tau)$
are the mirror weight functions associated with the tensor product $V_{\vec\La}$.

\end{thm}

\begin{thm}\label{tba2}
\cite{FV2} 
Let $V_{\vec\La}$ and $v_c$ be as in Theorem \ref{previous}, and
let $T_{V_{\vec\La}}(w)$ be the corresponding transfer matrix acting
on functions in ${\mathcal F} (V_{\vec\La}[0])$.
Then for any solution $(t_1,\dots,t_m)$ of the Bethe ansatz
equations
\bean\label{bae2}
\prod_{j:j\neq i}
\frac
{\theta(t_j-t_i-2\eta)}
{\theta(t_j-t_i+2\eta)}
\prod_{k=1}^n
\frac
{\theta(t_i-z_k-(1+\Lambda_k)\eta)}
{\theta(t_i-z_k-(1-\Lambda_k)\eta)}
=
e^{4\eta c},\qquad i=1,\dots,m,
\eean
such that, for all $i<j$, $t_i\neq t_j\mod \Z+\tau\Z$,
the vector $\psi=b(t_1)\cdots b(t_m)v_c\in{\mathcal F} (V_{\vec\La}[0])$
is a common eigenfunction
of all transfer matrices $T_{ V_{\vec\La}}(w)$ with eigenvalues
\be
\epsilon(w)=
e^{-2\eta c}\prod_{j=1}^m
\frac
{\theta(t_j-w-2\eta)}
{\theta(t_j-w)}
+
e^{2\eta c}\prod_{j=1}^m
\frac
{\theta(t_j-w+2\eta)}
{\theta(t_j-w)}
\prod_{k=1}^n
\frac
{\theta(w-z_k-(1-\Lambda_k)\eta)}
{\theta(w-z_k-(1+\Lambda_k)\eta)}.
\ee
Moreover $\psi(\lambda+1)=(-1)^me^c\psi(\lambda)$.
\end{thm}

Now we can reformulate the results of Sections \ref{8.2}, \ref{8.4}, \ref{8.5}
for eigenfunctions of the transfer matrices.

Let $\{E_J\}$ be the reduced basis in $V_{\vec\La}[0]$.
Let $\psi(\la)$ be the eigenfunction of the transfer matrix
corresponding to a solution $(t_1,\dots,t_m)$ of the Bethe ansatz equations
\Ref{bae2}. Write
the eigenfunction in the reduced bases, $
\psi(\lambda)\,=\,\sum_{J}\,
\psi_{J}(t_1,\dots,t_m,z,\lambda)\,
E_J\,$.

\begin{thm}\label{relat-2}
${}$

\begin{enumerate}
\item[I.] Let $n>1$, $1<j\leq n$. Let
 $m_1,...,m_{j-2},k,m_{j+1},...,m_n$
be non-negative integers such that $m_1+...+m_{j-2}+k+m_{j+1}+...+m_n=m$. Let 
$a,b$ be integers such that $a\neq b,\, 0\leq a,b \leq k$. Let
$M=(m_1,...,m_{j-2},k-a,a,m_{j+1},...,m_n)$ and
$L=(m_1,...,m_{j-2},k-b,b,m_{j+1},...,m_n)$. Then 
\bean\label{v1-2}
\psi_{M}(t_1,...,t_m, z,\,2\eta (\La_j - a-b\,+\,\sum_{l=j+1}^{n}(\La_l-2m_l)))
\,=\,
\\
\psi_{L}(t_1,...,t_m, z,\,2\eta (\La_j - a-b\,+\,\sum_{l=j+1}^{n}(\La_l-2m_l))).
\notag
\eean
\item[II.] 
Let $n>1$. Let  $k, m_2,...,m_n$
be non-negative integers such that $k+m_2+...+m_n=m$. Let 
$a,b$ be such that $a\neq b,\, 0\leq a,b \leq k$. Let
$M=(a,m_2,...,m_{n-1},k-a)$ and
$L=(b,m_2,...,m_{n-1},k-b)$. Then 
\bean\label{v2-2}
\psi_{M}(t_1,...,t_m, z,\,2\eta (a-b))
\,=\,
\\
\psi_{L}(t_1,...,t_m, z,\,2\eta (b-a)).
\notag
\eean
\item[III.] Let $n=1$, $\vec \Lambda = 2m$. Then
$\psi=\psi_{m} E_m $ and for all
 $a = 1,...,m$ we have
\bean\label{}
\psi(t_1,...,t_m, z,\, 2\eta a)\,=\,
\psi(t_1,...,t_m, z,\,-\,2\eta a)\,.
\eean
\end{enumerate}
\end{thm}
Notice that the resonance relations of this Theorem differ from the resonance
relations of Theorem \ref{relat} because the coordinates of the eigenfunction in this Theorem
are mirror weight functions while the coordinates of the eigenfunction in Theorem \ref{relat}
are ordinary weight functions.

Now assume that $\vec\La=(\La_1,...,\La_n)$ is a vector of natural numbers,
$\La_1+...+\La_n=2m$. Let $L_{\vec\La}=L_{\La_1}\otimes ... 
\otimes L_{\La_n}$ be the corresponding 
product of the finite dimensional spaces and $\pi : V_{\vec\La} \to
L_{\vec\La}$ the canonical projection. 

Let $\psi(\la)$ be  the Bethe ansatz eigenfunction  with values in 
$V_{\vec\Lambda}[0]$ and corresponding to a solution 
of the Bethe ansatz equations \Ref{bae2}.
Then the function
$$
\tilde\psi(\la)\,=\, \pi \psi (\la)\,=\,
\sum_{\text{adm}\,M,\,{}|M|=m}\, \psi_M(\la)\,E_M
$$
is an eigenfunction of the transfer matrix $T_{L_{\vec\La}}(w)$
acting on ${\mathcal F} (L_{\vec\La}[0])$, $ T_{L_{\vec\La}}(w)\tilde\psi(\la)= 
\epsilon(w)\tilde \psi(\la)$, where the eigenvalue $\epsilon(w)$ are defined in Theorem 
\ref{tba2}.

Let $S\, :\, {\mathcal F} (L_{\vec\La}[0])\,\to\,
{\mathcal F} (L_{\vec\La}[0])$ be the Weyl reflection.
The Weyl reflection commutes with  the transfer matrix $T_{L_{\vec\La}}(w)$.
Hence $S\tilde \psi(\la)$ is also an eigenfunction
of the transfer matrix  with the same eigenvalue.

Consider the anti-symmetrization,
$$
A\tilde \psi(\la)\,=\, \tilde \psi (\la)\,- S\tilde \psi (\la)\,=\, 
\tilde \psi(\la)\,- \,s_{\La_1}\otimes ...\otimes s_{\La_n}\,{}\,
\tilde \psi (-\la).
$$
Write $A\tilde \psi (\la)$ in coordinates, 
$$
A\tilde \psi(\la)\,=\,\sum_{\text{adm}\,M,\,{}|M|=m}\,
A\tilde \psi_M(\la)\, E_M.
$$
The vanishing conditions for this anti-symmetric eigenfunction are described by
the $sl_2$ modified fusion rules. 

Recall that if $(a_1,...,a_n)$ is a vector of integers,
then we say that $(a_1,...,a_n)$ satisfies the $sl_2$ fusion rules with respect to
$(\La_1,...,\La_n)$ if all $a_j$ are non-negative and for every $j$ the triple
$a_{j-1}, a_j, \La_j$ obeys the $sl_2$ fusion rules (where $a_0=a_n$).
We shall say that $(a_1,...,a_n)$ satisfies {\it the 
$sl_2$ modified fusion rules with respect to
$(\La_1,...,\La_n)$} if all $a_j$ are non-negative and for every $j$ the triple
$a_{j}, a_{j+1}, \La_j$ obeys the $sl_2$ fusion rules (where $a_{n+1}=a_1$).
Notice that if $\La_1=\La_2=...=\La_n$, then the two fusion rules coincide.

\begin {thm}\label{vc2} {\bf (Vanishing Conditions)}
${}$

Let $M=(m_1,...,m_n)$ be an admissible index and
 $w=(w_1,...,w_n)$, $w_j=\La_j-2m_j$,  the $\h$-weight of the
basis vector $E_M$. Form the sums $\tilde\Sigma^j=\sum_{i=j}^n w_j$.
Let $k$ be a non-negative integer. Then
$A\tilde \psi_M(2\eta k\,)\,=\,0\,$ unless the vector
$$
(-\tilde\Sigma^1 +k-1, -\tilde\Sigma^2 +k-1,...,-\tilde\Sigma^n +k-1)
$$
satisfies the $sl_2$ modified fusion rules with respect to $(\La_1,...,
\La_n)$ and
$A\tilde \psi_M(-\,2\eta k\,)\,=\,0\,$ unless the vector
$$
(\tilde\Sigma^1 +k-1, \tilde\Sigma^2 +k-1,...,\tilde\Sigma^n +k-1)
$$
satisfies the $sl_2$ modified fusion rules with respect to $(\La_1,...,
\La_n)$. In particular, $A\tilde \psi_M(\,0\,)\,=\,0\,$ since
$\tilde\Sigma^1=0$.

\end{thm}

The fact that this Theorem involves the modified fusion rules while
the corresponding Corollary \ref{vc} uses the ordinary fusion rules is
due to the fact that the coordinates of the eigenfunction in this Theorem
are given in terms of mirror weight functions while the coordinates of the
eigenfunction in Corollary \ref{vc} are given in terms of ordinary weight functions.

{\bf Remark.} If $n=1$, $\vec \La = 2m$ and $\psi(\la)=\psi_{m}(\la) \,e_m$ 
is  the Bethe ansatz eigenfunction of the transfer matrix with values in 
$V_{\vec\Lambda}[0]$. Then the corresponding Weyl anti-symmetric eigenfunction
$A\psi(\la)=(\psi_m(\la)-\psi_m(-\la))\, e_m$ satisfies
$$
A \psi(2\eta k)\,=\,0\,
$$
for all  $k=-m, -m+1,...,m-1, m$.

Now we shall assume  
that the parameter $\eta $ has the form $\eta = 1/2N$ where $N$ is a 
natural number and the parameter $c$ in the Bethe ansatz equations \Ref{bae2}
obeys the relation $e^{2c}=1$. Under these assumptions the vanishing conditions for
anti-symmetric Bethe eigenfunctions are described in terms of the 
$U_{e^{2\pi i/N}}(sl_2)$ modified fusion rules.

Recall that if $(a_1,...,a_n)$ is a vector of integers,
then we say that $(a_1,...,a_n)$ satisfies the $U_{e^{2\pi i/N}}(sl_2)$
 fusion rules with respect to
$(\La_1,...,\La_n)$ if all $a_j$ are non-negative and for every $j$ the triple
$a_{j-1}, a_j, \La_j$ obeys the $U_{e^{2\pi i/N}}(sl_2)$
 fusion rules (where $a_0=a_n$).
We shall say that $(a_1,...,a_n)$ satisfies {\it the 
$U_{e^{2\pi i/N}}(sl_2)$ modified fusion rules with respect to
$(\La_1,...,\La_n)$} if all $a_j$ are non-negative and for every $j$ the triple
$a_{j}, a_{j+1}, \La_j$ obeys the $U_{e^{2\pi i/N}}(sl_2)$
 fusion rules (where $a_{n+1}=a_1$).
If $\La_1=\La_2=...=\La_n$, then the two fusion rules coincide.

\begin{thm}\label{49}
Let $\eta = 1/2N$ and $e^{2c}=1$. Let $A\psi(\la)$ be
the anti-symmetrization of a Bethe ansatz eigenfunction of the transfer matrix.
Then 
$$
A\psi(2\eta k)=(-1)^m e^c A\psi(2\eta (k+N)).
$$
Moreover,
for any   admissible index $M$ and $k\in \{0, 1, 2,...,N-1\}$ 
we have $A\psi_M(2\eta k)=0$ unless the vector
$$
(-\tilde\Sigma^1 +k-1, -\tilde\Sigma^2 +k-1,...,-\tilde\Sigma^n +k-1)
$$
satisfies the $U_{e^{2\pi i/N}}(sl_2)$ modified 
fusion rules with respect to $(\La_1,...,\La_n)$.
\end{thm}

\section{Restricted Interaction-round-a-face Models}\label{RIRF}

\subsection{Discrete models, \cite{FV2}}
The construction of the transfer matrix admits the following
variation. The transfer matrix $T_{V_{\vec\La}}(w)$
and the difference operators $b(w)$
shift the argument of functions by $\pm 2\eta$. Therefore we may replace
${\mathcal F} (V_{\vec\La})$  by the space ${\mathcal F}_\mu(V_{\vec\La})$ 
of all functions from the set
$C_\mu=\{2\eta(\mu+ j)\,|\,j\in\Z\}$ to $V_{\vec\La}$. The transfer matrix and the
difference operators are then well-defined
on ${\mathcal F}_\mu(V_{\vec\La})$ if $\mu$ is generic. Also, it follows from
Theorem \ref{previous} that the restriction to $C_\mu$ of the
Bethe ansatz eigenfunctions is well defined for all $\mu$.
We thus have:

\begin{corollary}
Suppose $t_1,\dots,t_m$ is a solution to the Bethe ansatz equations
\Ref{bae2}. Then,
for generic $\mu$, the restriction to $C_\mu$ of $b(t_1)\cdots b(t_m)v_c$
is a common eigenfunction of the operators
$T_{V_{\vec\La}}(w)\in\End({\mathcal F}_\mu(V_{\vec\La}[0]))$.
\end{corollary}

\subsection{Unrestricted interaction-round-a-face models, \cite{FV2}}\label{fourths}
In this section, we consider a special case of the above construction,
and relate $T(w)$ to the transfer matrix of the (unrestricted) interaction-round-a-face
(IRF) (also called solid-on-solid) 
models of \cite{Ba, ABF}. Therefore, our formulae give, in particular,
eigenvectors of transfer matrices of IRF models, extending the results
of \cite{BR}.

Let $V=\C^2$. Denote the vectors $e_0, e_1$ of the standard basis 
in $\C^2$ by $e[1], e[-1]$, respectively. Let $V= V[1] \oplus V[-1]= \C e[1]\oplus
\C e[-1]$ be the weight decomposition of $V$. 
Let $R(z, \la)\in \End (\C^2\otimes \C^2)$ be the fundamental R-matrix
defined in Section \ref{ssevm}.

Let $W=V^{\otimes n}$ be the $E_{\tau,\eta}(sl_2)$ representation with 
the L-operator $L(w,\la)\in \End\,(V \otimes W) $ given by
\be
L(w,\lambda)=R^{(01)}(w-z_1, \lambda) \otimes R^{(02)}(w-z_2, \lambda-2\eta h^{(1)}) \otimes
.. \otimes R^{(0,n)}(w-z_{n},\lambda-2\eta \sum_{j=1}^{n-1} h^{(j)}). 
\ee 
This representation is isomorphic to the tensor
product of the evaluation fundamental representations
\be
V(z_n)\otimes ... \otimes V(z_1),
\ee
with the isomorphism
$u_1\otimes\cdots\otimes u_n\mapsto u_n\otimes\cdots\otimes u_1$.

We assume that $n$ is even, and therefore $W$ has a nontrivial zero weight subspace
$W[0]$.

The corresponding transfer matrix $T(w)\in\End({\mathcal F}(W[0]))$
has the form
\bean\label{81}
T(w)\,=\, \sum_\nu \text{Tr}_{V[\nu]^{(0)}}\, L(w,\la)\,{}\,\Gamma_\nu
\eean
where $(\Gamma_\nu f)(\la)=f(\la-2\eta \nu)$.

We define  ``Boltzmann weights''
$w(a,b,c,d;z)$, depending on complex parameters $a,b,c,d,z$, such that
$a-b,b-c,c-d,a-d\in\{1,-1\}$, by 
\bean
R(z, 2\eta d)\, e[{d-c}]\otimes e[{c-b}]
=\sum_a w(a,b,c,d;z) \, e[{a-b}]\otimes e[{d-a}], 
\notag
\eean
(the sum is over one or two allowed values of $a$). 
We set $w(a,b,c,d;z)=0$ if $a, b, c, d$ do not satisfy
$a-b,b-c,c-d,a-d\in\{1,-1\}$.

The dynamical quantum Yang--Baxter
equation translates into {\it the star-triangle equation}
\bean\label{STE} 
&\sum_g
w(a,b,g,f;z_2-z_3) 
w(b,c,d,g;z_1-z_3) 
w(g,d,e,f;z_1-z_2)& 
\notag
\\ 
&= 
\sum_g 
w(b,c,g,a;z_1-z_2)
w(a,g,e,f;z_1-z_3)
w(g,c,d,e;z_2-z_3).& 
\notag
\eean

\begin{figure}
\begin{picture}(350,100)
\put(10,48){\makebox{$w$}}
\put(20,50){\line(1,0){310}}
\put(30,60){\makebox{$b_n$}}
\put(30,35){\makebox{$a_n$}}
\put(47,20){\makebox{$z_n$}}
\put(50,30){\line(0,1){40}}
\put(60,60){\makebox{$b_1$}}
\put(60,35){\makebox{$a_1$}}
\put(77,20){\makebox{$z_{1}$}}
\put(80,30){\line(0,1){40}}
\put(160,60){\makebox{$\cdots$}}
\put(160,35){\makebox{$\cdots$}}
\put(237,20){\makebox{$z_{n-2}$}}
\put(240,30){\line(0,1){40}}
\put(244,60){\makebox{$b_{n-1}$}}
\put(244,35){\makebox{$a_{n-1}$}}
\put(267,20){\makebox{$z_{n-1}$}}
\put(270,30){\line(0,1){40}}
\put(280,60){\makebox{$b_n$}}
\put(280,35){\makebox{$a_n$}}
\put(297,20){\makebox{$z_{n}$}}
\put(300,30){\line(0,1){40}}
\put(310,60){\makebox{$b_1$}}
\put(310,35){\makebox{$a_1$}}
\end{picture}
\caption{Graphical representation of the row-to-row transfer matrix
of an IRF model. Each crossing represents a Boltzmann weight $w$
whose arguments are the labels of the adjoining regions and
the difference of the parameters associated to the lines.}
\label{fig1}
\end{figure}
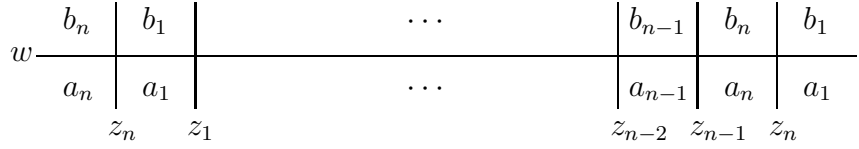

Introduce a basis $|a_1,\dots,a_n\rangle$ of
${\mathcal F}_\mu(W[0])$ labeled by $a_i\in \mu+\Z$ with $a_i-a_{i+1}\in\{1,-1\}$,
$i=1,\dots,n-1$, and $a_n-a_1\in\{1,-1\}$. We let $\delta(\lambda)=
1$ if $\lambda=0$ and $0$ otherwise. Then we define
\be
|a_1,\dots,a_n\rangle(\lambda)
=\delta(\lambda\,-\,2\eta a_1)e[{a_1-a_2}]\otimes e[{a_2-a_3}]\otimes
\cdots\otimes e[{a_n-a_1}].
\ee
If $\Gamma$ is the shift operator $\Gamma f(\lambda)=f(\lambda-2\eta)$,
then $\Gamma|a_1,\dots,a_n\rangle=|a_1\!+\!1,\dots,a_n\!+\!1\rangle$.
Using this, and the fact that $h^{(j)}|a_1,\dots,a_n\rangle=
(a_{j+1}\!-\!a_{j})|a_1,\dots,a_n\rangle$, we get
\bean\label{row-row}
T(w)|a_1,\dots,a_n\rangle
=\sum_{b_1,\dots,b_n}\prod_{j=1}^n
w(b_{j+1},a_{j+1},a_{j},b_{j};w\!-\!z_j)|b_1,\dots,b_n\rangle,
\eean
with the understanding that $b_{n+1}=b_1$, $a_{n+1}=a_1$. The
(finite) sum is over the values of the indices $b_i$ for which
the Boltzmann weights are nonzero.
Comparing with \cite{Ba}, we see that $T(w)$,
in this basis, is the row-to-row transfer matrix of the
(unrestricted) interaction-round-a-face model
associated to the solution $w(a,b,c,d;z)$ of the star-triangle equation,
see \cite{Ba}. The situation is best visualized by 
looking at the graphical representation of Fig.\ \ref{fig1}.

Applying results of Section \ref{t-mat} to the representation $W$, we get the following Corollary:

\begin{corollary}\label{cor51}
Let $c$ be a complex number. Consider the function
$$
\psi(\lambda)=\sum
e^{c\lambda}
\omega_{j_1,...,j_n}(t_1,\dots,t_{n/2},z_1,...,z_n,\lambda)\,
e_{j_1}\otimes ... \otimes e_{j_n}\,{}\,\in \,{}\, {\mathcal F} (W[0])
$$
where the sum is over all $j_1,...,j_n$ such that $j_k=0,1,\, j_1+...+j_n =n/2$
and $\omega_{j_1,...,j_n}$ $(t_1,\dots,t_{n/2},z_1,...,z_n,\lambda)$ are the weight functions
defined in \Ref{w.f} for $a_1=...=a_n=\eta$.
Then for any solution $(t_1,\dots,t_{n/2})$ of the Bethe ansatz
equations 
\begin{equation}\label{bae3}
\prod_{l=1}^{n}
\frac{\theta(t_j-z_l+\eta)}
{\theta(t_j-z_l-\eta)}
\prod_{k:k\neq j}
\frac{\theta(t_j-t_k-2\eta)}{\theta(t_j-t_k+2\eta)}
=e^{-4\eta c},\qquad j=1,\dots,{n/2},
\end{equation}
such that, for all $i<j$, $t_i\neq t_j\mod \Z+\tau\Z$,
the function $\psi(\la)$ is a common eigenfunction of all transfer matrices $T(w)$,
$T(w)\psi(\la)=\epsilon(w)\psi(\la)$, where
\bean\label{value}
\epsilon(w)=
e^{-2\eta c}\prod_{j=1}^{n/2}
\frac
{\theta(t_j-w-\eta)}
{\theta(t_j-w+\eta)}
+
e^{2\eta c}\prod_{j=1}^{n/2}
\frac
{\theta(t_j-w+3\eta)}
{\theta(t_j-w+\eta)}
\prod_{k=1}^n
\frac
{\theta(w-z_k)}
{\theta(w-z_k-2\eta)}.
\eean
Moreover $\psi(\lambda+1)=(-1)^{n/2} e^c\psi(\lambda)$.
\end{corollary}

The Weyl reflection $S\, :\, {\mathcal F} (W[0])\,\to\,
{\mathcal F} (W [0])$ commutes with the transfer matrix.
Hence $S\psi (\la)$ is also an eigenfunction of the transfer matrix with the same
eigenvalue. In particular the Weyl anti-symmetric function $A\psi(\la)=\psi(\la)- S\psi(\la)$
is an eigenfunction of the transfer matrix with the same eigenvalue.

The restriction of the eigenfunctions $\psi(\la)$ and $A\psi(\la)$ to the set $C_\mu$
is well defined for any $\mu$. The $\la$-poles of the L-operator $L(w,z)$ have the form $\la=2\eta k$
mod $\Z + \tau \Z$, where $k$ is an integer. Therefore the transfer matrix 
$T(w)$ restricted to the space
${\mathcal F}_\mu (W [0])$ is well defined for $\mu\neq 0$ mod $\Z + 
{1\over 2\eta} \Z + {\tau\over 2\eta} \Z$. Thus, for those $\mu$,
the restriction of the functions $\psi(\la), A\psi(\la)$ to $C_\mu$ defines common
eigenfunctions of all transfer matrices acting in ${\mathcal F}_\mu (W [0])$
and hence common eigenfunctions of all  row-to-row transfer matrices of the 
corresponding unrestricted interaction-round-a-face models. 
Notice that the eigenvalues of the restrictions do not depend on $\mu$.

\subsection{Infinite restricted interaction-round-a-face models}\label{9.1}
In this section we assume that $\eta$ is a generic complex number.
Let $|a_1,...,a_n\rangle$ be a  delta function such that the numbers $a_1, ..., a_n$ are integers and
for all $j=1,...,n$ we have $|a_j-a_{j+1}|=1$ (with $a_{n+1}=a_1$ ).
Such a delta function is an element of ${\mathcal F}_{\mu=0} (W [0])$.
The delta function will be called {\it positive} (resp. {\it negative}) if $a_1, ..., a_n$
are positive (resp. negative). The delta function will be called {\it neutral} if the set
$a_1, ..., a_n$ contains zero.

Introduce the subspace ${\mathcal F}^+ (W [0]) \subset {\mathcal F}_{\mu=0} (W [0])$ as the subspace
generated by (infinite) linear combinations of all positive delta functions.

The transfer matrix $T(w)$ is not defined on ${\mathcal F}_{\mu=0} (W [0])$ since
some of the Boltzmann weights $w(a,b,c,d;z)$ are not defined when $a,b,c,d $ belong to $ \Z$.
Nevertheless it follows from the formulae for the fundamental R-matrix in Section
\ref{ssevm}  that the Boltzmann weights are well defined if all $a, b, c, d$ are positive integers
 (and $|a-b|=|b-c|=|c-d|=|d-a|=1$ ). 
Moreover, the weights $w(1,2,1,0;z), w(0,1,2,1;z)$ are well defined and $w(0,1,2,1;z)=0$.

 Consider formula \Ref{row-row} for $T(w)|a_1,...,a_n\rangle$. 
If $|a_1,...,a_n\rangle$ is a positive delta function, then in
this formula only $|b_1,...,b_n\rangle$ with non-negative integer coordinates
can appear. Moreover, if $b_j=0$ for some $j$, then $b_{j-1}=b_{j+1}=a_j =1$ and 
$a_{j-1}=a_{j+1}=2$. If there is such $b_j$, then the coefficient of 
$|b_1,...,b_n\rangle$ in $T(w)|a_1,...,a_n\rangle$
contains the product $w(0,1,2,1;w-z_{j-1})w(1,2,1,0;z-z_j)$ which is well defined and equal to zero.
Thus, if $|a_1,...,a_n\rangle$ is a positive delta function, then all terms in the formula for
$T(w)|a_1,...,a_n\rangle$ are well defined and 
\bean
T(w)|a_1,\dots,a_n\rangle
=\sum_{b_1,\dots,b_n}\prod_{j=1}^n
w(b_{j+1},a_{j+1},a_{j},b_{j};w\!-\!z_j)|b_1,\dots,b_n\rangle,
\notag
\eean
where the (finite) sum is over only  the 
positive delta functions $|b_1,\dots,b_n\rangle$ for which the Boltzmann
weights are nonzero. This formula induces a well defined  operator $T^+(w)$ on 
${\mathcal F}^+ (W [0])$ which is 
the row-to-row transfer matrix of the infinite restricted model.

The restricted infinite model, like the unrestricted one, 
satisfies the star-triangular equation, see \cite{ABF}. This equation
holds for all positive integers $a,b,c,d,e,f$ such that $|a-b|=
|b-c|=c-d|=|d-e|=|e-f|=|f-a|=1$ under the condition that the summations are
over all positive integers $g$. This fact easily follows from the unrestricted
star-triangular equation and the explicit formulae
for the matrix coefficients of the fundamental R-matrix. 
Namely, it is enough to prove the restricted star-triangular equation
under the assumption that the set $a,b,c,d,e,f $ contains 1. If the set $a,b,c,d,e,f $ contains 1 and 3,
then the restricted star-triangular equation coincides with the unrestricted one.
Thus the only remaining cases to check are the cases when the set $a,b,c,d,e,f $ is the set
1,2,1,2,1,2 or 2,1,2,1,2,1.  Under
these boundary conditions, it is easy to check that if in the unrestricted equation we
have $g=0$, then 
the corresponding term is well defined and equal to zero. Thus,
the restricted star-triangular equation coincides with the unrestricted one and hence holds.
As a consequence of this result we conclude that
the row-to-row transfer matrices $T^+(w)$ of the infinite restricted model
commute for different values of $w$.

\subsection{Eigenfunctions of the infinite restricted model}
In this section we show that the Weyl anti-symmetrization of a Bethe eigenfunction
of Corollary \ref{cor51} gives a common eigenvector of all row-to-row transfer matrices
of the infinite restricted model.

Let  $T(w)\in\End({\mathcal F}(W[0]))$ be the transfer matrix of the
$E_{\tau,\eta}(sl_2)$ module $W$ defined in \Ref{81}. Let $\psi$ be its Bethe 
eigenfunction corresponding to a complex number $c$ and a
solution $(t_1,...,t_{n/2})$ of the Bethe ansatz
equations \Ref{bae3}. Let $\epsilon (w)$ be the eigenvalue defined in \Ref{value}, 
$T(w)\psi=\epsilon(w)\psi$. Let $A\psi=\psi- S\psi$ 
be the corresponding Weyl anti-symmetric eigenfunction  of $T(w)$ ( with
the same eigenvalue ). Let $A\psi_0 \in {\mathcal F}_{\mu=0} (W [0])$ 
be the restriction of $A\psi$ to the set $C_{\mu=0}$. As an element of
${\mathcal F}_{\mu=0} (W [0])$, it has the form
$$
A\psi_0 \,=\, \sum_{|a_1,...,a_n\rangle} \, (A\psi)_{|a_1,...,a_n\rangle} \, |a_1,...,a_n\rangle,
$$
where $ (A\psi)_{|a_1,...,a_n\rangle}$ are suitable numbers and
the sum is over all $|a_1,...,a_n\rangle$ such that $a_1,...,a_n$
are integers satisfying the condition $|a_1-a_2|=...=|a_{n-1}-a_n|=|a_n-a_1|=1$.
Define an element $A\psi^+ \in {\mathcal F}^+ (W [0])$ by
\bean\label{eigen-plus}
A\psi^+ \,=\, \sum_{|a_1,...,a_n\rangle\,>\,0} \, (A\psi)_{|a_1,...,a_n\rangle} \, |a_1,...,a_n\rangle.
\eean
where the sum is over all positive delta functions.
\begin{thm}\label{restric}
The vector $A\psi^+ \in {\mathcal F}^+ (W [0])$  is a common eigenvector ( with the same
eigenvalue )
of all  row-to-row transfer matrices $T^+(w)$ of the infinite restricted model,
$T^+(w) A\psi^+= \epsilon(w) A\psi^+$.
\end{thm}
\begin{proof}
The Weyl anti-symmetry of $A\psi_0$ means 
\bean\label{ant-sy}
(A\psi)_{|a_1,...,a_n\rangle}\,=\,-\,(A\psi)_{|-a_1,...,-a_n\rangle}
\eean
for all $|a_1,...,a_n\rangle$.
The vanishing conditions of Corollary \ref{vc}  imply
\bean\label{vanish}
(A\psi)_{|a_1,...,a_n\rangle}=0
\notag
\eean 
for all neutral delta functions $|a_1,...,a_n\rangle$.

For generic $\mu$, let 
$$
A\psi_\mu \,=\, \sum_{|\mu+a_1,...,\mu+a_n\rangle} \, (A\psi)_{|\mu+a_1,...,\mu+a_n\rangle} \, 
|\mu+a_1,...,\mu+a_n\rangle
$$
be the restriction of $A\psi$ to $C_\mu$. Here the summation is over integers $a_1,...,a_n$
such that $|a_1-a_2|=...=|a_{n-1}-a_n|=|a_n-a_1|=1$. Let
$$
T(w) A\psi_\mu \,=\, \sum_{|\mu+b_1,...,\mu+b_n\rangle} \, (T(w)A\psi)_{|\mu+b_1,...,\mu+b_n\rangle} \, 
|\mu+b_1,...,\mu+b_n\rangle
$$
be its image under the action of the transfer matrix. Then we have
\bean\label{coeffi}
&(T(w)A\psi)_{|\mu+b_1,...,\mu+b_n\rangle} \,=\,&
\\
&\sum_{a_1,\dots,a_n}\prod_{j=1}^n
w(\mu+b_{j+1},\mu+a_{j+1},\mu+a_{j},\mu+b_{j};w\!-\!z_j)\,
(A\psi)_{|\mu+a_1,...,\mu+a_n\rangle}\,,&
\notag
\eean
as well as 
\bean\label{coe-1}
(T(w)A\psi)_{|\mu+b_1,...,\mu+b_n\rangle} \,=\,\epsilon (w)
(A\psi)_{|\mu+b_1,...,\mu+b_n\rangle}.
\eean
If $|b_1,...,b_n\rangle$ is a positive delta function, then in formula \Ref{coeffi}
only $|\mu+a_1,...,\mu+a_n\rangle$ with nonzero $a_1,...,a_n$ could appear.
If $a_j=0$ for some $j$, then $a_{j-1}=a_{j+1}=b_j =1$ and 
$b_{j-1}=b_{j+1}=2$. If there is such $a_j$, then the coefficient of 
$(A\psi)_{|\mu+a_1,...,\mu+a_n\rangle}$ in $(T(w)A\psi)_{|\mu+b_1,...,\mu+b_n\rangle}$
contains the product $w(\mu+1,\mu+0,\mu+1,\mu+2;w-z_{j-1})
w(\mu+2,\mu+1,\mu+0,\mu+1;z-z_j)$.

For a positive $|b_1,...,b_n\rangle$, consider the limit of 
$(T(w)A\psi)_{|\mu+b_1,...,\mu+b_n\rangle}$ if $\mu$ tends to zero.
According to explicit formulae, the factors 
$w(\mu+b_{j+1},\mu+a_{j+1},\mu+a_{j},\mu+b_{j};w\!-\!z_j)$
with positive $b_{j+1},a_{j+1},a_{j},b_{j}$ are well defined for $\mu=0$
as well as the factors $w(\mu+1,\mu+0,\mu+1,\mu+2;w-z_{j-1})$,
$w(\mu+2,\mu+1,\mu+0,\mu+1;z-z_j)$. At the same time the limit
of the coefficients
$(A\psi)_{|\mu+a_1,...,\mu+a_n\rangle}$ is zero if $a_1,...,a_n$ are non-negative
and at least one of them is zero. Thus, taking the limit of formulae \Ref{coeffi} and
\Ref{coe-1}, we get $T^+(w) A\psi^+= \epsilon(w) A\psi^+$. 

\end{proof}

\subsection{Finite interaction-round-a-face models}\label{2nd-mot}
In this section we assume that $\eta=1/2N$ where $N$ is a natural number.
A delta function $|a_1,...,a_n\rangle$ with integers $a_1,...,a_n$
obeying conditions  $|a_1-a_2|=...=|a_{n-1}-a_n|=|a_n-a_1|=1$ and 
$ 0 < a_k < N $ for all $k$ will be called {\it admissible }.
A delta function $|a_1,...,a_n\rangle$ with at least one
of the integers $a_1,...,a_n$ divisible by $N$ will be called {\it bad}.
Introduce the finite dimensional subspace ${\mathcal F}^{\{N\}} (W [0])$
of ${\mathcal F}_{\mu=0} (W [0])$ as the subspace generated by the
admissible delta functions.

The linear operator $T^{\{N\}}(w)\in\End({\mathcal F}^{\{N\}}(W[0]))$
defined by
\bean
T^{\{N\}}(w)|a_1,\dots,a_n\rangle
=\sum_{b_1,\dots,b_n}\prod_{j=1}^n
w(b_{j+1},a_{j+1},a_{j},b_{j};w\!-\!z_j)|b_1,\dots,b_n\rangle,
\notag
\eean
where the sum is over admissible delta functions, is called 
the row-to-row transfer matrix of the (finite) restricted interaction-round-a-face model.
Considerations, similar to those in Section \ref{9.1}, show that the transfer
matrix is well defined, and the transfer matrices commute for different values
of $w$, see \cite{ABF}.

Let  $T(w)\in\End({\mathcal F}(W[0]))$ be the transfer matrix of the
$E_{\tau,\eta}(sl_2)$ module $W$ defined in \Ref{81}. Let $\psi$ be its Bethe 
eigenfunction corresponding to a complex number $c$ and a
solution $(t_1,...,t_{n/2})$ of the Bethe ansatz
equations \Ref{bae3}. In this section we always assume that $e^{2c}=1$.
 Let $\epsilon (w)$ be the eigenvalue defined in \Ref{value}, 
$T(w)\psi=\epsilon(w)\psi$. Let $A\psi=\psi- S\psi$ 
be the corresponding Weyl anti-symmetric eigenfunction  of $T(w)$. 
Let $A\psi_0 \in {\mathcal F}_{\mu=0} (W [0])$ 
be the restriction of $A\psi$ to the set $C_{\mu=0}$, $A\psi_0 \,=\, 
\sum_{|a_1,...,a_n\rangle} \, (A\psi)_{|a_1,...,a_n\rangle} \, |a_1,...,a_n\rangle$.
Define an element $A\psi^{\{N\}} \in {\mathcal F}^{\{N\}} (W [0])$ by
\bean\label{eigen-plus}
A\psi^{\{N\}} \,=
\, \sum_{\text{adm}\,|a_1,...,a_n\rangle} \, (A\psi)_{|a_1,...,a_n\rangle} \, |a_1,...,a_n\rangle.
\eean
where the sum is over all admissible delta functions.
\begin{thm}\label{finite-restric}
The vector $A\psi^{\{N\}} \in {\mathcal F}^{\{N\}} (W [0])$  is a common eigenvector ( with the same
eigenvalue ) of all  row-to-row transfer matrices $T^{\{N\}}(w)$ of the (finite) restricted model,
$T^{\{N\}}(w) A\psi^{\{N\}}= \epsilon(w) A\psi^{\{N\}}$.
\end{thm}
\begin{proof}
The Weyl anti-symmetry of $A\psi_0$ gives the relation \Ref{ant-sy}.
The condition $e^{2c}=1$ implies
\bean\label{}
(A\psi)_{|a_1+N,...,a_n+N\rangle}\,=\,(-1)^{n/2}e^c\,(A\psi)_{|a_1,...,a_n\rangle}
\notag
\eean
for all $|a_1,...,a_n\rangle$.
The vanishing conditions of  Theorem \ref{49} imply
\bean\label{vanish}
(A\psi)_{|a_1,...,a_n\rangle}=0
\notag
\eean 
for all bad  $|a_1,...,a_n\rangle$.

Now the Theorem is proved similarly to the proof of Theorem \ref{restric}
considering the function $A\psi_\mu$ and taking the limit 
of the equation $T(w)A\psi_\mu=\epsilon(w) A\psi_\mu$ as $\mu$ tends to zero.
\end{proof}

{\bf Remark.} The eigenvector $A\psi^{\{N\}}$ has an additional symmetry
$$
(A\psi)_{|a_1,...,a_n\rangle}\,=\,(-1)^{n/2+1}e^c\,(A\psi)_{|N-a_1,...,N-a_n\rangle}.
$$

\end{document}